\documentstyle[amscd, amstex]{amsart}

\newtheorem{thm}{Theorem}[section]
\newtheorem{lem}[thm]{Lemma}

\newtheorem{pr}[thm]{Proposition}

\theoremstyle{definition}
\newtheorem{rem}[thm]{Remark}

\newtheorem{defn}[thm]{Definition}

\newcommand{\psx}{{\bf P}_{\Sigma_X}}

\newcommand{\ts}{{\bf T}_\sigma}
\newcommand{\mult}{{\rm mult}}
\newcommand{\im}{{\rm im}}
\newcommand{\key}{\bibitem}
\newcommand{\psd}{{{\bf P}_{\Sigma_D}}}
\newcommand{\ps}{{{\bf P}_{\Sigma}}}

\newcommand{\gr}{{\rm Gr}}
\newcommand{\res}{{\rm Res}}
\newcommand\hidot{{\raise1pt\hbox{$\scriptscriptstyle\bullet$}}}
\newcommand\lodot{{\raise.3pt\hbox{$\scriptscriptstyle\bullet$}}}
\newcommand{\dd}{{\rm d}}
\newcommand{\pp}{{\bf P}}
\newcommand{\ttt}{{\bf T}}

\begin{document}

\title 
{On  the chiral ring of Calabi-Yau hypersurfaces in toric varieties}
\author{Anvar R. Mavlyutov}
\address {Department of Mathematics \& Statistics,
University of Massachusetts, Amherst, MA, 01003, USA}
\email{anvar@@math.umass.edu}

\maketitle

\tableofcontents

\setcounter{section}{-1}

\section{Introduction}

In nonlinear sigma models (see Appendix B.2 in \cite{ck} on physical theories) there are two twisted 
theories, the A-model
and the B-model. Mirror symmetry is an isomorphism between the A-model and the B-model for a pair of two 
distinct Calabi-Yau threefolds $V$ and $V^\circ$ with K\"{a}hler structures. 
One consequence of  mirror symmetry is an isomorphism between the quantum cohomology on 
$\oplus_{p,q}H^{p,q}(V)$ and
the chiral ring of the B-model $\oplus_{p,q} H^p(V^\circ,\wedge^q{\cal T}_{V^\circ})$, 
which implies the equality of the corresponding correlation functions (Yukawa couplings).
 These correlation functions describe  interactions between strings.
{From} a mathematical point of view,  knowledge about the B-model Yukawa coupling and the equality with the 
A-model Yukawa coupling of the mirror manifold produces enumerative information on this mirror manifold.
One important construction widely used in physics and mathematics is the Batyrev mirror construction in 
toric varieties (see \cite{b1}).

In this paper we study the chiral ring $\oplus_pH^p(X,\wedge^p{\cal T}_X)$ (this is actually a subring of the whole 
chiral ring) for quasismooth hypersurfaces $X$
in complete simplicial toric varieties.
In particular, we completely describe the chiral ring $\oplus_pH^p(X,\wedge^p{\cal T}_X)$ 
in the case of 3-dimensional Calabi-Yau  hypersurfaces. 
This applies to the mirror symmetric hypersurfaces in Batyrev's construction.

The following is an outline of the paper. We begin in Section~\ref{s:sem} with a review of  notation and general
facts {from} toric geometry. For complete toric varieties, the notions of semiample, nef (numerically effective) and
generated by global sections are equivalent for invertible sheaves (divisors). Geometry and intersection theory associated with
big (the self-intersection number is positive) and nef divisors was studied in \cite{m}.  
Here,  we generalize those results to all semiample divisors. Such divisors on complete toric varieties
 naturally produce a surjective morphism of the ambient space onto another complete toric variety. 
Moreover, this construction is unique with certain conditions relating semiample divisors to ample divisors on the new
toric variety. We also show that a proper  birational   morphism of toric varieties
induces a natural graded homomorphism of the coordinate rings. For a semiample divisor,
this gives an isomorphism of rings in the degree of the divisor.

 Section~\ref{s:pr} uses the results of Section~\ref{s:sem} to  describe the geometry of
semiample regular (transversal to orbits) hypersurfaces in complete toric varieties.
 We get a stratification of such hypersurfaces in terms of nondegenerate affine hypersurfaces
cohomology of which has been studied in \cite{b2}.  
We also review some facts about hypersurfaces in complete simplicial toric varieties. 
In particular, we recall {from} \cite{m} the relationship
between the Jacobian ring $R(f)$ (resp.,  $R_1(f)$) and the middle cohomology $H^{d-1}(X)$ of a quasismooth 
(resp., big and nef regular) hypersurface $X$ in a $d$-dimensional complete simplicial toric variety.
The ring $R_1(f)$ has been used in \cite{m} to describe the middle cohomology of a 3-dimensional big and nef
regular hypersurface completely. 

In Section~\ref{s:p}, we introduce the {\it (Zariski) $p$-th  exterior power $\wedge^p{\cal T}_X$ 
of the tangent sheaf} for an arbitrary orbifold $X$, which is defined 
similarly to the sheaf $\Omega_X^p$ of Zariski $p$-forms (see \cite[A.3]{ck}). Then we show that for a 
quasismooth hypersurface $X$ of degree $\beta$ there is a ring homomorphism 
$R(f)_{*\beta}\rightarrow H^*(X,\wedge^*{\cal T}_X)$ (the latter is our notation for $\oplus_pH^p(X,\wedge^p{\cal T}_X)$). 
Also, with respect to this homomorphism the map
between $R(f)$ and the middle cohomology of a quasismooth hypersurface is a morphism of modules.
In the Calabi-Yau case  the situation is especially nice because we get an injective 
ring homomorphism $R_1(f)_{*\beta}\rightarrow H^*(X,\wedge^*{\cal T}_X)$ 
(we call $R_1(f)_{*\beta}$ the polynomial
part of the chiral ring because its graded piece in $H^1(X,{\cal T}_X)$ should correspond to polynomial
infinitesimal deformations for a minimal Calabi-Yau $X$ (see \cite{ck})). 
 
According to the above terminology, in Section~\ref{s:non}, we study the non-polynomial part of the chiral
ring complementary to the polynomial part. We construct new elements 
 in $H^*(X,\wedge^*{\cal T}_X)$ for 
a big and nef quasismooth hypersurface $X$, and
in the case of a minimal Calabi-Yau these elements in $H^1(X,{\cal T}_X)$ should correspond to 
non-polynomial deformations. The new elements  are represented by a map {from} some quotient $R^\sigma(f)$
of the Jacobian ring  to $H^*(X,\wedge^*{\cal T}_X)$, and this map is actually a  morphism of modules with respect to 
the ring homomorphism $R(f)_{*\beta}\rightarrow H^*(X,\wedge^*{\cal T}_X)$. 
We also calculate some vanishing cup products of the new elements.
 The new part of $H^*(X,\wedge^*{\cal T}_X)$ has its analogue in 
the middle cohomology $H^{d-1}(X)$ of the hypersurface. This is also given by a map {from} certain graded pieces of $R^\sigma(f)$
to  $H^{d-1}(X)$. We show that this map is morphism of modules with respect to 
$R(f)_{*\beta}\rightarrow H^*(X,\wedge^*{\cal T}_X)$. 

In Section~\ref{s:tr}, we describe the {\it toric part} of cohomology of a semiample regular hypersurface. This part 
is the image of cohomology of the ambient space, while its complement, called the {\it residue part},  comes 
{from} the residues of rational differential forms with poles along the hypersurface. We show that the cohomology of 
a semiample regular hypersurface is a direct sum of its toric and residue parts. 

Section~\ref{s:cohreg} studies  the middle cohomology of a big and nef regular hypersurface. 
We  provide a better and more general description of the middle cohomology than the one given
for 3-dimensional hypersurfaces in \cite{m}.
Here, we use a new ring $R_1^\sigma(f)$, analogous to the ring $R_1(f)$.
An algebraic description of the middle cohomology can be used in the Calabi-Yau case to compute 
the product structure on the chiral ring.

In Section~\ref{s:cal3}, we consider semiample anticanonical regular hypersurfaces.
Such hypersurfaces are Calabi-Yau, 
implying that their chiral ring  is isomorphic to the middle cohomology. 
Using the description of Section~\ref{s:cohreg}, we have a partial description of the space $H^*(X,\wedge^*{\cal T}_X)$
in terms of  $R_1(f)$ and $R_1^\sigma(f)$. We show that this part is a subring of the chiral ring.
This subring is the whole $H^*(X,\wedge^*{\cal T}_X)$ in the case of Calabi-Yau threefolds.
 The product structure of the polynomial part $R_1(f)$  is in Section~\ref{s:p}, 
while the product of two different elements {from} $R_1(f)$ and $R_1^\sigma(f)$ is in Section~\ref{s:non}.
We describe the nontrivial product structure on the spaces $R_1^\sigma(f)$ in terms of triple products.
Since $H^*(X,\wedge^*{\cal T}_X)$ and the described subring 
have a nondegenerate pairing, induced by the cup product on the middle cohomology,
one can recover the chiral ring structure completely on these spaces.

{\it Acknowledgments.} I am very grateful to  David Cox for his support and valuable
comments. 
Parts of this work were inspired by some notes of David Cox and David Morrison  whose results we
include in Section~\ref{s:p}.

\section{Semiampleness}\label{s:sem}

In this section we first review some basic facts and notation,  and then generalize the  geometric construction of \cite{m}
associated with semiample divisors on complete toric varieties.
We show that a semiample divisor naturally produces a surjective morphism of the ambient space onto 
another complete toric variety. This construction is unique with certain conditions which relate
the semiample divisor to an ample divisor on the new toric variety. 
At the end of this section we show that a  proper birational morphism of toric varieties
gives a natural graded homomorphism of the homogeneous coordinate rings of the  varieties.
We apply this to the maps associated with semiample divisors.
 
Let $M$ be a lattice of rank $d$, then $N=\text{Hom}(M,{\Bbb Z})$ is the 
dual lattice; $M_{\Bbb R}$ (resp.~$N_{\Bbb R})$ denotes
the $\Bbb R$-scalar extension of $M$ (resp.~of~$N$).
The symbol ${\bf P}_{\Sigma}$ stands for a  $d$-dimensional 
toric variety associated with   
a finite rational   fan $\Sigma$ 
in $N_{\Bbb R}$. 
A toric variety $\ps$ is a disjoint union of its orbits by the action of the 
torus 
${\bf T}=N\otimes{\Bbb C}^*$ that sits naturally inside 
$\ps$. Each orbit ${\bf T}_\sigma$ is a torus corresponding to a cone 
$\sigma\in\Sigma$. The closure of each orbit ${\bf T}_\sigma$ is
again a toric variety denoted $V(\sigma)$.    

We use $\Sigma(k)$  for the set of all $k$-dimensional cones in $\Sigma$; in particular, 
$\Sigma(1)=\{\rho_1,\dots,\rho_n\}$ is the set of $1$-dimensional
cones in $\Sigma$ with the minimal integral generators 
$e_1,\dots,e_n$, respectively.
Each 1-dimensional cone $\rho_i$ corresponds to
a torus invariant divisor $D_i$ in ${\bf P}_\Sigma$. 

A torus invariant Weil divisor $D=\sum_{i=1}^{n}a_iD_i$ determines
a convex  polyhedron
$$\Delta_D=\{m\in M_{\Bbb R}:\langle m,e_i\rangle\geq-a_i
\text{ for all } i\}\subset M_{\Bbb R}.$$
Each Weil divisor $D$ gives a reflexive sheaf $O_{\ps}(D)$, whose sections over 
$U\subset \ps$ are the rational functions $f$ such that ${\rm div}(f)+D\ge0$ on $U$.
When $D=\sum_{i=1}^{n}a_iD_i$ is Cartier, 
there is a support function 
$\psi_D : N_{\Bbb R}\rightarrow\Bbb R$ that is linear on each
cone $\sigma\in\Sigma$ and determined by some $m_\sigma\in M$:
$$\psi_D(e_i)=\langle m_\sigma, e_i\rangle=-a_i\text{ for all }
e_i\in\sigma.$$  

When  ${\bf P}_\Sigma$ is complete, the polyhedron $\Delta_D$ of a torus invariant Weil divisor is bounded and called polytope.
Also, the line bundle $O_{\ps}(D)$, corresponding to a Cartier divisor $D$, is generated 
by global sections   if and only if $\psi_D$ is  convex.

We call a Cartier divisor $D$ on a complete toric variety $\ps$  {\it semiample} if
$O_{\ps}(D)$ is generated by global sections.

\begin{rem} This definition is consistent with the one in \cite[\S~5]{ev}  used in a non-toric context for
projective varieties, because
an invertible sheaf ${\cal L}$ on a complete toric variety is generated by global sections iff some positive 
power ${\cal L}^k$  is generated by global sections.
\end{rem}

Theorem~1.6  in \cite{m}
shows: $O_{\ps}(D)$ is generated by global sections is equivalent to the condition that the divisor $D$ is
nef (numerically effective). Therefore, the notions of semiample and nef are equivalent for divisors on complete
toric varieties.

Following \cite[\S~5]{ev}, a semiample divisor $D$  on $\ps$ also has the Iitaka dimension:
$$\kappa(D):=\kappa(O_{\ps}(D))=\dim \phi_D(\ps),$$
where $\phi_D:\ps@>>>{\Bbb P}(H^0(\ps,O_{\ps}(D)))$ is the map defined by the sections of the line bundle $O_{\ps}(D)$. 
The possible values for this characteristic are  $\kappa(D)=0,\dots,\dim\ps$. Moreover, the
Exercise on page~73 in \cite[Section~3.4]{f1} shows that $\kappa(D)$ for a torus invariant $D$
 is exactly the dimension of the associated polytope $\Delta_D$.
It will be convenient for us to introduce the following notion.

\begin{defn} A semiample divisor $D$ on a complete toric variety $\ps$
is called {\it $i$-semiample} if the  Iitaka dimension $\kappa(D)=i$.
\end{defn}

\begin{rem}\label{r:inm} In \cite{m} we called a  Cartier divisor 
$D$  semiample if $O_{\ps}(D)$ is
generated by global sections and the intersection number $(D^d)>0$. 
In fact, such divisors have  the maximal Iitaka dimension $\kappa(D)=\dim\ps$.
In the common terminology, they correspond to big ($(D^d)>0$) and nef, and, according to the above definition,
we should call them  $d$-semiample with $d=\dim\ps$. 
\end{rem}

All ample divisors on $\ps$  are semiample and have the  Iitaka dimension equal to $\dim\ps$.
Our goal is to show that semiample divisors give rise to a natural geometric construction connected with
ample divisors.
Let $D=\sum_{k=1}^{n}a_kD_k$ be an $i$-semiample divisor on $\ps$ with the convex support function $\psi_D$.
For each $d$-dimensional cone $\sigma\in\Sigma$ there is a unique $m_\sigma\in M$ such that
$\psi_D(v)=\langle m_\sigma, v\rangle$ for all $v\in\sigma$. 
Glue together the maximal dimensional cones in $\Sigma$ with the same value
$m_\sigma$. The glued set $\tau(m_\sigma)$ is a convex rational polyhedral cone. Indeed, 
let $v$ be in the convex hull of  $\tau(m_\sigma)$, then
 $\psi_D(v)\le\langle m_\sigma, v\rangle$, by the convexity of the support function. 
On the other hand,
$v$ is lying in some $d$-dimensional cone, where the value of $\psi_D$ is determined
by $m'\in M$. Hence, $\psi_D(e_k)=\langle m_\sigma, e_k\rangle\le\langle m', e_k\rangle$
for all generators $e_k$ {from} the set $\tau(m_\sigma)$. Since $v$ is a positive  linear
 combination of some generators lying in $\tau(m_\sigma)$, we get 
$\langle m_\sigma, v\rangle\le\langle m', v\rangle=\psi_D(v)$. Therefore, the glued set $\tau(m_\sigma)$
coincides with its convex hull. The new cones $\tau(m_\sigma)$ are not necessarily strongly convex,
but they all contain the same linear subspace 
\begin{equation}\label{e:sub}
\tau(m_\sigma)\cap(-\tau(m_\sigma))=\{v\in N_{\Bbb R}:\psi_D(-v)=-\psi_D(v)\}.
\end{equation} 
To see the equality note that $\langle m_\sigma, w\rangle\ge\psi_D(w)$ for any $w$, by the convexity of the support function. 
Therefore, for $v$ in the right-hand side of (\ref{e:sub}),
we have $\langle m_\sigma, -v\rangle\ge\psi_D(-v)=-\psi_D(v)\ge-\langle m_\sigma, v\rangle$ implying that 
$v\in\tau(m_\sigma)\cap(-\tau(m_\sigma))$. The other way is obvious. {From} here we get that
the linear space in (\ref{e:sub}) consists of $v\in N_{\Bbb R}$ such that $\langle m_\sigma, v\rangle$ is the same for all 
$\sigma$. Since $O_{\ps}(D)$ is
generated by global sections, the polytope $\Delta_D$ is the convex hull of $m_\sigma$.
Therefore, the dimension of (\ref{e:sub}) is exactly $d-i$.
If $\Delta_D$ contains the origin, this linear space can be obtained as the orthogonal complement of the polytope.

Denote by $N'=\{v\in N:\psi_D(-v)=-\psi_D(v)\}$ a sublattice of $N$, we also get the quotient lattice $N_D:=N/N'$.
Then the $i$-dimensional linear
space $N'_{\Bbb R}$ is a support of a complete fan  $\Sigma'$ filled up by the cones of the fan $\Sigma$ contained in $N'_{\Bbb R}$.
The quotient sets $\tau(m_\sigma)/N'_{\Bbb R}$ in ${(N_D)}_{\Bbb R}$ are strongly convex polyhedral cones and form another
complete fan $\Sigma_D$. Thus, we get the following picture: there is a natural exact sequence of lattices 
$$0@>>>N'@>>>N@>>>N_D@>>>0$$
compatible with the fans $\Sigma'$, $\Sigma$ and $\Sigma_D$, giving rise to toric morphisms
$$\pp_{\Sigma'} @>\nu>> \ps @>\pi>>\psd.$$
Let us note that linearly equivalent semiample divisors $D$ produce the same construction. 
The complete toric variety $\pp_{\Sigma'}$ is mapped into an open toric subvariety 
$\pp_{\tilde{\Sigma}'}\subset\ps$ given
by the subfan ${\tilde{\Sigma}'}\subset\Sigma$ of all cones that lie in $N'_{\Bbb R}$. 
Section~2.1 in \cite{f1} shows that the above sequence of toric morphisms  
induces a trivial fibration over the maximal dimensional torus $\ttt_{\Sigma_D}:=N_D\otimes{\Bbb C}^*$ of $\psd$:
$$\pp_{\Sigma'} @>\nu>>\pp_{\tilde{\Sigma}'} @>\pi>>\ttt_{\Sigma_D}.$$

We next show that the above construction is unique in a certain sense.
Using a standard description of a toric morphism, we can see that the toric subvarieties $V(\gamma)\subset\ps$ of dimension $i$,
such that 
$\gamma\in\Sigma(d-i)$ and $\gamma\subset N'_{\Bbb R}$, map birationally onto $\psd$. As in \cite[Section]{f1},
let us restrict the semiample divisor $D=\sum_{k=1}^{n}a_kD_k$ to  $V(\gamma)$. Using the linear equivalence,
we can assume that the origin is one of the vertices of the polytope $\Delta_D$. In this case, equation  (\ref{e:sub}) implies that
$a_k=0$ for $\rho_k\subset N'_{\Bbb R}$ (equivalently, $\psi_D=0$ on $N'_{\Bbb R}$),
 whence $V(\gamma)$ is not contained in the support of $D$. 
Therefore, we get a Weil divisor $D\cdot V(\gamma)$ in the Chow group $A_{i-1}(V(\gamma))$ 
representing the Cartier divisor $D|_{V(\gamma)}$.
Its support function $\psi_{D\cdot V(\gamma)}$ is represented by $\psi_D$ which descends to the quotient space 
${(N_D)}_{\Bbb R}=N_{\Bbb R}/N'_{\Bbb R}$. 
The lattice $M_D:={N'}^\perp\cap M$ is the dual of $N_D$, and the polytope $\Delta_D$ contained in ${(M_D)}_{\Bbb R}$ is exactly
the polytope of the Weil divisor $D\cdot V(\gamma)$. 
By construction, the function $\psi_{D\cdot V(\gamma)}$ is strictly convex with respect to the fan $\Sigma_D$.
Now the arguments of \cite[Section~1]{m} show that $\Sigma_D$ is the normal fan of $\Delta_D$,
and the pushforward $\pi_*(D\cdot V(\gamma))$ is an ample divisor. We also get a commutative diagram (see \cite{f2}):
\[ \begin{CD}
A_{i-1}(V(\gamma))@>\pi_*>>A_{i-1}({\bf P}_{\Sigma_D})\\
@AAA  @AAA \\
\text{Pic}({\bf P}_\Sigma) @<<\pi^*< \text{Pic}({\bf P}_{\Sigma_D}),
\end{CD} \]
where the right vertical arrow is injective and the left is the composition 
$\text{Pic}({\bf P}_\Sigma)@>>>\text{Pic}(V(\gamma))@>>>A_{i-1}(V(\gamma))$ of the restriction map
and the inclusion.
Since the support function of
$\pi_*(D\cdot V(\gamma))$ is induced by $\psi_D$, we have the equality
$\pi^*\pi_*[D\cdot V(\gamma)]=[D]$ in the Chow group $A_{d-1}({\bf P}_{\Sigma})$.  

Now we prove that the conditions on the divisor $D$ deduced in the previous paragraph uniquely determine the constructed morphism.
 Let $p:\pp_{\Sigma}@>>>\pp_{\Sigma_1}$ be a surjective morphism of complete toric varieties 
arising {from} a surjective homomorphism of lattices $\tilde p:N@>>>N_1$ which maps the fan $\Sigma$ into $\Sigma_1$.
The kernel of $\tilde p$ is a sublattice $N_2\subset N$. It is not difficult to see that a cone of $\Sigma$
is either lying in the space ${(N_2)}_{\Bbb R}$ or its relative interior has no intersection with this space.
Hence, the space ${(N_2)}_{\Bbb R}$ is a support of a complete fan $\Sigma_2$ filled up by those cones of $\Sigma$ lying in
${(N_2)}_{\Bbb R}$.
The toric subvarieties $V(\gamma)$ corresponding to $\gamma\in\Sigma(d-k)$ ($k:=\dim\pp_{\Sigma_1}$), contained in ${(N_2)}_{\Bbb R}$,
 are the only ones mapping birationally
onto $\pp_{\Sigma_1}$.
 Suppose now that we have
an $i$-semiample (torus invariant) divisor $D$ on $\pp_{\Sigma}$ such that
$p_*[D\cdot V(\gamma)]$ is ample and  
$p^*p_*[D\cdot V(\gamma)]=[D]$ for some $V(\gamma)$, $\gamma\in\Sigma(d-k)$, which maps birationally onto $\pp_{\Sigma_1}$.
Then the polytope of the divisor $p_*(D\cdot V(\gamma))$ has dimension equal to $\dim\pp_{\Sigma_1}$.
On the other hand, the support function of $D$ is induced by the support function of $p_*(D\cdot V(\gamma))$, 
implying that the polytopes of these divisors is the same set in $M\cap N_2^\perp$. Therefore,
the dimension of $\pp_{\Sigma_1}$ is $i$, and  the fan ${\Sigma_1}$ coincides with $\Sigma_D$ constructed before.
Thus, we proved the following.

\begin{thm}\label{t:fun}
Let $[D]\in A_{d-1}(\ps)$ be an $i$-semiample divisor class on  a complete toric variety $\bold P_\Sigma$ of dimension $d$.
Then, there exists a unique complete toric variety $\psd$ with a surjective morphism
$\pi:{\bf P}_\Sigma@>>>{\bf P}_{\Sigma_D}$, corresponding to a map of $\Sigma$ into $\Sigma_D$, 
such that $\pi_*[D\cdot V(\gamma)]$ is ample and $\pi^*\pi_*[D\cdot V(\gamma)]=[D]$ for some closed toric subvariety
$V(\gamma)\subset\ps$, $\gamma\in\Sigma$, which maps birationally onto $\psd$. Moreover,  $\dim\psd=i$, and the fan 
$\Sigma_D$ is the normal fan of $\Delta_D$ for a torus invariant $D$.
\end{thm}

\begin{rem}\label{r:can}
The fan $\Sigma_D$ is canonical with respect to the equivalence relation on the divisors.
Therefore, it will sometimes be convenient for us to use the notation $\Sigma_\beta:=\Sigma_D$ for a semiample divisor
class $\beta=[D]\in A_{d-1}(\ps)$. 
\end{rem}

While a restriction of a semiample divisor $D$ on $\ps$ to a closed toric subvariety is again a semiample divisor,
the Iitaka dimension of the restricted divisor may change. Let us investigate this problem. 
If $D$ is an $i$-semiample divisor  on $\ps$, then, by the above theorem, we have a unique toric morphism 
$\pi:\ps@>>>\psd$, arising {from} a homomorphism 
$\tilde\pi:N_{\Bbb R}@>>>{(N_D)}_{\Bbb R}$ mapping $\Sigma$ into $\Sigma_D$. 
This morphism encodes information about the structure of the variety $\ps$.
The Iitaka dimension of the semiample divisor  $D\cdot V(\sigma)$ on $V(\sigma)$, 
$\sigma\in\Sigma$, can be determined in the following way.
The complete toric variety $V(\sigma)$ is mapped onto a closed subvariety $V(\sigma_0)\subset\psd$ such that 
the cone $\sigma_0\in\Sigma_D$ is the smallest that contains $\tilde\pi(\sigma)$. We claim that this induced map
$\pi:V(\sigma)@>>>V(\sigma_0)$ is exactly the one associated with the semiample divisor $D\cdot V(\sigma)$.
To prove this we will verify  the conditions which uniquely determine  such a morphism.
As in the theorem above, let $V(\gamma)$ be  such that $\pi_*[D\cdot V(\gamma)]$ is ample and $\pi^*\pi_*[D\cdot V(\gamma)]=[D]$,
and let 
$V(\gamma')\subset V(\sigma)$  be a closed toric subvariety mapping birationally onto $V(\sigma_0)$.
By the projection formula (see \cite{f2}), we get
\begin{multline*}\pi_*[D\cdot V(\gamma')]=\pi_*[(\pi^*\pi_*[D\cdot V(\gamma)])\cdot V(\gamma')]
\\
=\pi_*[D\cdot V(\gamma)]\cdot \pi_*[V(\gamma')]=
\pi_*[D\cdot V(\gamma)]\cdot V(\sigma_0)
\end{multline*}
in the Chow group of the toric variety $V(\sigma_0)$. Since $\pi_*[D\cdot V(\gamma)]$ is ample, the divisor class 
$\pi_*[D\cdot V(\gamma')]$ is ample as well. The other condition for the semiample divisor $D\cdot V(\sigma)$ also follows:
$$\pi^*\pi_*[D\cdot V(\gamma')]=\pi^*[\pi_*[D\cdot V(\gamma)]\cdot V(\sigma_0)]=
\pi^*\pi_*[D\cdot V(\gamma)]\cdot V(\sigma)=[D\cdot V(\sigma)],$$
where we used the commutative diagram
\[ \begin{CD}
\text{Pic}({\bf P}_{\Sigma_D})@>\pi^*>> \text{Pic}({\bf P}_\Sigma)\\
@VVV  @VVV \\
\text{Pic}(V(\sigma_0))@>\pi^*>> \text{Pic}(V(\sigma)).
\end{CD} \]
Thus, by the uniqueness part of  Theorem~\ref{t:fun},  we get the next result.
 
\begin{pr}\label{p:comp1} Let $[D]\in A_{d-1}(\ps)$ be an $i$-semiample divisor class on  $\bold P_\Sigma$ with the associated
morphism $\pi:\ps@>>>\psd$ arising {from} a map of the fan $\Sigma$ into $\Sigma_D$. 
Then, for $\sigma\in\Sigma$, the restriction 
 $[D\cdot V(\sigma)]$ is a $k$-semiample divisor class on $V(\sigma)$ with
$k=i-\dim(\sigma_0)=\dim V(\sigma_0)$, where $\sigma_0\in\Sigma_D$ is the smallest cone that contains the image of $\sigma$. 
Moreover, the induced map $\pi:V(\sigma)@>>>V(\sigma_0)$ is the one associated with the semiample divisor class
$[D\cdot V(\sigma)]$.
\end{pr}

This proposition says that the maps associated with the semiample divisors are compatible with the restrictions.

Any toric variety $\ps$ has a homogeneous coordinate ring
$S(\Sigma)={\Bbb C}[x_1,\dots,x_n]$ with variables $x_1,\dots,x_n$
corresponding to the irreducible torus invariant divisors   $D_1,\dots,D_n$.
This ring is graded by the Chow group $A_{d-1}(\ps)$, assigning $[\sum_{i=1}^n a_i D_i]$ 
to $\deg(\prod_{i=1}^n x_i^{a_i})$. 
For a Weil divisor $D$ on $\ps$, there is an isomorphism 
$H^0(\ps, O_\ps(D))\cong S(\Sigma)_\alpha$, where 
$\alpha=[D]\in A_{d-1}(\ps)$. If $D$ is torus invariant, the monomials in $S(\Sigma)_\alpha$ correspond to
the lattice points of the associated polyhedron $\Delta_D$.
 
Now consider a proper birational morphism $\pi:\pp_{\Sigma_1}@>>>\pp_{\Sigma_2}$ of toric varieties, associated with 
a subdivision
$\Sigma_1$ of $\Sigma_2$. 
In this situation, the 1-dimensional cones of the two fans are related by
$\Sigma_2(1)\subset\Sigma_1(1)$,
and there is a natural relation of the coordinate rings $S(\Sigma_1)={\Bbb C}[x_k:\rho_k\in\Sigma_1(1)]$ 
and $S(\Sigma_2)={\Bbb C}[y_k:\rho_k\in\Sigma_2(1)]$ of the toric varieties. 
For $\alpha=[O_{\pp_{\Sigma_1}}(D)]\in A_{d-1}(\pp_{\Sigma_1})$ we have a commutative diagram:
$$\minCDarrowwidth0.7cm
\begin{CD}
S(\Sigma_1)_{\alpha}  @.\cong @. H^0(\pp_{\Sigma_1},O_{\pp_{\Sigma_1}}(D))\\
@VVV       @.                    @VVV \\
S(\Sigma_2)_{\pi_*\alpha}@.\cong @. H^0(\pp_{\Sigma_2},O_{\pp_{\Sigma_2}}(\pi_*D)),
\end{CD}$$
where the left vertical arrow sends a monomial $\prod_{\rho_k\in\Sigma_1(1)}x_k^{a_k+\langle m_,e_k\rangle}$ in
$S(\Sigma_1)_{\alpha}$ to $\prod_{\rho_k\in\Sigma_2(1)}y_k^{a_k+\langle m_,e_k\rangle}$, and the right vertical arrow
is induced by the natural morphism of sheaves $\pi_*O_{\pp_{\Sigma_1}}(D)@>>>O_{\pp_{\Sigma_2}}(\pi_*D)$.
This gives a graded ring homomorphism 
$\pi_*:S(\Sigma_1)@>>>S(\Sigma_2)$ which sends $x_k$ to $y_k$, if $\rho_k\in\Sigma_2(1)$, and sends $x_k$ to 1, otherwise. 

We now apply the above to semiample divisors.
Let $D$ be a semiample (torus invariant) divisor on a complete toric variety $\ps$ in degree $\beta\in A_{d-1}(\ps)$.
on a complete toric variety $\ps$, 
By Theorem~\ref{t:fun}, we get
the associated toric morphism $\pi:{\bf P}_\Sigma\rightarrow{\bf P}_{\Sigma_D}$ such that 
$\pi_*[D\cdot V(\gamma)]$ is ample and $\pi^*\pi_*[D\cdot V(\gamma)]=[D]$ for some closed toric subvariety
$V(\gamma)\subset\ps$, $\gamma\in\Sigma$, which maps birationally onto $\psd$.
In this situation, there is the following natural diagram:
$$\minCDarrowwidth0.7cm
\begin{CD}
S(\Sigma)_{p\beta} @>>> {S(V(\gamma))}_{p\bar\beta} @>>> S(\Sigma_D)_{p\pi_*\bar\beta}\\
@VVV        @VVV             @VVV \\
 H^0(\pp_{\Sigma},O_{\pp_{\Sigma}}(pD))@>>> H^0(V(\gamma),O_{V(\gamma)}(pD_\gamma))@>>> 
H^0(\psd,O_\psd(p\pi_*D_\gamma)),
\end{CD}$$
where $\bar\beta=[D_\gamma]$, $D_\gamma:=D\cdot V(\gamma)$,
 in the Chow group of $V(\gamma)$, and the vertical arrows are isomorphisms.
Since the monomials in $S(\Sigma)_{p\beta}$ and $S(\Sigma_D)_{p\pi_*\bar\beta}$ correspond 
to the lattice points of the same polytope $p\Delta_D$, we get the  isomorphisms 
$$S(\Sigma)_{p\beta}\cong{S(V(\gamma))}_{p\bar\beta}\cong S(\Sigma_D)_{p\pi_*\bar\beta}.$$

\section{Toric hypersurfaces}\label{s:pr}

Here, we  apply the results of the previous section to semiample hypersurfaces in a complete toric variety $\ps$,
which have only transversal intersections with the torus-orbits.  
We also review some  results about  hypersurfaces
in complete simplicial toric varieties.  As a reference we use \cite{m} and \cite{bc}.

A hypersurface $X\subset\ps$ is called {\it $\Sigma$-regular} (or simply {\it regular}) if $X\cap\ts$ is empty or 
a smooth subvariety of codimension 1 in each torus $\ts$ for  $\sigma\in\Sigma$.

By \cite[Proposition~6.8]{d}, a generic hypersurface $X\subset\ps$ of a given semiample degree is $\Sigma$-regular.

\begin{lem}\label{l:conn} Let $X$ be an $i$-semiample hypersurface in a complete toric variety $\ps$ with $i>1$.
 Then $X$  is connected, and $X$  is irreducible if $X$ is $\Sigma$-regular. 
\end{lem}

\begin{pf} The arguments are the same as for Lemma~2.3 in \cite{m}.
\end{pf} 

\begin{rem}\label{r:zer} Let us note that a $0$-semiample hypersurface is always empty because its divisor class is trivial.
\end{rem}

\begin{pr}\label{p:se} 
Let $X$ be a $\Sigma$-regular semiample  hypersurface in a complete toric variety $\ps$, and let
$\pi:\ps@>>>\psx$ be the associated morphism  for $[X]\in A_{d-1}(\ps)$,
then  $Y=\pi(X)$ is a $\Sigma_X$-regular ample hypersurface, and $X=\pi^{-1}(Y)$. 
\end{pr}

\begin{pf}
Start with the case of an $i$-semiample hypersurface with $i>1$. 
{From} Theorem~\ref{t:fun} we have a closed toric subvariety
$V(\gamma)\subset\ps$, for $\gamma\in\Sigma$, which maps birationally onto $\psx$ 
such that $\pi_*[X\cdot V(\gamma)]$ is ample and $\pi^*\pi_*[X\cdot V(\gamma)]=[X]$.
Since $X$ is transversal to the orbits of $\ps$, the divisor class of the hypersurface $X\cap V(\gamma)$ in $V(\gamma)$ is
exactly $[X\cdot V(\gamma)]$. Proposition~\ref{p:comp1} implies that $[X\cdot V(\gamma)]$ is an $i$-semiample divisor
class in $A_{i-1}(V(\gamma))$. The value $i$ is the maximum for the possible Iitaka dimensions of semiample divisors
on the toric variety  $V(\gamma)$. Applying Remark~\ref{r:inm} of the previous section and Proposition~2.4 in  \cite{m},
we get that $\pi(X\cap V(\gamma))$ is a  $\Sigma_X$-regular ample hypersurface.  
On the other hand, by Lemma~\ref{l:conn}, the hypersurface $X$ is irreducible. 
Therefore, its image $Y=\pi(X)$ is also irreducible. Since $\dim \pi(X)\le i$ and $\pi(X\cap V(\gamma))\subset\pi(X)$,
the hypersurface $\pi(X)$ coincides with $\pi(X\cap V(\gamma))$. The hypersurface $Y$ is ample regular and does not intersect
the 0-dimensional orbits. Together with the facts that $X$ and $\pi(X)$ are irreducible this implies
the property $X=\pi^{-1}(Y)$. 

The case of a 1-semiample  hypersurface is special because such a hypersurface is not necessarily connected.
In this situation, we have a closed toric subvariety
$V(\gamma)\subset\ps$, for $\gamma\in\Sigma$, which maps isomorphically onto $\psx\cong {\Bbb P}^1$ 
such that $\pi_*[X\cdot V(\gamma)]$ is ample and $\pi^*\pi_*[X\cdot V(\gamma)]=[X]$.
It follows {from} Proposition~\ref{p:comp1} and Remark~\ref{r:zer} that 
the image $Y=\pi(X)$ is contained in the 1-dimensional torus of $\psx\cong {\Bbb P}^1$.
The preimage $\pi^{-1}(Y)$ of this finite set  
can be easily seen {from} the description of the toric morphism $\pi$ in Section~\ref{s:sem}. This morphism is a trivial
fibration over the 1-dimensional torus of $\psx$, and  
each point of $\pi(X)$ gives exactly one irreducible component  of $\pi^{-1}(Y)$ which is actually a complete toric variety.

On the other hand, each point of $\pi(X)$ came {from} an  irreducible component of $X\subset\pi^{-1}(Y)$.
Hence, $X=\pi^{-1}(Y)$. This gives an isomorphism $\pi:X\cap V(\gamma)\cong\pi(X)$. Thus, $\pi(X)$ is a
$\Sigma_X$-regular ample hypersurface.  
\end{pf}

\begin{rem}            
If, in addition, we assume in this proposition that $X$ is an anticanonical hypersurface in $\ps$, 
then $X$ is big and $\psx$ is a Fano toric variety associated to a reflexive polytope, and this corresponds 
to the construction in \cite{b1}.
\end{rem}

Let $Y$ be an  ample regular hypersurface in a complete toric variety $\ps$.
A hypersurface in the  torus $\ttt\subset\ps$ isomorphic to 
the affine hypersurface $Y\cap{\bf T}$ in  $\ttt$
is   called {\it nondegenerate}. Cohomology of such  hypersurfaces has been studied in \cite{dk} and \cite{b1}.

\begin{lem} \cite{dk}\label{l:lef1}
Let $Z$ be a nondegenerate affine hypersurface in the torus $\ttt$, then
the natural map 
$H^i(\ttt)\rightarrow H^i(Z)$, induced by the inclusion, is an isomorphism  of Hodge structures for $i<\dim\ttt-1$ and 
an injection for $i=\dim\ttt-1$. 
\end{lem}

Using the standard description of a toric morphism, {from} Proposition~\ref{p:se} we get a stratification
of an $i$-semiample regular hypersurface $X\subset\ps$ in terms of nondegenerate affine
hypersurfaces:
\begin{equation}\label{e:str1}
X\cap\ts\cong (\pi(X)\cap{\bf T}_{\sigma_0})\times({\Bbb C}^*)^{l},
\end{equation}
where $\pi:\ps@>>>\psx$ is the associated morphism, $l=d-i+\dim\sigma_0-\dim\sigma$, $d=\dim\ps$,
and $\sigma_0\in\Sigma_X$ is the smallest cone containing the image of $\sigma\in\Sigma$.

{From} here on, we assume that $\pp:=\ps$ denotes a complete simplicial toric variety.
In this case, \cite{bc} shows that homogeneous polynomials in $S:=S(\Sigma)$ determine hypersurfaces in $\pp$.
In  terms of the coordinate ring $S$, a regular 
hypersurface in $\pp$ defined by a homogeneous polynomial $f\in S_\beta$ 
is characterized by the condition that 
$x_1(\partial f/\partial x_1),\dots,x_n(\partial f/\partial x_n)$ do not vanish simultaneously on $\pp$
(see \cite[Proposition~5.3]{c2}). 
A more general class of hypersurfaces in $\pp$ called {\it quasismooth} is defined by a similar condition 
that $\partial f/\partial x_1,\dots,\partial f/\partial x_n$ do not vanish simultaneously on $\pp$ 
(see \cite{bc}).

We also like to mention the following fact.

\begin{pr}\label{p:cal}  An anticanonical quasismooth hypersurface $X$ in a Gorenstein complete simplicial 
toric  variety $\pp$ is  Calabi-Yau.
\end{pr}

\begin{pf} 
A quasismooth hypersurface  is an orbifold (see \cite{bc}),  
and for a $(d-1)$-dimensional orbifold $X$ Calabi-Yau means 
that $\Omega_X^{d-1}\simeq  O_X$ and $H^i(X, O_X)=0$ for $i=1,\dots,d-2$ (see \cite[A.2]{ck}). 
The arguments of the proof that anticanonical implies Calabi-Yau are the same as in \cite{c3}: use the adjunction formula 
$\Omega_X^{d-1}\simeq\Omega^d_\pp(X)\otimes O_X$, the isomorphism $O_\pp(-X)\simeq\Omega^d_\pp$ and
the exact sequence
$0@>>>O_\pp(-X)@>>>O_\pp@>>>O_X@>>>0$.
\end{pf} 

\begin{defn}\cite{bc}\label{d:om}
Fix an integer basis $m_1,\dots,m_d$ for the lattice $M$. Then given subset
$I=\{i_1,\dots,i_d\}\subset\{1,\dots,n\}$, denote 
$\det(e_I)=\det(\langle m_j,e_{i_k}\rangle_{1\le j,k\le d})$, 
$\dd x_I=\dd x_{i_1}\wedge\dots\wedge \dd x_{i_d}$ and $\hat{x}_I=\prod_{i\notin I}x_i$.
Define the $d$-form $\Omega$ by the formula
 $$\Omega=\sum_{|I|=d}\det(e_I)\hat{x}_I \dd x_I,$$
where the sum is over all $d$ element subsets $I\subset\{1,\dots,n\}.$
\end{defn}

Let $X\subset\pp$ be a quasismooth (not necessarily Cartier) hypersurface 
defined by $f\in S_\beta$. For $A\in S_{(a+1)\beta-\beta_0}$ 
(here, $\beta_0=\sum_{i=1}^n\deg(x_i)$), 
consider a rational $d$-form 
$$\omega_A:=A\Omega/f^{a+1}\in H^0(\pp,\Omega^d_{\pp}((a+1)X)).$$ 
This form gives a class in $H^d(\pp\setminus X)$, and by the residue map 
$$\res: H^d(\pp\setminus X)\rightarrow H^{d-1}(X)$$
we get $\res(\omega_A)\in H^{d-1}(X)$. 

\begin{rem}\label{r:resr}
 The residue map and the residues of  rational differential forms with poles along a regular hypersurface
 are well defined even if the toric variety is not simplicial (see the proof of Theorem~3.7 in \cite{dk} and 
Remark~6.4 in \cite{b1}).
\end{rem}

\begin{defn}\cite{bc}\label{d:r1}
Given  $f\in S_\beta$, we have the {\it Jacobian ideal}
$J(f)$ in $S$  
generated by the partial derivatives $\partial f/\partial x_1,\dots, \partial f/\partial x_n$, the ideal
$$J_0(f)=\langle x_1(\partial f/\partial x_1),\dots,x_n(\partial f/\partial x_n)\rangle$$
and the ideal quotient (see \cite[p.~193]{clo})
$J_1(f)=J_0(f):x_1\cdots x_n$.
These give the {\it Jacobian ring} $R(f)=S/J(f)$, $R_0(f)=S/J_0(f)$  and $R_1(f)=S/J_1(f)$ 
graded by the Chow group $A_{d-1}(\pp)$.
\end{defn}

In \cite{m} we have shown that the induced maps 
$$\res(\omega_{\_})^{d-1-q,q}:R(f)_{(q+1)\beta-\beta_0}\rightarrow H^q(X,\Omega^{d-1-q}_X)$$
(sending $A$ to the Hodge component $\res(\omega_{A})^{d-1-q,q}$)
for a quasismooth  hypersurface $X\subset\pp$ and, respectively,
$$\res(\omega_{\_})^{d-1-q,q}:R_1(f)_{(q+1)\beta-\beta_0}\rightarrow H^q(X,\Omega^{d-1-q}_X)$$
for a big and nef regular hypersurface are well defined.
There we also studied the relationship between the multiplicative structure on $R(f)$ (resp., $R_1(f)$)
and the cup product on the middle cohomology of a quasismooth (resp., big and nef regular) hypersurface in $\pp$.
{From} Theorem 4.4 \cite{m} we have the following description of the middle cohomology of  big and nef
regular hypersurfaces $X\subset\ps$:
\begin{equation}\label{e:mid}
H^{d-1-q,q}(X)\cong R_1(f)_{(q+1)\beta-\beta_0}\bigoplus\biggl(\sum_{i=1}^n{\varphi_i}_!H^{d-2-q,q-1}
(X\cap D_i)\biggr),
\end{equation}
where ${\varphi_i}_!$ are the Gysin maps for $\varphi_i:X\cap D_i\hookrightarrow X$. In the case, 
when the dimension of the ambient space is 4 we have (see \cite[Theorem 5.2]{m}):

\begin{thm}\label{t:mid}
Let $X\subset\ps$, $\dim\ps=4$, be  a  big and nef regular hypersurface
defined by $f\in S_\beta$.  
Then there is a natural isomorphism
$$H^{3-q,q}(X)\cong R_1(f)_{(q+1)\beta-\beta_0}
\bigoplus\biggl(\bigoplus\begin{Sb}\sigma
\in\Sigma_X(2)\end{Sb}(R_1(f_\sigma)_{q\beta^\sigma-\beta_0^\sigma})
^{n(\sigma)}\biggr),$$
where $n(\sigma)$ is the number of cones $\rho_i$ such that $\rho_i\subset\sigma$ and $\rho_i\notin
\Sigma_X(1)$, and where $f_\sigma$ is the polynomial of degree $\beta^\sigma$, defining the ample hypersurface
$\pi(X)\cap V(\sigma)\subset V(\sigma)$ (here, $\pi:\ps@>>>\psx$ is the associated morphism),  and $\beta_0^\sigma$
is the degree of the anticanonical divisor on the $2$-dimensional toric variety $V(\sigma)$. 
\end{thm}

\section{Polynomial part of the chiral ring}\label{s:p}

Here we show that for a quasismooth hypersurface $X$ of degree $\beta$ there is a homomorphism between 
$R(f)_{*\beta}$ and the chiral ring $H^*(X,\wedge^*{\cal T}_X)$. 
We will also show that $R_1(f)_{*\beta}$ is  a subring of the chiral ring for a semiample anticanonical regular 
hypersurface $X\subset\pp$ (which is Calabi-Yau). This subring may be called  ``polynomial'' because 
its graded piece in $H^1(X,{\cal T}_X)$ should correspond to polynomial infinitesimal deformations of $X$
performed in the toric variety $\pp$ (see \cite{ck}).

Let $\Omega^p_X$ be the sheaf of {\it Zariski $p$-forms} on an orbifold $X$ (see Appendix A.3 in \cite{ck}).
We can also define $\wedge^p{\cal T}_X:=(\Omega^p_X)^*={Hom}_{O_X}(\Omega^p_X,O_X)$ for an orbifold $X$. 
We call this the {\it (Zariski) $p$-th  exterior power of the tangent sheaf of $X$}.
For $p=1$ this sheaf is isomorphic to  the usual tangent sheaf $\Theta_X$, by Proposition A.4.1 in \cite{ck}. 
When  $X$ is smooth, $\wedge^p{\cal T}_X$   coincides with the standard exterior power sheaf.
Moreover, if $j:X_o\subset X$ is the inclusion of
the smooth locus of $X$, then the argument in the proof of  Proposition 3.10 in \cite{o} shows that 
$j_*(\bigwedge^p\Theta_{X_o})=\wedge^p{\cal T}_X$. One can  use the same argument 
to prove $\Omega^p_X\simeq(\wedge^p{\cal T}_X)^*$ and that $\Omega^p_X$ is isomorphic to the  dual
$(\bigwedge^p{\Theta}_X)^*$ of the usual $p$-th  exterior power of $\Theta_X$, whence 
$\wedge^p{\cal T}_X\simeq(\bigwedge^p{\Theta}_X)^{**}$.
In particular, we also  have the natural maps of sheaves 
$\wedge^p{\cal T}_X\otimes\wedge^q{\cal T}_X\rightarrow\wedge^{p+q}{\cal T}_X$ and
$\wedge^p{\cal T}_X\otimes\Omega^q_X\rightarrow\Omega^{q-p}_X$.

Let $X\subset\pp$ be a quasismooth hypersurface defined by $f\in S_\beta$, which is an orbifold 
as we know {from} \cite{bc}.  By definition of quasismooth, we get an open cover
${\cal U}=\{U_i\}_{i=1}^n$ of $\pp$, where $U_i=\{x\in\pp: f_i(x)\ne0\}$ and $f_i$ denotes 
the partial derivative $\partial f/\partial x_{i}$.

\begin{defn}\label{d:pol} Denote $\partial_{i_0\dots i_p}=\frac{\partial}{\partial x_{i_0}}\wedge\dots\wedge
\frac{\partial}{\partial x_{i_p}}$ for an ordered subset $\{i_0,\dots,i_p\}$ in $\{1,\dots,n\}$. 
Then given $A\in S_{p\beta}$, set
$$(\gamma_A)_{i_0\dots i_p}={\fracwithdelims\{\}{(-1)^{p^2/2}A\langle\partial_{i_0\dots i_p},\dd f\rangle}
{f_{i_0}\cdots f_{i_p}}}_{i_0\dots i_p},$$
where $\langle\,,\,\rangle$ denotes the contraction (the extra factor of $(-1)^{p^2/2}$ which is 
$\sqrt{-1}$ for odd $p$ is added to make convenient commutative diagrams later).
\end{defn}

This  defines a \v{C}ech cocycle, giving its class in $\check H^p({\cal U}|_X,\wedge^p{\cal T}_X)$. Indeed, 
$(\gamma_A)_{i_0\dots i_p}$ 
is homogeneous of degree 0 and  is a cochain in $C^p({\cal U}|_X,\wedge^p{\cal T}_X)$ by the exact sequence
$$0\rightarrow\wedge^p{\cal T}_{X_o} @>>> i^*\wedge^p{\cal T}_{\pp_o} @>i^*\dd f>>
\wedge^{p-1}{\cal T}_{X_o}\otimes O_{X_o}(X)$$
(where $i:X\subset\pp$ is the inclusion, 
$\pp_o$ is the smooth locus of $\pp$ such that $X_o=\pp_o\cap X$ (see \cite[\S4,~p.~55]{hi})),
and because of $\langle\langle\partial_{i_0\dots i_p},\dd f\rangle,\dd f\rangle=0$ since $\dd f\wedge\dd f=0$.
On the other hand, it is straightforward to verify that 
$$(\gamma_A)_{i_0\dots i_p}=\biggl\{\frac{(-1)^{p^2/2}A}{f_{i_0}\cdots f_{i_p}}
\sum_{j=0}^{p}(-1)^j f_{i_j}\partial_{i_0\dots\widehat{i_j}\dots i_p}\biggr\}_{i_0\dots i_p}$$
 vanishes under the \v{C}ech coboundary
map $C^p({\cal U}|_X,\wedge^p{\cal T}_X)\rightarrow C^{p+1}({\cal U}|_X,\wedge^p{\cal T}_X)$. One can 
actually show that  $(\gamma_A)_{i_0\dots i_p}$ is a coboundary in $C^p({\cal U}|_X,i^*\wedge^p{\cal T}_\pp)$.
 
For $A\in S_{p\beta}$ let $\gamma_A\in H^p(X,\wedge^p{\cal T}_X)$ be the image of the \v{C}ech cocycle 
$(\gamma_A)_{i_0\dots i_p}$ under the natural map 
$\check H^p({\cal U}|_X,\wedge^p{\cal T}_X)\rightarrow H^p(X,\wedge^p{\cal T}_X)$.
And we get a well defined map $\gamma_{\_}:R(f)_{*\beta}\rightarrow H^*(X,\wedge^*{\cal T}_X)$
because of the following statement.

\begin{lem}\label{l:cob} If $A\in J(f)_{p\beta}$, then the cocycle $(\gamma_A)_{i_0\dots i_p}$ is a \v{C}ech coboundary in
$C^p({\cal U}|_X,\wedge^p{\cal T}_X)$.
\end{lem}

\begin{pf} If $A\in J(f)_{p\beta}$, then we can assume that $A$ is a multiple of one of the partial 
derivatives $f_k=\partial f/\partial x_{k}$. We have
\begin{multline*}
f_k\frac{\langle\partial_{i_0\dots i_p},\dd f\rangle}{f_{i_0}\cdots f_{i_p}}=
\sum_{j=0}^{p}(-1)^j\frac{f_k \partial_{i_0\dots\widehat{i_j}\dots i_p}}
{f_{i_0}\cdots\widehat{f_{i_j}}\cdots f_{i_p}}
=\sum_{j=0}^{p}(-1)^j\frac{\langle\partial_k,\dd f\rangle\partial_{i_0\dots\widehat{i_j}\dots i_p}}
{f_{i_0}\cdots\widehat{f_{i_j}}\cdots f_{i_p}}
\\
-\sum_{j=0}^{p}(-1)^j\frac{\partial_k\wedge\langle\partial_{i_0\dots\widehat{i_j}\dots i_p},\dd f\rangle}
{f_{i_0}\cdots\widehat{f_{i_j}}\cdots f_{i_p}}
=\sum_{j=0}^{p}(-1)^j\frac{\langle\partial_{k i_0\dots\widehat{i_j}\dots i_p},\dd f\rangle}
{f_{i_0}\cdots\widehat{f_{i_j}}\cdots f_{i_p}},
\end{multline*}
where the second sum after the second equality is identically  zero. 
Hence, it follows that $(\gamma_A)_{i_0\dots i_p}$ is in the image of the \v{C}ech coboundary map
$C^{p-1}({\cal U}|_X,\wedge^p{\cal T}_X)\rightarrow C^{p}({\cal U}|_X,\wedge^p{\cal T}_X)$.
\end{pf}

We now study the compatibility of the multiplication in the Jacobian ring $R(f)$ and the cohomology ring 
$H^*(X,\wedge^*{\cal T}_X)$. The cocycle  $(\gamma_A)_{i_0\dots i_p}$ (up to an extra factor) and the calculations in 
the next two theorems are essentially due to D.~Cox and D.~Morrison.

\begin{thm}\label{t:hom} Let $X\subset\pp$ be a quasismooth hypersurface  defined by $f\in S_\beta$.
The map $R(f)_{*\beta}\rightarrow H^*(X,\wedge^*{\cal T}_X)$, assigning $\gamma_A$ to a polynomial $A$, is 
a ring homomorphism.
\end{thm}

\begin{pf} We need to show that $\gamma_A\cup\gamma_B=\gamma_{AB}$ for 
$A\in S_{p\beta}$ and $B\in S_{q\beta}$.
Similar to \cite[page~63]{cg}, the cup product $\gamma_A\cup\gamma_B$ is represented by the \v{C}ech cocycle 
$$(-1)^{pq}{\fracwithdelims\{\}{(-1)^{p^2/2}A\langle\partial_{i_0\dots i_p},\dd f\rangle\wedge
(-1)^{q^2/2}B\langle\partial_{i_p\dots i_{p+q}},\dd f\rangle}
{f_{i_0}\cdots f_{i_{p+q}}\cdot f_{i_p}}}_{i_0\dots i_{p+q}}.$$ 
Note that 
\begin{multline}\label{e:cup}
\langle\partial_{i_0\dots i_p},\dd f\rangle\wedge\langle\partial_{i_p\dots i_{p+q}},\dd f\rangle
=\sum_{j=0}^{p}(-1)^jf_{i_j} \partial_{i_0\dots\widehat{i_j}\dots i_p}\wedge\sum_{l=p}^{p+q}(-1)^{l-p} f_{i_l}
\partial_{i_p\dots\widehat{i_l}\dots i_{p+q}}\\
=f_{i_p}\sum_{j=0}^{p+q}(-1)^jf_{i_j}  \partial_{i_0\dots\widehat{i_j}\dots i_p}
=f_{i_p}\langle\partial_{i_0\dots i_{p+q}},\dd f\rangle,
\end{multline}
where we used $\partial_{i_p}\wedge\partial_{i_p}=0$. Hence the result follows.
\end{pf}

The middle cohomology of a quasismooth hypersurface $X\subset\pp$ is 
a module over $H^*(X,\wedge^*{\cal T}_X)$ with respect to the natural cup product
$$H^p(X,\wedge^p{\cal T}_X)\otimes H^q(X,\Omega^{d-1-q}_X) @>\cup>>H^{p+q}(X,\Omega^{d-1-p-q}_X).$$
{From} the previous section we know that there is a natural map 
$$\res(\omega_{\_})^{d-1-q,q}:R(f)_{(q+1)\beta-\beta_0}@>>> H^q(X,\Omega^{d-1-q}_X).$$ 
We  normalize this map as
$[\omega_A]=(-1)^{q/2}q!\res(\omega_A)^{d-1-q,q}$ (where we assume $(-1)^{q/2}=(\sqrt{-1})^q$) 
to show that this gives a morphism of modules 
$R(f)_{(*+1)\beta-\beta_0}\rightarrow H^*(X,\Omega^{d-1-*}_X)$.

\begin{thm}\label{t:mod} Let $X\subset\pp$ be a quasismooth hypersurface  defined by $f\in S_\beta$.
Then the diagram 
$$\minCDarrowwidth0.7cm
\begin{CD}
\qquad\qquad R(f)_{p\beta}\otimes R(f)_{(q+1)\beta-\beta_0} @>>>  R(f)_{(p+q+1)\beta-\beta_0} \\
@V{\scriptstyle\gamma_{\_}\otimes[\omega_{\_}]}VV            @V{\scriptstyle[\omega_{\_}]}VV    \\
H^p(X,\wedge^p{\cal T}_X)\otimes H^q(X,\Omega^{d-1-q}_X)@>\cup>>H^{p+q}(X,\Omega^{d-1-p-q}_X)
\end{CD}
$$
commutes, where the top arrow is induced by the multiplication. When $X\subset\pp$ is a $d$-semiample regular 
hypersurface the same diagram commutes with $R_1(f)_{(*+1)\beta-\beta_0}$ in place of $R(f)_{(*+1)\beta-\beta_0}$.
\end{thm}

\begin{pf}
{From} Theorem 3.3 in \cite{m} we know that $[\omega_B]=(-1)^{q/2}q!\res(\omega_B)^{d-1-q,q}$, for 
$B\in S_{(q+1)\beta-\beta_0}$,
is represented by the \v{C}ech cocycle 
$$(-1)^{d-1+(q(q+2)/2)}\fracwithdelims\{\}{BK_{i_q}\cdots K_{i_0}\Omega}{f_{i_0}\cdots f_{i_q}}
_{i_0\dots i_q}\in\check{H}^q({\cal U}|_X,\Omega_X^{d-1-q}),$$
where $K_i$ is the contraction operator 
$(\partial/\partial x_{i})\lrcorner$.
Therefore, for $A\in S_{p\beta}$ the cup product 
$\gamma_A\cup[\omega_B]$ is represented by the \v{C}ech cocycle
$${\biggl\{\frac{(-1)^{p^2/2}A\langle\partial_{i_0\dots i_p},
\dd f\rangle}{f_{i_0}\cdots f_{i_p}}\lrcorner
\frac{(-1)^{d-1+(q(q+2)/2)}BK_{i_{p+q}}\cdots K_{i_p}\Omega}{f_{i_p}\cdots f_{i_{p+q}}}\biggr\}}
_{i_0\dots i_{p+q}}.$$
But note that 
\begin{multline*}  \langle\partial_{i_0\dots i_p},\dd f\rangle\lrcorner K_{i_{p+q}}\cdots K_{i_p}\Omega
=\sum_{j=0}^{p}(-1)^jf_{i_j}  \partial_{i_0\dots\widehat{i_j}\dots i_p}\lrcorner K_{i_{p+q}}\cdots K_{i_p}\Omega
\\
=(-1)^pf_{i_p} \partial_{i_0\dots i_{p-1}}\lrcorner K_{i_{p+q}}\cdots K_{i_p}\Omega=
(-1)^{p+pq}f_{i_p}K_{i_{p+q}}\cdots K_{i_0}\Omega.
\end{multline*}
Since $(-1)^{p^2/2}\cdot(-1)^{q(q+2)/2}\cdot(-1)^{p+pq}=(-1)^{(p+q)(p+q+2)/2}$ we obtain 
$\gamma_A\cup[\omega_B]=[\omega_{AB}]$, whence the diagram commutes.
\end{pf}
     
For an anticanonical quasismooth hypersurface $X$ in a Gorenstein toric variety $\pp$ 
(by Proposition~\ref{p:cal}, $X$ is Calabi-Yau) the situation is especially nice. 
In this case the natural product $\wedge^p{\cal T}_X\otimes\Omega^{d-1}_X\rightarrow\Omega^{d-1-p}_X$
induced by the contraction is an isomorphism since $\Omega^{d-1}_X\simeq O_X$ and 
$\Omega^{d-1-p}_X\simeq Hom_{O_X}(\Omega^p_X,\Omega^{d-1}_X)$ (see \cite[A.3]{ck}), so that the cup product with
$[\omega_1]$ corresponding to $1\in S_0$ ($\beta=\beta_0$ because of anticanonical) gives
\begin{equation}
\cup[\omega_1]:H^p(X,\wedge^p{\cal T}_X)\cong H^{p}(X,\Omega^{d-1-p}_X).\label{f:simeq}
\end{equation}
For regular hypersurfaces this implies: 

\begin{thm}\label{t:subr} Let $X\subset\pp$ be a semiample anticanonical regular  hypersurface  defined by $f\in S_{\beta}$.
Then the map $\gamma_{\_}:R_1(f)_{*\beta}\rightarrow H^*(X,\wedge^*{\cal T}_X)$
 is  an injective ring homomorphism.
\end{thm}

\begin{pf} The map is a well defined  ring homomorphism by Theorems \ref{t:hom}, \ref{t:mod}  and 
(\ref{f:simeq}), while the injectivity follows {from} Theorem 4.4 in \cite{m}.
\end{pf}

Later we will need to use the following result.

\begin{lem}\label{l:proj} Let $j:L@>>> K$ be a morphism of  orbifolds, and let $a\in H^p(K,\wedge^q{\cal T}_K)$ be such that
$\widetilde{j^*}a=\eta^*\tilde a$ for some $\tilde a$ under the maps
$$H^p(K,\wedge^q{\cal T}_K)@>\widetilde{j^*}>>H^p(L,j^*\wedge^q{\cal T}_K)@<\eta^*<<H^p(L,\wedge^q{\cal T}_L).$$
Then $j^*(a\cup b)=\tilde a\cup j^* b$ for $b\in H^{r}(K,\Omega^{s}_K)$.
\end{lem}

\begin{pf} The map $j^*$ decomposes as 
$H^{*}(K,\Omega^{*}_K)@>\widetilde{j^*}>>H^{*}(L,j^*\Omega^{*}_K)@>\eta_*>>H^{*}(L,\Omega^{*}_L)$.
Therefore, 
$$j^*(a\cup b)=\eta_*\widetilde{j^*}(a\cup b)=\eta_*(\widetilde{j^*}a\cup\widetilde{j^*} b)=
\eta_*(\eta^*\tilde a\cup\widetilde{j^*} b)=\tilde a\cup\eta_*\widetilde{j^*} b=\tilde a\cup j^* b,
$$
where we use the projection formula.
\end{pf}

\section{Non-polynomial part of the chiral ring}\label{s:non}

This section studies the non-polynomial  part of the chiral ring which is complementary to the polynomial part.
We will construct new cocycles representing elements in $H^*(X,\wedge^*{\cal T}_X)$ for a big and nef quasismooth
hypersurface $X\subset\ps$. In Section~\ref{s:cal3} we will 
see that these elements with $R_1(f)_\beta$   span $H^1(X,{\cal T}_X)$ for a semiample anticanonical 
regular hypersurface
$X\subset\ps$ ($\dim\ps\ne1,3$). This means that we have found all cocycles corresponding to non-polynomial infinitesimal
deformations for a minimal Calabi-Yau $X$ (see \cite{ck}).

Let $X$ be a $d$-semiample quasismooth hypersurface, defined by $f\in S_\beta$, 
in a complete simplicial toric variety $\ps$ of dimension    
$d$. Then, {from} 
Proposition~\ref{p:se}  we get the associated toric morphism  
$\pi:{\bf P}_\Sigma\rightarrow{\bf P}_{\Sigma_X}$.
Take a 2-dimensional cone $\sigma\in\Sigma_X$ with at least one 1-dimensional cone $\rho_i\subset\sigma$
such that $\rho_i\notin\Sigma_X(1)$. 
Using such a cone $\sigma$ we can form a new cover of the toric variety $\ps$ by the open sets
$$U_{\sigma'}=\biggl\{x\in\ps: \prod_{\rho_k\subset\sigma\setminus\sigma'}x_k\ne0\biggr\}$$
for all 2-dimensional cones $\sigma'\in\Sigma$ that lie in $\sigma$.
Let us fix one order for this open cover corresponding to as the cones lie inside $\sigma$:
\begin{equation}
\setlength{\unitlength}{1cm}
\begin{picture}(8,3.5)
\put(4.8,3.3){$\rho_{l_0}$}
\put(2,1.9){\line(2,1){2.7}}
\put(4.3,2.9){$\sigma_1$}
\put(4.8,2.8){$\rho_{l_1}$}
\put(2,1.9){\line(3,1){2.7}}
\put(4.3,2.5){$\sigma_2$}
\put(2,1.9){\line(6,1){2.7}}
\put(4.8,2.34){$\rho_{l_2}$}
\multiput(4.5,2.2)(0,-.1){3}{\circle*{0.01}}
\put(4.8,1.88){$\rho_{l_{k-1}}$}
\put(4.3,1.7){$\sigma_{k}$}
\put(2,1.9){\line(1,0){2.8}}
\put(4.8,1.4){$\rho_{l_k}$}
\put(2,1.9){\line(6,-1){2.8}}
\put(4,1.3){$\sigma_{k+1}$}
\put(2,1.9){\line(3,-1){2.8}}
\put(4.8,1){$\rho_{l_{k+1}}$}
\put(5.8,2.4){\LARGE$\sigma$}
\multiput(4.5,1)(0,-.1){3}{\circle*{0.01}}
\put(4.78,0.52){$\rho_{l_{n(\sigma)}}$}
\put(2,1.9){\line(2,-1){2.7}}
\put(2,1.9){\line(5,-3){2.73}}
\put(4.8,0.1){$\rho_{l_{n(\sigma)+1}}$}
\end{picture}\label{e:pic}
\end{equation}
where $n(\sigma)$ is the number of cones $\rho_i$ such that $\rho_i\subset\sigma$ and $\rho_i\notin
\Sigma_X(1)$.

Now we take a refinement $U_{i,\sigma_j}=U_i\cap U_{\sigma_j}$ 
of this open cover and the open cover ${\cal U}=\{U_i\}_{i=1}^n$ {from} the previous section.
Denote the refined cover ${\cal U}^\sigma$, considering the order on this cover as the lexicographic order for 
the pairs of indices  $({i,j})$.

\begin{defn}\label{d:non} Given $\rho_i\subset\sigma\in\Sigma_X(2)$ such that $\rho_i\notin\Sigma_X(1)$,
then, as in (\ref{e:pic}), $i=l_k$ for some $k$, and we set 
$$\partial^i_{k}=\frac{x_{l_{k-1}}\partial_{l_{k-1}}}{\mult(\sigma_k)},\quad
\partial^i_{k+1}=-\frac{x_{l_{k+1}}\partial_{l_{k+1}}}{\mult(\sigma_{k+1})},
\text{ and }\partial^i_{j}=0\text{  for }j\ne k,k+1.$$
For $A\in S_{\beta^\sigma_1}$ (here, $\beta^\sigma_1:=\sum_{\rho_k\subset\sigma}\deg(x_k)$),
define 
$$
(\gamma^i_A)_{(i_0,{j_0}),(i_1,{j_1})}={\biggl\{\frac{A}{\prod_{\rho_k\subset\sigma} x_k}
\biggl(\frac{\langle\partial_{i_1}\wedge\partial^i_{j_1},\dd f\rangle}{f_{i_1}}-
\frac{\langle\partial_{i_0}\wedge\partial^i_{j_0},\dd f\rangle}{f_{i_0}}\biggr)\biggr\}}
_{(i_0,{j_0}),(i_1,{j_1})}.
$$ 
\end{defn}

\begin{lem}\label{l:non} In the definition, $(\gamma^i_A)_{(i_0,{j_0}),(i_1,{j_1})}$ is a \v{C}ech cocycle in
$C^1({\cal U^\sigma}|_X,{\cal T}_X)$.
\end{lem}

\begin{pf}
By the arguments after Definition~\ref{d:pol}, $(\gamma^i_A)_{(i_0,{j_0}),(i_1,{j_1})}$ is
 a  cocycle class in $\check H^1({\cal U^\sigma}|_X,{\cal T}_X)$.
The only thing that we need to check in addition is that it is well defined on the given cover, which
follows easily {from} the following two observations. Let $X$ be equivalent to a torus invariant divisor
$D=\sum_{k=1}^na_kD_k$ with the associated polytope $\Delta_D$ and the support function $\psi_D$.
Since $\psi_D$ is linear on $\sigma$ and determines $a_k$,
a monomial $\prod_{k=1}^nx_k^{a_k+\langle m,e_k\rangle}$ (in $x_{l_j}f_{l_j}$)
with $m\in\Delta_D$ is divisible by 
$x_{l_j}$ implies that $a_k+\langle m,e_k\rangle>0$ for all
$\rho_k\subset\sigma$ such that $\rho_k\notin\Sigma_X(1)$. 
In particular, such a  monomial is divisible by $x_i$.
On the other hand, we  have an identity on $\ps$:
\begin{equation}\label{e:eur}
\frac{x_{l_{k-1}}\partial_{l_{k-1}}}{\mult(\sigma_k)}+
\frac{x_{l_{k+1}}\partial_{l_{k+1}}}{\mult(\sigma_{k+1})}
=\frac{\mult(\sigma_k+\sigma_{k+1})}{\mult(\sigma_k)\mult(\sigma_{k+1})}x_{l_k}\partial_{l_k},
\end{equation}
where  $\sigma_k$ and $\sigma_{k+1}$ are  the two cones contained in $\sigma$ and  containing $\rho_i$
(the identity corresponds to an Euler vector field (see \cite[Remark~3.10]{bc}) coming {from} the relation of the cone
generators $\mult(\sigma_{k+1})e_{l_{k-1}}+\mult(\sigma_k)e_{l_{k+1}}=\mult(\sigma_k+\sigma_{k+1})e_{l_k}$
(see \cite[Section~8.2]{d})). 
\end{pf} 

\begin{rem}
Finding  the above cocycle is far {from} obvious, but Propositions~\ref{p:gysin}, \ref{p:relres}, \ref{p:2iso} with 
Theorem~\ref{t:nonm}~and~equation (\ref{f:simeq}) show  how this comes up in the case of 
Calabi-Yau threefolds {from} the description of the middle cohomology in Theorem~\ref{t:mid}.
\end{rem}

Next we generalize the cocycles {from} Definition~\ref{d:non}.

\begin{defn}\label{d:nonc} Let $\rho_i\subset\sigma\in\Sigma_X(2)$ be such that $\rho_i\notin\Sigma_X(1)$.
Given $A\in S_{(p-1)\beta+\beta^\sigma_1}$, $\beta^\sigma_1=\sum_{\rho_k\subset\sigma}\deg(x_k)$, and 
an index set 
$I=\{(i_0,{j_0}),\dots,(i_{p},{j_{p}})\}$,
define 
$$
(\gamma^i_A)_I={\Biggl\{\frac{(-1)^{(p-1)^2/2}A}{\prod_{\rho_k\subset\sigma} x_k}
\sum_{\tilde{I}=I\setminus\{(i_k,j_k)\}}(-1)^k\frac{\langle\partial_{\tilde i_0\dots\tilde i_{p-1}}
\wedge\partial^i_{\tilde j_{p-1}},\dd f\rangle}
{f_{\tilde i_0}\cdots f_{\tilde i_{p-1}}}\Biggr\}}_I,
$$ 
where the sum is over the ordered sets 
$$\tilde{I}=\{(\tilde i_0,{\tilde j_0}),\dots,
(\tilde i_{p-1},{\tilde j_{p-1}})\}=\{(i_0,{j_0}),\dots,\widehat{(i_k,{j_k})},\dots,(i_{p},{j_{p}})\}.
$$
\end{defn}

 Similar to the proof of Lemma~\ref{l:non}, this also determines a  cocycle class in 
$\check H^p({\cal U^\sigma}|_X,\wedge^p{\cal T}_X)$. Denoting its image in
$H^p(X,\wedge^p{\cal T}_X)$ by $\gamma^i_A$, we get a map 
$$\gamma^i_{\_}:S_{(p-1)\beta+\beta^\sigma_1}\rightarrow 
H^p(X,\wedge^p{\cal T}_X),$$
when $\rho_i\setminus\{0\}$ lies in the relative interior of a 2-dimensional cone $\sigma\in\Sigma_X$.

\begin{lem}\label{l:cobi} If $A\in\langle J(f), x_i\rangle_{(p-1)\beta+\beta^\sigma_1}$ and  $p>1$ or 
$A\in\langle x_i\rangle_{\beta^\sigma_1}$, then $\gamma^i_A=0$.
\end{lem}

\begin{pf} If $A$ is divisible by $x_i$, then $(\gamma^i_A)_I$ is clearly a \v{C}ech coboundary, by Definition~\ref{d:nonc}.
Assume $p>1$ and $A\in J(f)_{(p-1)\beta+\beta^\sigma_1}$ is a multiple of one of the partial derivatives
$f_s$. Similar to the proof of Lemma~\ref{l:cob},  we have 
\begin{multline*}
\frac{f_s\langle\partial_{\tilde i_0\dots\tilde i_{p-1}}
\wedge\partial^i_{\tilde j_{p-1}},\dd f\rangle}
{(\prod_{\rho_k\subset\sigma} x_k)f_{\tilde i_0}\cdots f_{\tilde i_{p-1}}}
\equiv
\sum_{l=0}^{p-1}(-1)^l \frac{\langle\partial_{s\tilde i_0\dots\widehat{\tilde i_l}\dots\tilde i_{p-1}}
\wedge\partial^i_{\tilde j_{p-1}},\dd f\rangle}
{(\prod_{\rho_k\subset\sigma} x_k) f_{\tilde i_0}\cdots\widehat{f_{\tilde i_l}}\cdots f_{\tilde i_{p-1}}}
\\
\equiv
\sum_{\Tilde{\Tilde{I}}=\Tilde I\setminus\{(\Tilde i_l,\Tilde j_l)\}}
(-1)^l \frac{\langle\partial_{s{\Tilde{\Tilde i}}_0\dots\Tilde{\Tilde i}_{p-2}}
\wedge\partial^i_{\Tilde{\Tilde j}_{p-2}},\dd f\rangle}
{(\prod_{\rho_k\subset\sigma} x_k) f_{\Tilde{\Tilde i}_0}\cdots f_{\Tilde{\Tilde i}_{p-2}}}
\end{multline*}
(the sum is over the ordered sets $\Tilde{\Tilde I}=\{(\tilde i_0,{\tilde j_0}),\dots,\widehat{(\tilde i_l,{\tilde j_l})},\dots,
(\tilde i_{p-1},{\tilde j_{p-1}})\}$)
 modulo well defined expressions on the open set 
$U_{\tilde i_0,\sigma_{\tilde j_0}}\cap\cdots\cap U_{\tilde i_{p-1},\sigma_{\tilde j_{p-1}}}$, 
because $\langle\partial^i_{\tilde j_{p-1}},\dd f\rangle$ is divisible by $x_i$ and because of equation
(\ref{e:eur}).
On the other hand, there is 
an identity
$$\sum_{\tilde{I}=I\setminus\{(i_k,j_k)\}}(-1)^k\sum_{\Tilde{\Tilde{I}}=\Tilde I\setminus\{(\Tilde i_l,\Tilde j_l)\}}(-1)^l
\frac{\langle\partial_{s\Tilde{\Tilde i}_0\dots\Tilde{\Tilde i}_{p-2}}
\wedge\partial^i_{\Tilde{\Tilde j}_{p-2}},\dd f\rangle}
{f_{\Tilde{\Tilde i}_0}\cdots f_{\Tilde{\Tilde i}_{p-2}}}=0,$$
since the square of a coboundary map is zero. This shows that $(\gamma^i_A)_I$ is a \v{C}ech coboundary for $A\in J(f)$.
\end{pf} 

\begin{defn}
Given  $f\in S_\beta$, let
$J^i(f)$ be the ideal  in $S$  
generated by the Jacobian ideal $J(f)$ and $x_i$. Then we get the quotient ring  $R^i(f)=S/J^i(f)$  
graded by the Chow group $A_{d-1}(\ps)$.
\end{defn}

Lemma~\ref{l:cobi} shows that there are well defined maps 
$\gamma^i_{\_}:R^i(f)_{(p-1)\beta+\beta^\sigma_1}\rightarrow 
H^p(X,\wedge^p{\cal T}_X)$, for $p>1$, and 
$\gamma^i_{\_}:(S/\langle x_i\rangle)_{\beta^\sigma_1}\rightarrow H^1(X,{\cal T}_X)$. 
Note, however, that  a monomial 
$\prod_{\rho_l\not\subset\sigma} x_l\prod_{l=1}^nx_l^{(p-1)a_l+\langle m,e_l\rangle}$ 
in $\langle x_k\rangle_{(p-1)\beta+\beta^\sigma_1}$ (with $\rho_k\subset\sigma$)
 corresponds to $m\in M$ satisfying the inequalities  
$(p-1)a_l+\langle m,e_l\rangle\ge-1$ for $\rho_l\subset\sigma$, $l\ne k$, and $(p-1)a_k+\langle m,e_k\rangle\ge0$.
Since the support function, corresponding to $\beta=[\sum_{i=1}^na_iD_i]$, is
linear on $\sigma$ and determines $a_i$, it follows {from} a relation of the cone generators that 
$(p-1)a_i+\langle m,e_i\rangle\ge0$ and, consequently, the above monomial is divisible by $x_i$, 
for all
$\rho_i\subset\sigma$ such that $\rho_i\notin\Sigma_X(1)$.
Therefore, for all such $\rho_i$ the ideal  $J^i(f)$ is the same as 
$$J^\sigma(f):=\langle J(f), x_k:\rho_k\subset\sigma\rangle$$
 in the degree ${(p-1)\beta+\beta^\sigma_1}$.
Hence, we define $R^\sigma(f)=S/J^\sigma(f)$.

The cocycle $(\gamma^i_A)_I$ in Definition~\ref{d:nonc} came {from} the proof of the following theorem.

\begin{thm}\label{t:ncom}
Let $X\subset\ps$ be a $d$-semiample quasismooth hypersurface defined by $f\in S_\beta$. Then, for $q>1$,
the diagram
$$
\minCDarrowwidth0.7cm
\begin{CD}
\qquad R(f)_{p\beta}\otimes R^\sigma(f)_{(q-1)\beta+\beta^\sigma_1} @>>>  
R^\sigma(f)_{(p+q-1)\beta+\beta^\sigma_1} \\
@V{\scriptstyle\gamma_{\_}\otimes\gamma^i_{\_}}VV            @V{\scriptstyle\gamma^i_{\_}}VV    \\
H^p(X,\wedge^p{\cal T}_X)\otimes H^q(X,\wedge^q{\cal T}_X)@>\cup>>H^{p+q}(X,\wedge^{p+q}{\cal T}_X)
\end{CD}
$$
commutes, where $\beta^\sigma_1=\sum_{\rho_k\subset\sigma}\deg(x_k)$ and 
the top arrow is induced by the multiplication. For $q=1$ the diagram commutes with 
$(S/\langle x_i\rangle)_{\beta^\sigma_1}$
in place of $R^\sigma(f)_{\beta^\sigma_1}$. 
\end{thm}

\begin{pf}  For simplicity, we just  show that if 
$A\in S_{p\beta}$ and 
$B\in S_{\beta^\sigma_1}$,
then $\gamma_A\cup\gamma^i_B=\gamma^i_{AB}$ (the general case is similar though more complicated to write out).
For such $A$ and $B$ the cup product $\gamma_A\cup\gamma^i_B$ is represented by the \v{C}ech cocycle
$$
(-1)^p{\Biggl\{\frac{(-1)^{p^2/2}AB}{\prod_{\rho_k\subset\sigma} x_k}\frac{\langle\partial_{i_0\dots i_p},
\dd f\rangle}{f_{i_0}\cdots f_{i_p}}\wedge 
\biggl(\frac{\langle\partial_{i_{p+1}}\wedge\partial^i_{j_{p+1}},\dd f\rangle}{f_{i_{p+1}}}-
\frac{\langle\partial_{i_{p}}\wedge\partial^i_{j_{p}},\dd f\rangle}{f_{i_{p}}}\biggr)\Biggr\}}_I,
$$ 
where $I=\{(i_0,{j_0}),\dots,(i_{p+1},{j_{p+1}})\}$.
Compute
\begin{multline*}
\langle\partial_{i_0\dots i_p},\dd f\rangle\wedge
\langle\partial_{i_{p+1}}\wedge\partial^i_{j_{p+1}},\dd f\rangle
=\sum_{k=0}^p(-1)^k f_{i_k}\bigl(\partial_{i_0\dots\widehat{i_k}\dots i_p}\wedge
\langle\partial_{i_{p+1}}\wedge\partial^i_{j_{p+1}},\dd f\rangle
\\
+(-1)^p\langle\partial_{i_0\dots\widehat{i_k}\dots i_p},\dd f\rangle
\wedge\partial_{i_{p+1}}\wedge\partial^i_{j_{p+1}}\bigr)
=(-1)^p\sum_{k=0}^p(-1)^k f_{i_k}\langle\partial_{i_0\dots\widehat{i_k}\dots i_{p+1}}
\wedge\partial^i_{j_{p+1}},\dd f\rangle,
\end{multline*}
where the sum of the second terms in the first equality is identically equal to zero.
On the other hand, similar to (\ref{e:cup}),
$$\langle\partial_{i_0\dots i_p},\dd f\rangle\wedge
\langle\partial_{i_{p}}\wedge\partial^i_{j_{p}},\dd f\rangle=
f_{i_p}\langle\partial_{i_0\dots i_{p}}\wedge\partial^i_{j_{p}},\dd f\rangle.$$
Hence, the result follows easily.
\end{pf}
      
We next show when  the cup product of two cocycles  $(\gamma^i_A)_I$ and $(\gamma^j_B)_J$  vanishes.

\begin{lem}\label{l:van} 
The cup product $\gamma^i_{\_}\cup\gamma^j_{\_}=0$ if $\rho_i,\rho_j\subset\sigma\in\Sigma_X(2)$, $i\ne j$,
do not span   a $2$-dimensional cone of $\Sigma$.
\end{lem}
           
\begin{pf} For simplicity, 
we assume  that $\gamma^i_A$ and $\gamma^j_B$ are {from}
$H^1(X,{\cal T}_X)$.

The cup product $\gamma^i_A\cup\gamma^j_{B}$ is represented by the \v{C}ech cocycle
$$(-1){\biggl\{\frac{AB}{(\prod_{\rho_k\subset\sigma} x_k)^2}
\biggl(\frac{u_{i_1,j_1}^i}{f_{i_1}}-
\frac{u_{i_0,j_0}^i}{f_{i_0}}\biggr)\wedge
\biggl(\frac{u_{i_2,j_2}^j}{f_{i_2}}-
\frac{u_{i_1,j_1}^j}{f_{i_1}}\biggr)\biggr\}}_I,$$ 
where $u_{i_k,j_k}^s$ denotes $\langle\partial_{i_k}\wedge\partial^s_{j_k},\dd f\rangle$ for $s\in\{i,j\}$,
and $I=\{(i_0,{j_0}),(i_1,{j_1}),(i_2,{j_2})\}$.
Note that $u_{i_1,j_1}^i\wedge u_{i_1,j_1}^j=0$ because  either $\partial^i_{j_1}$ or $\partial^j_{j_1}$ vanishes
since the corresponding cone $\sigma_{j_1}\subset\sigma$ can not contain both $\rho_i$ and $\rho_j$, 
by the given condition. On the other hand, 
$$\frac{AB\langle\partial_{i_0}\wedge\partial^i_{j_0},\dd f\rangle\wedge
\langle\partial_{i_1}\wedge\partial^j_{j_1},\dd f\rangle}
{(\prod_{\rho_k\subset\sigma} x_k)^2f_{i_0}\cdot f_{i_1}}$$
is well defined on the open set $U_{i_0}\cap U_{\sigma_{j_0}}\cap U_{i_1}\cap U_{\sigma_{j_1}}$, again
because of the given condition.
Hence,  the above cocycle vanishes in the cohomology, being  the image of 
$$(-1){\biggl\{\frac{AB u_{i_0,j_0}^i\wedge u_{i_1,j_1}^j}
{(\prod_{\rho_k\subset\sigma} x_k)^2f_{i_0}\cdot f_{i_1}}\biggr\}}_{(i_0,{j_0}),(i_1,{j_1})}$$
under the \v{C}ech coboundary map 
$C^1({\cal U^\sigma}|_X,\wedge^2{\cal T}_X)\rightarrow C^2({\cal U^\sigma}|_X,\wedge^2{\cal T}_X)$.
\end{pf}

We created  the cocycles $(\gamma^i_A)_I$, now we define the corresponding elements in the middle cohomology 
$H^{d-1}(X)$
of a $d$-semiample quasismooth hypersurface $X$. 

\begin{defn}\label{d:nonm} Let $\rho_i\subset\sigma\in\Sigma_X(2)$ be such that $\rho_i\notin\Sigma_X(1)$.
Given $A\in S_{p\beta-\beta_0+\beta^\sigma_1}$ (where $\beta_0=\sum_{k=1}^n\deg(x_k)$,
 $\beta^\sigma_1=\sum_{\rho_k\subset\sigma}\deg(x_k)$) and 
an index set 
$I=\{(i_0,{j_0}),\dots,(i_{p},{j_{p}})\}$,
define 
$$
(\omega_A^i)_I={\Biggl\{\frac{(-1)^{d+((p-1)^2/2)}A}{\prod_{\rho_k\subset\sigma} x_k}
\sum_{\tilde{I}=I\setminus\{(i_k,j_k)\}}(-1)^k \frac{K_{\tilde i_{p-1}}\cdots K_{\tilde i_{0}}
(\partial^i_{\tilde j_{0}}\lrcorner\Omega)}
{f_{\tilde i_0}\cdots f_{\tilde i_{p-1}}}\Biggr\}}_I,
$$ 
where the sum is over the ordered sets 
$$\tilde{I}=\{(\tilde i_0,{\tilde j_0}),\dots,
(\tilde i_{p-1},{\tilde j_{p-1}})\}=\{(i_0,{j_0}),\dots,\widehat{(i_k,{j_k})},\dots,(i_{p},{j_{p}})\}.
$$
\end{defn}

This determines a \v{C}ech cocycle class  in $\check H^p({\cal U^\sigma}|_X,\Omega^{d-1-p}_X)$, whose image in 
$H^p(X,\Omega^{d-1-p}_X)$ is denoted by $\omega_A^i$.

\begin{lem}\label{l:cobi1} If $A\in J^i(f)_{p\beta-\beta_0+\beta^\sigma_1}$, then $\omega^i_A=0$.
\end{lem}

\begin{pf} If $A$ is divisible by $x_i$, then, by Definition~\ref{d:nonm}, $(\omega^i_A)_I$ is a \v{C}ech coboundary.
Assume that $A\in J(f)$ is a multiple of one of the partial derivatives
$f_s$.

First, consider the case  $p=1$. If $\rho_s\subset\sigma$ and $s\ne i$, then, by the argument
after Definition~\ref{d:non},  $f_s$ is divisible by $x_i$, implying $(\omega^i_A)_I$ is a \v{C}ech coboundary.
The case $f_s=f_i$ is impossible, because of  $S_{\beta_i-\beta_0+\beta^\sigma_1}=0$ ($\beta_i:=\deg(x_i)$),
 following {from} the completeness of the fan $\Sigma$.
The same is true if $\rho_s\not\subset\sigma$ and $\dim\ps>2$. Notice
\begin{multline*}
\frac{f_s{K_{\tilde i_{0}}
(\partial^i_{\tilde j_{0}}\lrcorner\Omega)}}
{(\prod_{\rho_k\subset\sigma} x_k)f_{\tilde i_0}}
=\frac{K_s{K_{\tilde i_{0}}
\partial^i_{\tilde j_{0}}\lrcorner(\dd f\wedge\Omega)}}
{(\prod_{\rho_k\subset\sigma} x_k)f_{\tilde i_0}}
-\frac{\langle\partial^i_{\tilde j_{0}},\dd f\rangle K_s{K_{\tilde i_{0}}\Omega}}
{(\prod_{\rho_k\subset\sigma} x_k)f_{\tilde i_0}}
+\frac{K_s(\partial^i_{\tilde j_{0}}\lrcorner\Omega)}
{(\prod_{\rho_k\subset\sigma} x_k)}.
\end{multline*}
Also, note that if $\dim\ps=2$ and $\rho_s\not\subset\sigma$, 
then $K_s(\partial^i_{\tilde j_{0}}\lrcorner\Omega)$ is a multiple of $x_i$, by the definition
of the form $\Omega$.
Since $\dd f\wedge\Omega\equiv0$ modulo multiples of $f$,
by equation (3) in \cite{m}, and since
$\langle\partial^i_{\tilde j_{0}},\dd f\rangle$ is divisible by $x_i$, it follows that 
$(\omega^i_A)_I$ is a \v{C}ech coboundary in this case.

The  case left is $p>1$. We have 
\begin{multline*}
f_s \frac{K_{\tilde i_{p-1}}\cdots K_{\tilde i_{0}}
(\partial^i_{\tilde j_{0}}\lrcorner\Omega)}
{(\prod_{\rho_k\subset\sigma} x_k)f_{\tilde i_0}\cdots f_{\tilde i_{p-1}}}
=(-1)^{p+1}\frac{K_sK_{\tilde i_{p-1}}\cdots K_{\tilde i_{0}}
\partial^i_{\tilde j_{0}}\lrcorner(\dd f\wedge\Omega)}
{(\prod_{\rho_k\subset\sigma} x_k)f_{\tilde i_0}\cdots f_{\tilde i_{p-1}}}
\\
+(-1)^{p}\frac{\langle\partial^i_{\tilde j_{0}},\dd f\rangle K_sK_{\tilde i_{p-1}}\cdots K_{\tilde i_{0}}
\Omega}
{(\prod_{\rho_k\subset\sigma} x_k)f_{\tilde i_0}\cdots f_{\tilde i_{p-1}}}
-
\sum_{l=0}^{p-1}(-1)^{p+l} \frac{K_sK_{\tilde i_{p-1}}\cdots\widehat{K_{\tilde i_l}}\cdots K_{\tilde i_{0}}
(\partial^i_{\tilde j_{0}}\lrcorner\Omega)}
{(\prod_{\rho_k\subset\sigma} x_k) f_{\tilde i_0}\cdots\widehat{f_{\tilde i_l}}\cdots f_{\tilde i_{p-1}}}
\\
\equiv
(-1)^{p+1}\sum_{\Tilde{\Tilde{I}}=\Tilde I\setminus\{(\Tilde i_l,\Tilde j_l)\}}
(-1)^l \frac{K_sK_{\Tilde{\Tilde i}_{p-2}}\cdots K_{\Tilde{\Tilde i}_0}
(\partial^i_{\Tilde{\Tilde j}_{0}}\lrcorner\Omega)}
{(\prod_{\rho_k\subset\sigma} x_k) f_{\Tilde{\Tilde i}_0}\cdots f_{\Tilde{\Tilde i}_{p-2}}}
\end{multline*}
(the sum is over the ordered sets $\Tilde{\Tilde I}=\{(\tilde i_0,{\tilde j_0}),\dots,\widehat{(\tilde i_l,{\tilde j_l})},\dots,
(\tilde i_{p-1},{\tilde j_{p-1}})\}$)
 modulo well defined expressions on the open set 
$U_{\tilde i_0,\sigma_{\tilde j_0}}\cap\cdots\cap U_{\tilde i_{p-1},\sigma_{\tilde j_{p-1}}}\cap X$, 
because $\dd f\wedge\Omega\equiv0$ modulo multiples of $f$,
$\langle\partial^i_{\tilde j_{p-1}},\dd f\rangle$ is divisible by $x_i$ and because of equation
(\ref{e:eur}).
And, we also have
an identity
$$\sum_{\tilde{I}=I\setminus\{(i_k,j_k)\}}(-1)^k\sum_{\Tilde{\Tilde{I}}=\Tilde I\setminus\{(\Tilde i_l,\Tilde j_l)\}}(-1)^l
\frac{K_sK_{{\Tilde{\Tilde i}}_{p-2}}\cdots K_{{\Tilde{\Tilde i}}_0}
(\partial^i_{{\Tilde{\Tilde j}}_0}\lrcorner\Omega)}
{f_{\Tilde{\Tilde i}_0}\cdots f_{\Tilde{\Tilde i}_{p-2}}}=0,$$
since the square of a coboundary map is zero. Hence, $(\gamma^i_A)_I$ is a \v{C}ech coboundary if $A\in J(f)$.
\end{pf}

The last lemma shows that there is a  well defined map
$$\omega^i_{\_}: R^i(f)_{p\beta-\beta_0+\beta^\sigma_1}\rightarrow  H^{p}(X,\Omega^{d-1-p}_X).$$ 
Since $p\beta$ is $d$-semiample, multiplying a monomial in 
$\langle x_k\rangle_{p\beta-\beta_0+\beta^\sigma_1}$ (for $\rho_k\subset\sigma$)
by $\prod_{\rho_l\not\subset\sigma} x_l$ and applying the argument in the proof of Lemma~\ref{l:non}, 
we get a monomial divisible by all $x_i$ corresponding to 
$\rho_i\subset\sigma$ such that $\rho_i\notin\Sigma_X(1)$. 
Therefore, for all such $\rho_i$ the ideal  $J^i(f)$ is the same as 
$J^\sigma(f)$ in the degree ${p\beta-\beta_0+\beta^\sigma_1}$.
 
The cocycle $(\omega_A^i)_I$ came {from} the proof of the following result.

\begin{thm}\label{t:nonm}
Let $X\subset\ps$ be a $d$-semiample quasismooth hypersurface defined by $f\in S_\beta$. Then, for $p>1$,
the diagram 
$$\minCDarrowwidth0.7cm
\begin{CD}
R^\sigma(f)_{(p-1)\beta+\beta^\sigma_1}\otimes R(f)_{(q+1)\beta-\beta_0} @>>> R^\sigma(f)_{(p+q)\beta-\beta_0+\beta^\sigma_1} \\
@V{\scriptstyle\gamma^i_{\_}\otimes[\omega_{\_}]}VV                @V{\scriptstyle \omega^i_{\_}}VV    \\
\quad H^p(X,\wedge^p{\cal T}_X)\otimes H^q(X,\Omega^{d-1-q}_X) @>\cup>> H^{p+q}(X,\Omega^{d-1-p-q}_X)
\end{CD}
$$
commutes, where the top arrow is the multiplication (for $p=1$ the diagram commutes with 
$(S/\langle x_i\rangle)_{\beta^\sigma_1}$
in place of $R^\sigma(f)_{\beta^\sigma_1}$).
\end{thm}

\begin{pf}
For simplicity, we only 
 show  that $\gamma^i_{A}\cup[\omega_{B}]=\omega^i_{AB}$ for $A\in S_{\beta^\sigma_1}$
and $B\in S_{(q+1)\beta-\beta_0}$ (as in the proof of Theorem~\ref{t:ncom},
the general case is similar, but more complicated to write out).
Similar to the proof of Theorem~\ref{t:mod}, the cup product 
$\gamma^i_{A}\cup[\omega_{B}]$ is represented by the \v{C}ech cocycle
$$
{\biggl\{\frac{(-1)^{d-1+(q(q+2)/2)}AB}{\prod_{\rho_k\subset\sigma} x_k}
\biggl(\frac{\langle\partial_{i_1}\wedge\partial^i_{j_1},\dd f\rangle}{f_{i_1}}-
\frac{\langle\partial_{i_0}\wedge\partial^i_{j_0},\dd f\rangle}{f_{i_0}}\biggr)\lrcorner
\frac{K_{i_{q+1}}\cdots K_{i_1}\Omega}{f_{i_1}\cdots f_{i_{q+1}}}\biggr\}}_I,
$$
where $I$ is the  index set $\{(i_0,{j_0}),\dots,(i_{q+1},{j_{q+1}})\}$, corresponding to the cover 
${\cal U^\sigma}|_X$.
Compute 
\begin{multline*}
\biggl(\frac{\langle\partial_{i_1}\wedge\partial^i_{j_1},\dd f\rangle}{f_{i_1}}-
\frac{\langle\partial_{i_0}\wedge\partial^i_{j_0},\dd f\rangle}{f_{i_0}}\biggr)\lrcorner
\frac{K_{i_{q+1}}\cdots K_{i_1}\Omega}{f_{i_1}\cdots f_{i_{q+1}}}\\
=(-1)^{q+1}\frac{K_{i_{q+1}}\cdots K_{i_1}(\partial^i_{j_1}\lrcorner\Omega)}{f_{i_1}\cdots f_{i_{q+1}}}
-(-1)^{q+1}\frac{K_{i_{q+1}}\cdots K_{i_1}(\partial^i_{j_0}\lrcorner\Omega)}{f_{i_1}\cdots f_{i_{q+1}}}
\\
+(-1)^{q+1}\frac{\langle\partial^i_{j_0},\dd f\rangle K_{i_{q+1}}\cdots K_{i_0}\Omega}{f_{i_0}\cdots f_{i_{q+1}}}.
\end{multline*}
Also, notice
\begin{multline*}
K_{i_{q+1}}\cdots K_{i_0}\partial^i_{j_0}\lrcorner(\dd f\wedge\Omega)=
\langle\partial^i_{j_0},\dd f\rangle K_{i_{q+1}}\cdots K_{i_0}\Omega\\
-f_{i_0}K_{i_{q+1}}\cdots K_{i_1}(\partial^i_{j_0}\lrcorner\Omega)+
\sum^{q+1}_{k=1}(-1)^{k-1} f_{i_k}
K_{i_{q+1}}\cdots\widehat{K_{i_k}}\cdots K_{i_1}K_{i_0}(\partial^i_{j_0}\lrcorner\Omega).
\end{multline*}
Since $\dd f\wedge\Omega\equiv0$ modulo multiples of $f$, as in Lemma~\ref{l:cobi1},
we can see that $\gamma^i_{A}\cup[\omega_{B}]$ is actually represented by the \v{C}ech cocycle
$$
{\Biggl\{\frac{(-1)^{d+(q^2/2)}AB}{\prod_{\rho_k\subset\sigma} x_k}
\biggl(\sum_{\tilde{I}=I\setminus\{(i_k,j_k)\}}(-1)^k \frac{K_{\tilde i_{q}}\cdots K_{\tilde i_{0}}
(\partial^i_{\tilde j_{0}}\lrcorner\Omega)}
{f_{\tilde i_0}\cdots f_{\tilde i_{q}}}\biggr)\Biggr\}}_I,
$$ 
where the sum is over the ordered sets 
$$
\tilde{I}=\{(\tilde i_0,{\tilde j_0}),\dots,
(\tilde i_{q},{\tilde j_{q}})\}=\{(i_0,{j_0}),\dots,\widehat{(i_k,{j_k})},\dots,(i_{q+1},{j_{q+1}})\}.
$$
\end{pf}

The next result (a proof of which is similar to the above) shows that 
the map $\omega^i_{\_}: R^\sigma(f)_{*\beta-\beta_0+\beta^\sigma_1}\rightarrow  H^{*}(X,\Omega^{d-1-*}_X)$
is a morphism of modules with respect to the ring homomorphism $R(f)_{*\beta}\rightarrow H^*(X,\wedge^*{\cal T}_X)$.

\begin{thm}\label{t:nonm1} Let $X\subset\ps$ be a $d$-semiample quasismooth hypersurface  defined by $f\in S_\beta$.
Then the diagram 
$$\minCDarrowwidth0.7cm
\begin{CD}
\qquad\qquad R(f)_{p\beta}\otimes R^\sigma(f)_{q\beta-\beta_0+\beta^\sigma_1} @>>>  R^\sigma(f)_{(p+q)\beta-\beta_0+\beta^\sigma_1} \\
@V{\scriptstyle\gamma_{\_}\otimes\omega^i_{\_}}VV            @V{\scriptstyle\omega^i_{\_}}VV    \\
H^p(X,\wedge^p{\cal T}_X)\otimes H^q(X,\Omega^{d-1-q}_X)@>\cup>>H^{p+q}(X,\Omega^{d-1-p-q}_X)
\end{CD}
$$
commutes, where the top arrow is induced by the multiplication. 
\end{thm}

Similar to Lemma~\ref{l:van}, we also
get  when  the cup product of two cocycles  $(\gamma^i_A)_I$ and $(\omega^j_B)_J$  vanishes.

\begin{lem}\label{l:van1} The cup product $\gamma^i_{\_}\cup\omega^j_{\_}=0$ if 
$\rho_i,\rho_j\subset\sigma\in\Sigma_X(2)$, $i\ne j$,
do not span   a 2-dimensional cone of $\Sigma$. 
\end{lem}

\section{Toric and residue parts of cohomology}\label{s:tr}

In this section we describe the toric part of cohomology 
of a  semiample  regular hypersurface in  a complete simplicial toric variety $\ps$. 
This part is the image of cohomology of the ambient space induced by the inclusion of the hypersurface. 
In this case, we also show that cohomology has a natural 
decomposition into a direct sum of the toric part and the residue part which comes from the residues of rational differential
forms with poles along the hypersurface.

Since $\ps$ is simplicial, we know {from} \cite{f1} that
the cohomology ring  $H^*(\ps)$ (with complex coefficients) is isomorphic to
$${\Bbb C}[D_1,\dots,D_n]/(P(\Sigma)+SR(\Sigma)),$$
where the generators correspond to the torus invariant divisors of $\ps$, 
and where
$$P(\Sigma)=\biggl\langle \sum_{i=1}^n \langle m,e_i\rangle D_i: m\in M\biggr\rangle,$$
$$SR(\Sigma)=\bigl\langle D_{i_1}\cdots D_{i_k}:\{e_{i_1},\dots,e_{i_k}\}\not\subset\sigma 
\text{ for all }\sigma\in\Sigma\bigr\rangle$$
($SR(\Sigma)$ is the {\it Stanley-Reisner} ideal of $\Sigma$).
The  {\it toric  part} $H^*_{\rm toric}(X)$ of  cohomology of a hypersurface $X$ in $\ps$ is defined as
the image of the restriction map $i^*: H^*(\ps)@>>>H^*(X)$ induced by the inclusion $i:X\subset\ps$.

\begin{thm}\label{t:toric}
Let $X$ be a semiample regular hypersurface in a complete simplicial toric variety $\ps$.
Then 
$$H^*_{\rm toric}(X)\cong H^*(\ps)/Ann([X])\cong {\Bbb C}[D_1,\dots,D_n]/I,$$
where $Ann([X])$ is the annihilator of the class $[X]\in H^2(\ps)$, and where $I=(P(\Sigma)+SR(\Sigma)):[X]$ is the ideal quotient.
\end{thm}

\begin{pf} We need to  show that $\ker(i^*:H^*(\ps)@>>>H^*(X))$ coincides with 
$\ker(\cup[X]: H^*(\ps)@>>>H^{*+2}(\ps))$.
 Since $\cup[X]=i_!i^*$ (where $i_!$ is the Gysin map), this is  equivalent to 
$\ker(i_!)\cap\im(i^*)=0$ in $H^p(X)$ for all $p$. Using  an induction on the dimension of the hypersurface,
 we will show  a stronger statement:
\begin{equation}\label{e:dec} 
H^p(X)=\im(i^*)\oplus\ker(i_!)\text{\quad for all }p.
\end{equation}

If $\dim X=0$, then $\ps={\Bbb P}^1$. In this case, the composition $H^0({\Bbb P}^1)@>i^*>> H^0(X)@>i_!>>H^2({\Bbb P}^1)$
is clearly an isomorphism, and (\ref{e:dec}) follows. 

Let $\dim X=d-1>0$. For all odd $p$, $H^p(X)=\ker(i_!)$ and equation (\ref{e:dec}) holds because $H^{odd}(\ps)$ vanishes.
So we can assume that $p$ is even.

We show  first that $H^p(X)=\im(i^*)+\ker(i_!)$. 
The Gysin spectral sequence (see \cite[Section~4]{m}) gives
an exact sequence  
$$\oplus_{k=1}^n H^{p-2}(X\cap D_k) @>>>H^p(X)@>>>\gr_p^W H^p(X\cap\ttt)@>>>0.$$
Also, by the Gysin exact sequence (see \cite[Theorem~3.7]{dk}), we get
\begin{equation}\label{e:gysin}
0@>>>H^{p+1}(\ps\setminus X)@>\res>>H^p(X)@>i_!>>H^{p+2}(\ps)
\end{equation}
for even $p$. Hence, $\res(H^{p+1}(\ps\setminus X))=\ker(i_!)$. We claim that
the composition 
\begin{equation}\label{e:comp}
H^{p+1}(\ps\setminus X)@>\res>>H^p(X)@>>>\gr_p^W H^p(X\cap\ttt)
\end{equation}
is a surjective map for $p>0$.
If $[X]$ is an $i$-semiample divisor class, then we get the  associated  morphism $\pi:\ps@>>>\psx$, and
the ample regular hypersurface $Y=\pi(X)$ in $\psx$, by Proposition~\ref{p:se}.
The statement is trivial for $p\ne i-1$ because, in this case,
\begin{equation}\label{e:grvan}
\gr_p^W H^p(X\cap\ttt)\cong\gr_p^W H^p((Y\cap\ttt_{\Sigma_X})\times ({\Bbb C}^*)^{d-i})=0
\end{equation}
(where $\ttt_{\Sigma_X}$ is the maximal torus of $\psx$), 
by equation (\ref{e:str1}) and the K\"unneth isomorphism theorem with Lemma~\ref{l:lef1}. 
For $p=i-1$, 
consider the following commutative diagram:
$$\minCDarrowwidth0.7cm
\begin{CD}
 H^{i}(\ps\setminus X)@>\res>>H^{i-1}(X)@>>>H^{i-1}(X\cap\ttt)\\
@AAA        @AAA       @AAA \\
 H^{i}(\psx\setminus Y)@>\res>>H^{i-1}(Y)@>>>H^{i-1}(Y\cap\ttt_{\Sigma_X}),
\end{CD}
$$
where the vertical arrows are induced by the morphism $\pi$.
The right vertical arrow descends to an isomorphism 
\begin{equation}\label{e:griso}
\pi^*:\gr_{i-1}^W H^{i-1}(Y\cap\ttt_{\Sigma_X})\cong\gr_{i-1}^W H^{i-1}(X\cap\ttt)
\end{equation}
which follows {from} equation (\ref{e:str1}), the  K\"unneth isomorphism and Lemma~\ref{l:lef1}.
On the other hand, the proof of Theorem~4.4 in \cite{m} and Remark~\ref{r:resr} show that the weight space 
$W_{i-1} H^{i-1}(Y\cap\ttt_{\Sigma_X})$ lies in the image of the  composition of maps on the bottom of the diagram.
Thus, we have shown that the composition  (\ref{e:comp}) is surjective  for all $p>0$. Hence,
$\ker(i_!)$ in  $H^p(X)$ maps onto $\gr_p^W H^p(X\cap\ttt)$.
Since $\gr_p^W H^p(\ttt)=0$ for $p>0$, we get the commutative diagram:
$$\minCDarrowwidth0.7cm
\begin{CD}
\bigoplus_{k=1}^n H^{p}(D_k) @>>>H^{p+2}(\ps)@>>>0@.  \\
@AA i_! A           @AA i_! A @. @.   \\
\bigoplus_{k=1}^n H^{p-2}(X\cap D_k) @>>>H^p(X)@>>>\gr_p^W H^p(X\cap\ttt)@>>>0 \\
@AA i^* A           @AA i^* A @.   \\
\bigoplus_{k=1}^n H^{p-2}(D_k) @>>>H^p(\ps)@>>>0, @.
\end{CD}
$$
where the rows are exact sequences arising 
 {from} the Gysin spectral sequence. 
Chasing this diagram and using the induction assumption (\ref{e:dec})
for the semiample regular hypersurfaces $X\cap D_k\subset D_k$, 
we can see that $H^p(X)$ is spanned by $\ker(i_!)$ and $\im(i^*)$ for 
all $p>0$. Let us show this in the case $p=0$.
If $X$ is connected, then 
$i^*:H^0(\ps)@>>> H^0(X)$ is an isomorphism of 1-dimensional spaces,
whence $H^0(X)=\im(i^*)$. 
By Lemma~\ref{l:conn}, we are left to consider the case when $X$ is a $1$-semiample hypersurface.
We use another commutative diagram:
$$\minCDarrowwidth0.7cm
\begin{CD}
H^0(\ps)@>i^*>> H^0(X)@>i_!>>H^2(\ps)\\
@AA\pi^*A @AA\pi^*A @AA\pi^*A \\
H^0(\psx)@>>> H^0(Y)@>>>H^2(\psx).
\end{CD}$$
The property $X=\pi^{-1}(Y)$ {from} Proposition~\ref{p:se} gives an isomorphism $\pi^*:H^0(Y)@>>>H^0(X)$.
Using the diagram and the fact $\psx\cong{\Bbb P}^1$, we deduce  $H^0(X)=\im(i^*)+\ker(i_!)$.

To prove (\ref{e:dec}) it suffices now to show that $\im(i^*)$ and  $\ker(i_!)$ have complementary dimensions
in $H^p(X)$. {From} equation (\ref{e:gysin}) we get $\dim \ker (i_!)=h^{p+1}(\ps\setminus X)$. 
The exact sequence of cohomology with compact supports 
$$H^p(\ps)@>i^*>>H^p(X)@>>>H_c^{p+1}(\ps\setminus X)@>>>0$$
 also gives $\dim\im(i^*)=h^p(X)-h_c^{p+1}(\ps\setminus X)$ for even $p$.
Since $H^p(X)=\im(i^*)+\ker(i_!)$, the inequalities 
\begin{equation}\label{e:ineq} 
h_c^{p+1}(\ps\setminus X)\le h^{p+1}(\ps\setminus X)
\end{equation}
 hold for all even $p$. By Poincar\'e duality, we have the equalities 
$h_c^{p+1}(\ps\setminus X)=h^{2d-p-1}(\ps\setminus X)$, $h^{p+1}(\ps\setminus X)=h_c^{2d-p-1}(\ps\setminus X)$.
Applying them to (\ref{e:ineq}), we get
$$h^{2d-p-1}(\ps\setminus X)\le h_c^{2d-p-1}(\ps\setminus X)$$
for all even $p$. Hence, all these inequalities are equalities, and equation (\ref{e:dec}) follows.
The  proof by induction is finished.
\end{pf}

\begin{rem} We should note that 
the above nontrivial result or its equivalent has been used without a proof for smooth Calabi-Yau hypersurfaces
(complete intersections) in  many papers (e.g., \cite[Proposition~8.1]{b3}, \cite[Section~3.4]{hly},  
\cite[Section~9]{st}; cup product induces a nondegenerate pairing on the toric part---\cite[Lemma~8.6.11]{ck}, \cite[Introduction]{gi}).
In the case of ample quasismooth hypersurfaces, this follows directly {from} the Hard-Lefschetz theorem.
It is an open question whether Theorem~\ref{t:toric} holds in general for smooth or quasismooth semiample hypersurfaces.
\end{rem}

\begin{rem} An interesting equality follows from the proof of Theorem~\ref{t:toric}:
$$h^{p}(\ps\setminus X)=h_c^{p}(\ps\setminus X)\text{ for odd }p.$$ If $X$ is ample, these Hodge numbers vanish for $p$ away from the middle dimension $d$.
But in the semiample case they are nontrivial in general.
\end{rem}

As a consequence of the above proof, we have a direct sum decomposition $H^p(X)=\im(i^*)\oplus\ker(i_!)$ for a semiample regular
hypersurface. By the Gysin exact sequence, the kernel of the Gysin map is exactly the image of the residue map.
Therefore, it is natural to introduce the following.

\begin{defn} The {\it residue part} $H^*_{\rm res}(X)$ of cohomology of a quasismooth hypersurface $X$ in a complete simplicial
toric variety $\ps$ is defined as the image of the residue map $\res:H^{*+1}(\ps\setminus X)@>>>H^*(X)$.
\end{defn} 

\begin{rem} The  residue part $H^*_{\rm res}(X)$ is isomorphic to the primitive cohomology $PH^*(X)$ defined
in \cite{bc}  by the exact sequence
$$H^*(\ps)@>>>H^*(X)@>>>PH^*(X)@>>>0.$$
\end{rem}

For a semiample regular hypersurface  $X$ in a complete simplicial
toric variety $\ps$ we have
$$H^*(X)=H^*_{\rm toric}(X)\oplus H^*_{\rm res}(X).$$
Theorem~\ref{t:toric} described the toric part. 
Note that 
$$H^*_{\rm toric}(X)\cup H^*_{\rm res}(X)\subset H^*_{\rm res}(X),$$
since $i_!(i^*a\cup b)=a\cup i_!b=0$ for $b\in\ker(i_!)$, by the projection formula.
Therefore, the residue part is a submodule of $H^*(X)$ over the ring $H^*_{\rm toric}(X)$.

Finally, we suggest an algorithmic approach to computing the residue part of cohomology.
As in the proof of Theorem~\ref{t:toric}, the Gysin spectral  sequence gives the commutative diagram:
\begin{equation}\label{e:resi}
\minCDarrowwidth0.5cm
\begin{CD}
 0  @. 0  @. 0  @.\\ 
 @AAA           @AAA @AAA   @.\\ 
\bigoplus_{k=1}^n H_{\rm res}^{p-2}(X\cap D_k) @>>>H_{\rm res}^p(X)@>>>\gr_p^W PH^p(X\cap\ttt)@>>>0 \\
 @AAA           @AAA @AAA   @. \\
\bigoplus_{k=1}^n H^{p-2}(X\cap D_k) @>>>H^p(X)@>>>\gr_p^W H^p(X\cap\ttt)@>>>0 \\
 @AAA           @AAA @AAA @.   \\
\bigoplus_{k=1}^n H^{p-2}(D_k) @>>>H^p(\ps)@>>>\gr_p^W H^p(\ttt)@>>>0,
\end{CD}
\end{equation}
where the columns and the rows are exact, and where $PH^p(X\cap\ttt)$ is defined,
as in \cite[Definition~3.13]{b2},  by the exact sequence
$$H^*(\ttt)@>>>H^*(X\cap\ttt)@>>>PH^*(X\cap\ttt)@>>>0.$$
The hypersurfaces $X\cap D_k$ in $D_k$ are semiample regular of lower dimension, and 
the space $\gr_p^W PH^p(X\cap\ttt)$ can be  described in terms of cohomology of a nondegenerate affine hypersurface,
again, using the proof of  Theorem~\ref{t:toric}. Therefore, this provides a way  to calculate  $H_{\rm res}^p(X)$.

\section{Cohomology of semiample regular hypersurfaces}\label{s:cohreg}

In this section we continue the study of the cohomology of  semiample regular hypersurfaces which
was initiated in \cite[Section~4]{m}.
Applying the algorithmic approach of the previous section, we
will  compute the residue  part of  the middle cohomology of  a big and nef regular hypersurface $X$.
In particular, we will generalize  the description in equation (\ref{e:mid}) and Theorem~\ref{t:mid}. 
An algebraic description of the middle cohomology is important  because, in the Calabi-Yau case, this is isomorphic to
the chiral ring $H^*(X,\wedge^*{\cal T}_X)$, by equation (\ref{f:simeq}).
In terms of this description, one should be able to compute the 
product structure of the chiral ring.
Here, we also compute the nontrivial cup products
$\gamma^i_{A}\cup\omega^j_{B}$ of elements constructed in Section~\ref{s:non}.

Let $X$ be a $d$-semiample regular hypersurface, defined by $f\in S_\beta$,
in a complete simplicial toric variety $\ps$.
Our goal is to relate   $\omega_A^i$, defined in Section~\ref{s:non}, to the description of
the middle cohomology of $X$ given in equation (\ref{e:mid}). First, we define  
new \v{C}ech cocycles, representing  elements in  $H^{d-3}(X\cap D_i)$. 

\begin{defn}\label{d:gysc} Given $\sigma\in\Sigma_X(2)$ with the ordered integral generators $e_{l_0}$ and $e_{l_{n(\sigma)+1}}$
as in (\ref{e:pic}),
introduce  a $(d-2)$-form
$$\Omega_\sigma=\frac{x_{l_0}x_{l_{n(\sigma)+1}}K_{l_{n(\sigma)+1}}K_{l_0}\Omega}{\mult(\sigma)\prod_{\rho_k\subset\sigma}x_k}.$$ 
Then, for $A\in S_{(p+1)\beta-\beta_0+\beta_1^\sigma}$ and $\rho_i\subset\sigma$ such that 
$\rho_i\notin\Sigma_X(1)$, define
$$(\widetilde\omega^i_A)_I=(-1)^{p^2/2}\fracwithdelims\{\}{AK_{i_{p}}\cdots K_{i_0}\Omega_\sigma}{f_{i_0}\cdots f_{i_{p}}}
_I,$$
where $I$ is the index set $\{i_0,\dots,i_p\}$.
\end{defn} 

Consider a rational $(d-2)$-form 
$$(A\Omega_\sigma/f^{p+1})\in H^0(D_i,\Omega^{d-2}_{D_i}((p+1)X_i)),$$
where $X_i:=X\cap D_i$ (we will use both notations).
By the residue map we get $\res(A\Omega_\sigma/f^{p+1})\in H^{d-3}(X\cap D_i)$.
The next statement shows that up to a constant, $(\widetilde\omega^i_A)_I$ is a \v{C}ech cocycle which represents this residue.

\begin{pr}\label{p:rel1} 
Let $X\subset\ps$ be a $d$-semiample regular hypersurface 
defined by $f\in S_\beta$. Given $\rho_i\subset\sigma\in\Sigma_X(2)$ such that 
$\rho_i\notin\Sigma_X(1)$, and $A\in S_{(p+1)\beta-\beta_0+\beta_1^\sigma}$,
then, under the natural map
$$\check{H}^p({\cal U}|_{X\cap D_i},\Omega_{X\cap D_i}^{d-3-p})\rightarrow 
H^p(X\cap D_i,\Omega_{X\cap D_i}^{d-3-p})\cong H^{d-3-p,p}(X\cap D_i),$$   
the Hodge component $\res(A\Omega_\sigma/f^{p+1})^{d-3-p,p}$ is represented by the \v{C}ech cocycle
$$\frac{(-1)^{d-3+(p(p+1)/2)}}{p!}\fracwithdelims\{\}{AK_{i_{p}}\cdots K_{i_0}\Omega_\sigma}{f_{i_0}\cdots f_{i_{p}}}_I
\in C^p({\cal U}|_{X\cap D_i},\Omega^{d-3-p}_{X\cap D_i}).$$
\end{pr}
 
\begin{pf}
The proof of this is similar to the proof of Theorem~3.3 in \cite{m} (see also \cite{cg}).
We only need to show that
\begin{equation}\label{e:mul}
{\rm d}f\wedge\Omega_\sigma\equiv0\mbox{ modulo multiples of } f \text{ and } x_i.
\end{equation} 
Note 
\begin{multline*}
\dd f\wedge\Omega_\sigma=\dd f\wedge\frac{x_{l_0}x_{l_{n(\sigma)+1}}K_{l_{n(\sigma)+1}}K_{l_0}\Omega}
{\mult(\sigma)\prod_{\rho_k\subset\sigma}x_k}=\frac{x_{l_0}x_{l_{n(\sigma)+1}}K_{l_{n(\sigma)+1}}K_{l_0}(\dd f\wedge\Omega)}
{\mult(\sigma)\prod_{\rho_k\subset\sigma}x_k}\\
-\frac{x_{l_0}f_{l_0}x_{l_{n(\sigma)+1}}K_{l_{n(\sigma)+1}}\Omega}
{\mult(\sigma)\prod_{\rho_k\subset\sigma}x_k}
+\frac{x_{l_0}x_{l_{n(\sigma)+1}}f_{l_{n(\sigma)+1}}K_{l_0}\Omega}
{\mult(\sigma)\prod_{\rho_k\subset\sigma}x_k}.
\end{multline*}
The first summand is divisible by $f$, because $\dd f\wedge\Omega\equiv0$ modulo multiples of $f$,  
as in Lemma~\ref{l:cobi1}, and because $f$ is not divisible by any variable $x_k$, corresponding 
to $\rho_k\subset\sigma$, since $X$ is   regular. 
The sum of the other two terms 
is a multiple of $x_i$, because, by the argument after
Definition~\ref{d:non},  $x_{l_j}f_{l_j}$
are divisible by all variables $x_k$, 
corresponding to the cones $\rho_k\subset\sigma$ not contained in $\Sigma_X(1)$, and 
because of an Euler identity similar to (\ref{e:eur}). Hence, equation (\ref{e:mul}) follows.

We also verify that $(\widetilde\omega^i_A)_I$  is a \v{C}ech cocycle.
The \v{C}ech coboundary of $(\widetilde\omega^i_A)_I$ is
$$(-1)^{p^2/2}\biggl\{A\sum_{k=0}^{p+1}\frac{(-1)^kf_{i_k}K_{i_{p+1}}\cdots\widehat{K_{i_k}}\cdots 
K_{i_0}\Omega_\sigma}{f_{i_0}\cdots f_{i_{p+1}}}\biggr\}_I.$$
On the other hand,
\begin{multline*}
\sum^{p+1}_{k=0}(-1)^k f_{i_k}
K_{i_{p+1}}\cdots\widehat{K_{i_k}}\cdots K_{i_0}\Omega_\sigma=
K_{i_{p+1}}\cdots K_{i_0}(\dd f\wedge\Omega_\sigma)\\
-(-1)^{p+2}\dd f\wedge K_{i_{p+1}}\cdots K_{i_0}\Omega_\sigma.
\end{multline*} 
Applying equation (\ref{e:mul}) and $\dd f=0$ on $X$, we can see that
the image of $(\omega^i_A)_I$ under the \v{C}ech coboundary map is zero.
\end{pf}

Denote by $\widetilde\omega^i_A$ the image of the cocycle $(\widetilde\omega^i_A)_I$ in $H^p(X\cap D_i,\Omega^{d-3-p}_{X\cap D_i})$.
In the next step we show a relation  between $\widetilde\omega^i_A$ and $\omega^i_A$.
 
\begin{pr}\label{p:gysin}
Let $X\subset\ps$ be a $d$-semiample regular hypersurface defined by $f\in S_\beta$. Then
${\varphi_i}_!\widetilde\omega^i_{A}=\omega^i_{A}$, where  ${\varphi_i}_!$ is the Gysin map for  
$\varphi_i:X\cap D_i\hookrightarrow X$. 
\end{pr}

\begin{pf}
It suffices 
to show that ${\varphi_i}_!\widetilde\omega^i_{A}$, for $A\in S_{p\beta-\beta_0+\beta^\sigma_1}$, is represented by the 
\v{C}ech cocycle $(\omega^i_A)_I$.

The Gysin map ${\varphi_i}_!$ we can compute, using the following commutative diagram 
$$
\minCDarrowwidth0.4cm
\begin{CD} 
0 @>>> C^p({\cal V}^\sigma,\Omega^{d-1-p}_{X}) @>>> 
C^p({\cal V}^\sigma,\Omega_{X}^{d-1-p}(\log{X_i})) 
@>\res>> C^{p}({\cal V}_i^\sigma,\Omega^{d-2-p}_{X\cap D_i}) \\
@.    @AAA                               @AAA             @AAA   \\
0 @>>> C^{p-1}({\cal V}^\sigma,\Omega^{d-1-p}_{X}) @>>> 
C^{p-1}({\cal V}^\sigma,\Omega_{X}^{d-1-p}(\log{X_i})) 
@>\res>> C^{p-1}({\cal V}_i^\sigma,\Omega^{d-2-p}_{X\cap D_i}), 
\end{CD}
$$
where the vertical arrows are the \v{C}ech coboundary maps, 
${\cal V}^\sigma$  denotes the open cover
 ${\cal U}^\sigma|_X$,   and  the cover ${\cal V}_i^\sigma$  is the restriction 
${\cal V}^\sigma|_{X_i}$, $X_i=X\cap D_i$.
By the residue map, the cocycle $(\widetilde\omega^i_A)_{\tilde I}$ is lifted to the cochain 
$$
\psi_{\tilde I}=(-1)^{(p-1)^2/2}
\biggl\{\frac{AK_{\tilde i_{p-1}}\cdots K_{\tilde i_0}\Omega_\sigma}{f_{\tilde i_0}\cdots f_{\tilde i_{p-1}}}\wedge
\sum_{k=1}^n\langle m_{\tilde j_0},e_k\rangle \frac{\dd x_k}{x_k}\biggr\}
_I$$
in $C^{p-1}({\cal V}^\sigma,\Omega_{X}^{d-1-p}(\log{X_i}))$, where $\tilde I$ is the  index set 
$\{(\tilde i_0,{\tilde j_0}),\dots,(\tilde i_{p-1},{\tilde j_{p-1}})\}$, corresponding to the cover 
${\cal V}^\sigma$, and where $m_{\tilde j_0}\in M_{\Bbb R}$,
for $\sigma_{\tilde j_0}\supset\rho_i$ generated by $e_i$ and $e_s$, satisfies $\langle m_{\tilde j_0},e_i\rangle=1$,
$\langle m_{\tilde j_0},e_s\rangle=0$, and $m_{\tilde j_0}=0$ in all other cases.
Appropriately, this can be obtained, using some affine open cover on $X$, where $X\cap D_i$ is given by 
$\prod_{k=1}^n x_k^{\langle m,e_k\rangle}=0$ up to some multiplicity (we omit the details).   

The image of $\psi_{\tilde I}$ under the \v{C}ech coboundary map should represent ${\varphi_i}_!\widetilde\omega^i_{A}$.
Using the diagram, we can see that changing of  $\psi_{\tilde I}$ by a cochain in 
$C^{p-1}({\cal V}^\sigma,\Omega^{d-1-p}_{X})$ does not affect the image.
Notice that $\psi_{\tilde I}$ is equivalent to
$$
(-1)^{(p-1)^2/2}\biggl\{\frac{AK_{\tilde i_{p-1}}\cdots K_{\tilde i_0}}{f_{\tilde i_0}\cdots f_{\tilde i_{p-1}}}
\biggl(\Omega_\sigma\wedge\sum_{k=1}^n\langle m_{\tilde j_0},e_k\rangle \frac{\dd x_k}{x_k}\biggr)\biggr\}_{\tilde I}
$$ 
modulo some cochain in $C^{p-1}({\cal V}^\sigma,\Omega^{d-1-p}_{X})$.
Assume for a moment that 
\begin{equation}\label{e:w}
\biggl(\Omega_\sigma\wedge\sum_{k=1}^n\langle m_{\tilde j_0},e_k\rangle\frac{\dd x_k}{x_k}\biggr)
-(-1)^{d}\frac{\partial^i_{\tilde j_{0}}\lrcorner\Omega}{\prod_{\rho_k\subset\sigma}x_k}
\end{equation} 
is well defined on $U_{\sigma_{\tilde j_0}}$. Then 
$\psi_{\tilde I}$ is actually equivalent to
$$
(-1)^{d+((p-1)^2/2)}\biggl\{\frac{AK_{\tilde i_{p-1}}\cdots K_{\tilde i_0}}{f_{\tilde i_0}\cdots f_{\tilde i_{p-1}}}
\biggl(\frac{\partial^i_{\tilde j_{0}}\lrcorner\Omega}{\prod_{\rho_k\subset\sigma}x_k}\biggr)\biggr\}_{\tilde I}
$$ 
modulo some cochain in $C^{p-1}({\cal V}^\sigma,\Omega^{d-1-p}_{X})$. The image
of this  under the \v{C}ech coboundary map is clearly $(\omega^i_A)_I$. 
We are left to show that (\ref{e:w})  is well defined on $U_{\sigma_{\tilde j_0}}$.
The case $\sigma_{\tilde j_0}\not\supset\rho_i$ is trivial because $m_{\tilde j_0}=0$ and $\partial^i_{\tilde j_{0}}=0$.
The cases left are $\tilde j_0=k,k+1$ for $i=l_k$ as in $(\ref{e:pic})$; we only check the case $\tilde j_0=k$ 
(then $\langle m_{\tilde j_0},e_{l_{k-1}}\rangle=0$, 
$\partial^i_{\tilde j_{0}}=x_{l_{k-1}}\partial_{l_{k-1}}/(\mult(\sigma_k)$), 
the other case is similar.
It is enough to verify that multiples of $(\dd x_i)/x_i$ cancel each other in the difference (\ref{e:w}).
We defined  $\Omega=\sum_{|I|=d}\det(e_I)\hat{x}_I \dd x_I$; note that the multiples of $(\dd x_i)/x_i$ 
in (\ref{e:w}) are
$$
\Biggl(\sum_{|J|=d-2}
\biggl(\frac{\det(e_{\{l_0,l_{n(\sigma)+1}\}\cup J})}{\mult(\sigma)}-(-1)^{d}
\frac{\det(e_{\{l_{k-1}\}\cup J\cup\{i\}})}{\mult(\sigma_k)}
\biggr)\frac{\hat{x}_J \dd x_J}{\prod_{\rho_k\subset\sigma}x_k}\Biggr)\wedge\frac{\dd x_i}{x_i},
$$
where the sum is over all $(d-2)$-element subsets $J\subset\{1,\dots,n\}.$
Interchanging $i=l_k$ with the ordered set $\{l_{k-1}\}\cup J$ in $\det(e_{\{l_{k-1}\}\cup J\cup\{i\}})$ and using 
the relations of the cone generators 
(see \cite[Section~8.2]{d})
$$
\frac{e_{l_{k}}}{\mult(\sigma_k)}=
\frac{\mult(\sigma_{0,k})e_{l_{k-1}}}{\mult(\sigma_k)\mult(\sigma_{0,k-1})}-
\frac{e_{l_{0}}}{\mult(\sigma_{0,k-1})},
$$
$$
\frac{\mult(\sigma)e_{l_{k-1}}}{\mult(\sigma_{0,k-1})\mult(\sigma_{k-1,n(\sigma)+1})}=
\frac{e_{l_{n(\sigma)+1}}}{\mult(\sigma_{k-1,n(\sigma)+1})}+
\frac{e_{l_{0}}}{\mult(\sigma_{0,k-1})},
$$
where $\sigma_{s,t}$  denotes the cone generated by $e_{l_s}$ and $e_{l_t}$, we get that the multiples of 
$(\dd x_i)/x_i$ in (\ref{e:w}) cancel each other. The proposition is proved.
\end{pf}

The last proposition shows  the relation of  $\omega_A^i$ to the description of
the middle cohomology of $X$ given in equation (\ref{e:mid}). But we also need to understand 
the relation of $\omega_A^i$ to the description of the  cohomology in Theorem~\ref{t:mid}. 
For this, we will have to consider some toric subvarieties of codimension 2 in $\ps$, and to study
the relation of some quotients of the homogeneous coordinate rings of these toric subvarieties and $\ps$.
This  work will culminate in Theorem~\ref{t:midn}, which generalizes Theorem~\ref{t:mid}.

As in \cite[Section~5]{m}, we consider a 2-dimensional cone $\sigma'\in\Sigma$ contained in $\sigma\in\Sigma_X(2)$ and
containing $\rho_i$ (in the notation of (\ref{e:pic}), we have $i=l_k$ and $\sigma'=\sigma_k$ or $\sigma_{k+1}$),
 and let $S(V(\sigma'))={\Bbb C}[x_{\gamma'}:\sigma'\subset{\gamma'}\in\Sigma(3)]$  
be  the coordinate ring of the $(d-2)$-dimensional 
complete simplicial toric variety $V(\sigma')\subset\ps$.
{From} Lemma 1.4 in \cite{m}, it follows that $X_{\sigma'}:=X\cap V(\sigma')$ (we will use both notations)
has a positive self-intersection number inside 
$V(\sigma')$,
implying $X_{\sigma'}$ is a big and nef hypersurface.
We have a natural commutative diagram:
\begin{equation*}
\begin{CD}
S_{*\beta} @.\cong @. H^0(\ps,O_\ps(*X))\\
@V{\scriptstyle{\scriptstyle\varphi_{\sigma'}^*}}VV @. @V{\scriptstyle{\scriptstyle\varphi_{\sigma'}^*}}VV  \\
S(V(\sigma'))_{*\beta^{\sigma'}} @. \cong @. H^0(V(\sigma'),O_{V(\sigma')}(*X_{\sigma'})),
\end{CD}
\end{equation*} 
where  
$\beta^{\sigma'}\in A_{d-3}(V(\sigma'))$ is the restriction of $\beta$, 
and the vertical arrows are the restriction maps induced by the inclusion $\varphi_{\sigma'}:V(\sigma')\subset\ps$.
To describe the vertical arrow on the left one first has to restrict a Cartier divisor 
$D=\sum_{k=1}^n a_k D_k$ (as in \cite[Section~5.1]{f1}, assuming that $a_k=0$ for $\rho_k\subset\sigma'$) in
degree $\beta$  to $V(\sigma')$:
 $$D|_{V(\sigma')}=\sum_{\gamma'}a_{k(\gamma')}\frac{\mult(\sigma')}{\mult(\gamma')}V(\gamma'),$$
where the sum is over all $\gamma'\in\Sigma(3)$ spanned by $\sigma'$ and a generator $e_{k(\gamma')}$.
Then a monomial $\prod_{k=1}^n x_k^{qa_k+\langle m,e_k\rangle}$
 in $S_{q\beta}$ with $m\in{\sigma'}^\perp$
is sent by the restriction map $\varphi_{\sigma'}^*$ to 
$\prod_{\gamma'} x_{\gamma'}^{qa_{\gamma'}+\langle m,e_{\gamma'}\rangle}$, where 
$a_{\gamma'}=a_{k(\gamma')}{\mult(\sigma')}/{\mult(\gamma')}$ and 
$e_{\gamma'}=e_{k(\gamma')}{\mult(\sigma')}/{\mult(\gamma')}$; if $m\notin{\sigma'}^\perp$, the monomial is 
sent to 0. Hence, 
we can see that the restriction map $S_{*\beta}@>>>S(V(\sigma'))_{*\beta^{\sigma'}}$ is surjective, and
its kernel is the ideal in $S_{*\beta}$ generated by all variables $x_k$ such that $\rho_k\subset\sigma$, 
by the argument in the proof of Lemma~\ref{l:non}.
Therefore, we have an isomorphism:
\begin{equation}\label{e:isor}
\varphi_{\sigma'}^*:(S/\langle x_k:\rho_k\subset\sigma\rangle)_{*\beta}\cong S(V(\sigma'))_{*\beta^{\sigma'}}.
\end{equation} 
If $X$ is defined by $f\in S_\beta$, then the restriction  of $f$, denoted by $f_{\sigma'}$, determines exactly the hypersurface  
$X_{\sigma'}\subset V(\sigma')$.  

We also have a natural map
\begin{equation}\label{e:resn}
S_{(p+1)\beta-\beta_0+\beta_1^\sigma}@>>>H^0(D_i,\Omega^{d-2}_{D_i}((p+1)X_i)),
\end{equation} 
sending $A$ to the rational $(d-2)$-form $(A\Omega_\sigma/f^{p+1})$ considered after Definition~\ref{d:gysc}.
Let us determine the restriction of this form with respect to the map 
$$H^0(D_i,\Omega^{d-2}_{D_i}((p+1)X_i))@>\varphi_{i,\sigma'}^*>>H^0(V(\sigma'),\Omega^{d-2}_{V(\sigma')}((p+1)X_{\sigma'})),$$
induced by the inclusion $\varphi_{i,\sigma'}:V(\sigma')\subset D_i$.
The form $\Omega$ in Definition~\ref{d:om} is determined up to $\pm1$, depending on the choice of the basis for the lattice $M$.
We have fixed one basis $m_1,\dots,m_d$, but it is always possible to find another basis $m_1^\sigma,\dots,m^\sigma_d$,
for $\sigma\in\Sigma_X(2)$,
so that the corresponding $\Omega$ is the same as before and 
$m_1^\sigma,\dots,m^\sigma_{d-2}$ form a basis for the lattice $M\cap\sigma^\perp$.
With the new choice of the basis,  the proof of Proposition 9.5 in \cite{bc} shows that 
$$\Omega=\prod_{k=1}^nx_k\biggl(\sum_{k=1}^n\langle m^\sigma_1,e_k\rangle\frac{\dd x_k}{x_k}\biggr)\wedge\cdots\wedge 
\biggl(\sum_{k=1}^n\langle m^\sigma_d,e_k\rangle\frac{\dd x_k}{x_k}\biggr).$$
Using this, we compute 
\begin{multline*}
{\Omega_\sigma}=
\frac{x_{l_0}x_{l_{n(\sigma)+1}}K_{l_{n(\sigma)+1}}K_{l_0}\Omega}
{\mult(\sigma)\prod_{\rho_k\subset\sigma}x_k}\\
=\frac{\prod_{\rho_k\not\subset\sigma}x_k}{\mult(\sigma)}e^{d-1,d}
\biggl(\sum_{k=1}^n\langle m^\sigma_1,e_k\rangle\frac{\dd x_k}{x_k}\biggr)
\wedge\cdots\wedge 
\biggl(\sum_{k=1}^n\langle m^\sigma_{d-2},e_k\rangle\frac{\dd x_k}{x_k}\biggr),
\end{multline*}
where 
$e^{d-1,d}$ denotes $(\langle m^\sigma_{d-1},e_{l_0}\rangle\langle m^\sigma_{d},e_{l_{n(\sigma)+1}}\rangle-
\langle m^\sigma_{d},e_{l_{0}}\rangle\langle m^\sigma_{d-1},e_{l_{n(\sigma)+1}}\rangle)$. By the properties of 
$\mult(\sigma)$ in \cite[Section~8]{d}, we can see that $e^{d-1,d}/\mult(\sigma)$ is $\pm1$. There were two
(reverse to each other) possibilities 
of labeling the generators of $\sigma$ when we chose the order in (\ref{e:pic}). {\it In further calculations we assume}
such a choice of $\rho_{l_0}$ and  $\rho_{l_{n(\sigma)+1}}$  that $e^{d-1,d}/\mult(\sigma)=1$.
Set $t^\sigma_j=\prod_{k=1}^nx_k^{\langle m^\sigma_j,e_k\rangle}$, then $t^\sigma_1,\dots,t^\sigma_{d-2}$ are 
the coordinates on the torus $\ttt_{\sigma'}$. In terms of the homogeneous coordinates $x_{\gamma'}$ on $V(\sigma')$,
the affine coordinates $t^\sigma_j$ are  identified with $\prod_{\gamma'}x_{\gamma'}^{\langle m^\sigma_j,e_{\gamma'}\rangle}$.
Hence, 
\begin{multline*}
\varphi_{i,\sigma'}^*\biggl(\frac{A\Omega_\sigma}{f^{p+1}}\biggr)=
\varphi_{i,\sigma'}^*\biggl(\frac{A\prod_{\rho_k\not\subset\sigma}x_k}{f^{p+1}}
\frac{\dd t^\sigma_1}{t^\sigma_1}\wedge\cdots\wedge\frac{\dd t^\sigma_{d-2}}{t^\sigma_{d-2}}\biggr)
\\=
\frac{\varphi_{\sigma'}^*(A\prod_{\rho_k\not\subset\sigma}x_k)}{f_{\sigma'}^{p+1}}
\biggl(\sum_{\gamma'}\langle m^\sigma_1,e_{\gamma'}\rangle\frac{\dd x_{\gamma'}}{x_{\gamma'}}\biggr)
\wedge\cdots\wedge 
\biggl(\sum_{\gamma'}\langle m^\sigma_{d-2},e_{\gamma'}\rangle\frac{\dd x_{\gamma'}}{x_{\gamma'}}\biggr)
\\
=\frac{\varphi_{\sigma'}^*(A\prod_{\rho_k\not\subset\sigma}x_k)}{(\prod_{\gamma'}x_{\gamma'})f_{\sigma'}^{p+1}}\Omega_{V(\sigma')},
\end{multline*}
where, as in Definition~\ref{d:om}, 
$\Omega_{V(\sigma')}$ is the $(d-2)$-form on the toric variety $V(\sigma')$, 
corresponding to the basis  $m_1^\sigma,\dots,m^\sigma_{d-2}$.
A monomial in  $S_{(p+1)\beta-\beta_0+\beta_1^\sigma}$ with $\beta=[\sum_{k=1}^n a_kD_k]$  corresponds to 
a lattice point $m$, satisfying the inequalities  $(p+1)a_k+\langle m,e_{k}\rangle\ge0$, for $\rho_k\subset\sigma$,
and $(p+1)a_k+\langle m,e_{k}\rangle\ge1$, for $\rho_k\not\subset\sigma$.
Then, by the earlier explicit description of $\varphi_{\sigma'}^*$, we can see that
the restriction $\varphi_{\sigma'}^*(A\prod_{\rho_k\not\subset\sigma}x_k)$ is a polynomial 
in $S(V(\sigma'))_{(p+1)\beta^{\sigma'}}$, divisible by $\prod_{\gamma'}x_{\gamma'}$.
Therefore, we get the following commutative diagram
$$
\begin{CD}
S_{(p+1)\beta-\beta_0+\beta_1^\sigma}@>>>H^0(D_i,\Omega^{d-2}_{D_i}((p+1)X_i)\\
@V{\scriptstyle{\scriptstyle\varphi_{\sigma'}^*}}VV  @V{\scriptstyle{\scriptstyle\varphi_{i,\sigma'}^*}}VV  \\
S(V(\sigma'))_{(p+1)\beta^{\sigma'}-\beta_0^{\sigma'}}@>>>H^0(V(\sigma'),\Omega^{d-2}_{V(\sigma')}((p+1)X_{\sigma'})),
\end{CD}
$$
where $\beta_0^{\sigma'}:=\deg(\prod_{\gamma'}x_{\gamma'})\in A_{d-3}(V(\sigma'))$ is the anticanonical degree, and 
the horizontal arrows are given  by (\ref{e:resn}) and a similar one sending a polynomial $A$ to the form
$(A\Omega_{V(\sigma')}/f_{\sigma'}^{p+1})$. 
Recall  {from} Section~\ref{s:pr} that for the hypersurface $X_{\sigma'}\subset V(\sigma')$ we have the residue map 
$$\res: S(V(\sigma'))_{(p+1)\beta^{\sigma'}-\beta_0^{\sigma'}}@>>>H^p(X,\Omega_{X\cap V({\sigma'})}^{d-3-p}),$$  
sending a polynomial $B$ to the Hodge component $\res(\omega_B)^{d-3-p,p}$. As in Section~\ref{s:p},  denote
$[\omega_B]=(-1)^{p/2}p!\res(\omega_B)^{d-3-p,p}$.
By  the naturality of the residue map and Proposition~\ref{p:rel1}, we obtain the following result.

\begin{pr}\label{p:relres} Let $X\subset\ps$ be a $d$-semiample regular hypersurface 
defined by $f\in S_\beta$. Given $\rho_i\subset\sigma\in\Sigma_X(2)$ such that 
$\rho_i\notin\Sigma_X(1)$, and given $\sigma'\in\Sigma(2)$ such that $\rho_i\subset\sigma'\subset\sigma$,
then we have a commutative diagram:
$$
\minCDarrowwidth2cm
\begin{CD}
S_{(p+1)\beta-\beta_0+\beta_1^\sigma}@>\tilde\omega^i_{\_}>>H^p(X\cap D_i,\Omega^{d-3-p}_{X\cap D_i})\\
@V{\scriptstyle{\scriptstyle\varphi_{\sigma'}^*}}VV  @V{\scriptstyle{\scriptstyle\varphi_{i,\sigma'}^*}}VV  \\
S(V(\sigma'))_{(p+1)\beta^{\sigma'}-\beta_0^{\sigma'}}@>(-1)^{d-3-p}[\omega_{\_}]>>
H^p(X\cap V(\sigma'),\Omega^{d-3-p}_{X\cap V(\sigma')}).
\end{CD}
$$
\end{pr}

{From} Section~\ref{s:pr} we know that the map
$$\res: R_1(f_{\sigma'})_{(p+1)\beta^{\sigma'}-\beta_0^{\sigma'}}@>>>H^p(X\cap V({\sigma'}),\Omega_{X\cap V({\sigma'})}^{d-3-p})$$  
 is  well defined. The map $\tilde\omega^i_{\_}$ should also be well defined on some quotient of the coordinate ring $S$. 
In Definition~\ref{d:r1} we had the rings
 $R_0(f)=S/J_0(f)$ and $R_1(f)=S/J_1(f)$.
Now introduce the following similar rings.

\begin{defn}\label{d:rs1}
Given  $f\in S_\beta$ of $d$-semiample degree $\beta\in A_{d-1}(\ps)$   and $\sigma\in\Sigma_\beta(2)$
(see Remark~\ref{r:can}), 
let
$J^\sigma_0(f)$ be the ideal  in $S$  
generated by the ideal $J_0(f)$ and all $x_k$ such that $\rho_k\subset\sigma$, and
let $J^\sigma_1(f)$ be the ideal quotient $J^\sigma_0(f):(\prod_{\rho_k\not\subset\sigma}x_k)$. 
Then we get the quotient rings $R_0^\sigma(f)=S/J_0^\sigma(f)$ and  $R_1^\sigma(f)=S/J_1^\sigma(f)$
graded by the Chow group $A_{d-1}(\ps)$.
\end{defn}

We have the toric morphism $\pi:\ps@>>>\psx$, associated with a $d$-semiample hypersurface $X\subset \ps$.
By the previous discussion, for $\sigma'\subset\sigma\in\Sigma_X(2)$,
$X_{\sigma'}=X\cap V(\sigma')$ is a  big and nef hypersurface, defined by $f_{\sigma'}$, in the 
toric variety $V(\sigma')$. It follows {from}  Proposition~\ref{p:comp1} that
the restriction of $\pi$ is the toric morphism $\pi_{\sigma'}: V(\sigma')@>>>V(\sigma)$, associated with 
the semiample divisor $X_{\sigma'}$. In particular, we have a ring homomorphism ${\pi_{\sigma'}}_*:S(V(\sigma'))@>>>S(V(\sigma))$
between the coordinate rings of the toric varieties. The image of $f_{\sigma'}$ is a polynomial 
$f_{\sigma}\in S(V(\sigma))_{\beta^\sigma}$, which determines
the ample hypersurface $Y_{\sigma}:=\pi_{\sigma'}(X_{\sigma'})$ in $V(\sigma)$.

\begin{pr}\label{p:2iso} Let $\beta\in A_{d-1}(\ps)$ be $d$-semiample 
and let $\beta_0^{\sigma'}=\deg(\prod_{\gamma'}x_{\gamma'})\in A_{d-3}(V(\sigma'))$,
$\beta_0^{\sigma}=\deg(\prod_{\gamma}y_{\gamma})\in A_{d-3}(V(\sigma))$ be the anticanonical degrees. 
Then, there are natural isomorphisms:
 
{\rm (i)} $R^\sigma_0(f)_{*\beta}\,{\cong}\, R_0(f_{\sigma'})_{*\beta^{\sigma'}}
\,{\cong}\,R_0(f_{\sigma})_{*\beta^{\sigma}}$,

{\rm (ii)} $R^\sigma_1(f)_{*\beta-\beta_0+\beta_1^\sigma}\,{\cong}\,R_1(f_{\sigma'})_{*\beta^{\sigma'}-\beta_0^{\sigma'}}
\,{\cong}\,R_1(f_{\sigma})_{*\beta^{\sigma}-\beta_0^{\sigma}}$.
\end{pr}

\begin{pf} (i) To show the first isomorphism, induced by $\varphi_{\sigma'}^*$, it suffices, 
because of equation (\ref{e:isor}), 
to check that the ideal $J_0(f)$ in $S$ is mapped onto the ideal  $J_0(f_{\sigma'})$ in $S(V(\sigma'))$.
By Proposition 5.3 in \cite{c2}, the ideal  $J_0(f)$ is generated by $f$ and 
$x_{i_1}\partial f/\partial x_{i_1},\dots,x_{i_d}\partial f/\partial x_{i_d}$ for linearly independent $e_{i_1},\dots e_{i_d}$. 
We can assume that $e_{i_1},\dots e_{i_d}$ are generators of some simplicial cone $\tau$, containing $\sigma'$, and
$e_{i_{d-1}}, e_{i_d}$ are generators of $\sigma'$. By the explicit description
of the restriction map $\varphi_{\sigma'}^*$, $f$ is sent to $f_{\sigma'}$, while $x_{i_{d-1}}\partial f/\partial x_{i_{d-1}}$
and $x_{i_{d}}\partial f/\partial x_{i_d}$ are sent to $0$. To understand the image of the other polynomials,
as in \cite{bc}, we
write $f=\sum_{m\in\Delta\cap M}a_m\prod_{k=1}^n x_k^{b_k+\langle m,e_{k}\rangle}$, where $\Delta$ is the polytope associated with
a torus invariant divisor $\sum_{k=1}^nb_kD_k$ (assuming $b_{i_{d-1}}=b_{i_d}=0$) in degree $\beta$.
Then 
$$x_{i_s}\frac{\partial f}{\partial x_{i_s}}=\sum_{m\in\Delta\cap M}a_m(b_{i_s}+\langle m,e_{i_s}\rangle)
\prod_{k=1}^nx_k^{b_k+\langle m,e_{k}\rangle}.$$
Applying the restriction map $\varphi_{\sigma'}^*$ to this,
we get, for $s\ne d-1,d$,
$$
\sum_{m\in\Delta\cap M\cap\sigma^\perp}a_m(b_{i_s}+\langle m,e_{i_s}\rangle)
\prod_{\gamma'}x_{\gamma'}^{b_{\gamma'}+\langle m,e_{\gamma'}\rangle}=
\frac{\mult(\gamma'_s)}{\mult(\sigma')}x_{\gamma'_s}\frac{\partial f}{\partial x_{\gamma'_s}},
$$
where the cone $\gamma'_s$ is spanned by $\sigma'$ and the generator
$e_{i_s}$, and where $b_{\gamma'}=b_{k(\gamma')}\mult(\sigma')/\mult(\gamma')$,
$e_{\gamma'}=e_{k(\gamma')}\mult(\sigma')/\mult(\gamma')$ correspond to the cone ${\gamma'}$ spanned by $\sigma'$ and a generator
$e_{k(\gamma')}$. Therefore, we get the first isomorphism  
$\varphi_{\sigma'}^*:R^\sigma_0(f)_{*\beta}{\cong} R_0(f_{\sigma'})_{*\beta^{\sigma'}}$.

For the second isomorphism, induced by 
${\pi_{\sigma'}}_*:S(V(\sigma'))_{*\beta^{\sigma'}}\cong S(V(\sigma))_{*\beta^{\sigma}}$ (see Section~\ref{s:pr}), 
it is enough to show 
that $J_0(f_{\sigma'})$ is mapped onto the ideal $J_0(f_{\sigma})$ in $S(V(\sigma))$.
This can be easily achieved by the argument in the previous paragraph.

(ii) 
By the construction of the maps $\varphi_{\sigma'}^*$ and ${\pi_{\sigma'}}_*$, we get the 
commutative diagram:
$$
\begin{CD}
R^\sigma_0(f)_{*\beta}@.{\cong}@. R_0(f_{\sigma'})_{*\beta^{\sigma'}}@.{\cong}@.R_0(f_{\sigma})_{*\beta^{\sigma}}\\
@AA\prod_{\rho_k\not\subset\sigma}x_kA @. @AA\prod_{\gamma'}x_{\gamma'}A @. @AA\prod_{\gamma}y_{\gamma}A\\
R^\sigma_1(f)_{*\beta-\beta_0+\beta_1^\sigma}@.\rightarrow @.R_1(f_{\sigma'})_{*\beta^{\sigma'}-\beta_0^{\sigma'}}
@.\rightarrow @.R_1(f_{\sigma})_{*\beta^{\sigma}-\beta_0^{\sigma}},
\end{CD}
$$
where the vertical arrows are injections, induced by the multiplication.
To show that the bottom arrows   are isomorphisms it suffices to check that the images
of the spaces {from} the bottom into the spaces on the top correspond to each other under the isomorphisms of part (i).
Note that these images are the ideals generated by $\prod_{\rho_k\not\subset\sigma}x_k$, $\prod_{\gamma'}x_{\gamma'}$ and
$\prod_{\gamma}y_{\gamma}$, respectively. By the explicit description of the maps $\varphi_{\sigma'}^*$ and ${\pi_{\sigma'}}_*$,
one can see that these are mapped onto each other.
\end{pf}

Finally, we can put all of the above together and generalize equation (\ref{e:mid}),
 describing  the middle cohomology of a big and nef regular hypersurface.

\begin{thm}\label{t:midn} 
Let $X\subset\ps$ be  a $d$-semiample regular hypersurface
defined by $f\in S_\beta$, $d=\dim\ps$.  
Then there is a natural isomorphism, for $p=d-1-q$:
$$H^{p,q}(X)\cong R_1(f)_{(q+1)\beta-\beta_0}
\oplus\Biggl(\bigoplus\begin{Sb}\sigma
\in\Sigma_X(2)\end{Sb}(R^\sigma_1(f)_{q\beta-\beta_0+\beta_1^\sigma})
^{n(\sigma)}\Biggr) \oplus H^{p,q}_{\rm toric}(X)\oplus C,$$
where  $C=\sum_{\tau\in\Sigma(2)} {\varphi_\tau}_!H_{\rm res}^{p-2,q-2}(X\cap V(\tau))$ (the Gysin maps ${\varphi_\tau}_!$
are induced by the inclusions $\varphi_\tau:X\cap V(\tau)\subset X$), and
the graded pieces of $R_1(f)$ and $R^\sigma_1(f)$ are embedded by the maps
$[\omega_{\_}]$ and  $\omega^i_{\_}$ for all $\rho_i\notin\Sigma_X$ contained in some $\sigma\in\Sigma_X(2)$
 ($n(\sigma)$ is the number of such cones).
Moreover, $R^\sigma_1(f)_{q\beta-\beta_0+\beta_1^\sigma}=0$ for $q=0,d-1$, and the cup product of any two elements {from}
the distinct summands of the above decomposition vanishes.
\end{thm}

\begin{pf} Theorem~4.4 in \cite{m}  combined with the diagram (\ref{e:resi}) gives an isomorphism:
\begin{equation*}
H^{d-1-q,q}(X)\cong R_1(f)_{(q+1)\beta-\beta_0}\oplus H^{d-1-q,q}_{\rm toric}(X)\oplus\sum_{i=1}^n {\varphi_i}_! 
H_{\rm res}^{d-2-q,q-1}(X\cap D_i),
\end{equation*}
where  ${\varphi_i}_!$ are the Gysin maps induced by the inclusions.
Applying (\ref{e:resi}) to the hypersurface $X\cap D_i$ in $D_i$, we get an exact sequence
\begin{equation}\label{e:exacts}
\bigoplus_{\rho_i\subset\tau\in\Sigma(2)} H_{\rm res}^{d-5}(X\cap V(\tau))@>>>PH^{d-3}(X\cap D_i)
@>>>\gr_{d-3}^W PH^{d-3}(X\cap\ttt_{\rho_i}).
\end{equation}
The space $\gr_{d-3}^W PH^{d-3}(X\cap\ttt_{\rho_i})$ vanishes, by equation (\ref{e:grvan}), unless 
$X\cap D_i$ is a $(d-2)$-semiample hypersurface in $D_i$. By Proposition~\ref{p:comp}, the latter happens only
 when  $\rho_i\notin\Sigma_X$  lies in some $\sigma\in\Sigma_X(2)$. 
In  this case, there is $\sigma'\in\Sigma(2)$ such that $\rho_i\subset\sigma'\subset\sigma$, and,
by equation (\ref{e:griso}), 
 we have isomorphisms
$$\gr_{d-3}^W H^{d-3}(X\cap\ttt_{\rho_i})\cong \gr_{d-3}^W H^{d-3}(\pi(X)\cap\ttt_{\sigma})\cong 
\gr_{d-3}^W H^{d-3}(X\cap\ttt_{\sigma'})$$
induced by the morphism $\pi:\ps@>>>\psx$.
The hypersurface $X\cap V(\sigma')$ in $V(\sigma')$ is $(d-2)$-semiample (big and nef). So, we can apply Theorem~4.4 in \cite{m}  
 to deduce that 
the composition
$$R_1(f_{\sigma'})_{q\beta^{\sigma'}-\beta_0^{\sigma'}}@>\res>>
H_{\rm res}^{d-2-q,q-1}(X\cap V({\sigma'}))@>>>H^{d-2-q,q-1}\bigl(PH^{d-3}(X\cap\ttt_{\sigma'})\bigr)$$  
is an isomorphism.
Using  Propositions~\ref{p:relres}~and~\ref{p:2iso}, we get that another composition
$$R^\sigma_1(f)_{q\beta-\beta_0+\beta_1^\sigma}@>\tilde\omega^i_{\_}>>H_{\rm res}^{d-2-q,q-1}(X\cap D_i)@>>>
H^{d-2-q,q-1}\bigl(PH^{d-3}(X\cap\ttt_{\rho_i})\bigr)$$
is also an isomorphism. Hence, by equation (\ref{e:exacts}), 
$$H_{\rm res}^{d-2-q,q-1}(X\cap D_i)\cong R^\sigma_1(f)_{q\beta-\beta_0+\beta_1^\sigma}\oplus
\sum_{\rho_i\subset\tau\in\Sigma(2)}{\varphi^i_\tau}_!H_{\rm res}^{d-3-q,q-2}(X\cap V(\tau))$$
for $\rho_i\notin\Sigma_X$   contained in some $\sigma\in\Sigma_X(2)$, and
$$H_{\rm res}^{d-2-q,q-1}(X\cap D_i)\cong 
\sum_{\rho_i\subset\tau\in\Sigma(2)}{\varphi^i_\tau}_!H_{\rm res}^{d-3-q,q-2}(X\cap V(\tau))$$
for all other $\rho_i$ (here, $\varphi^i_\tau:X\cap V(\tau)\subset X\cap D_i$ is the inclusion).
{From} (\ref{e:resi}) we have an exact sequence
$$\bigoplus_{\tau\in\Sigma(2)} H_{\rm res}^{d-5}(X\cap V(\tau))@>>>\bigoplus_{i=1}^n H_{\rm res}^{d-3}(X\cap D_i)
@>>>H_{\rm res}^{d-1}(X)$$
which shows that the kernel of the right arrow is included into the parts complementary to 
$R^\sigma_1(f)_{q\beta-\beta_0+\beta_1^\sigma}$ in $H_{\rm res}^{d-3}(X\cap D_i)$.
The  direct sum decomposition of the middle cohomology follows. 

The fact $R^\sigma_1(f)_{-\beta_0+\beta_1^\sigma}=0$ is obvious, while $R^\sigma_1(f)_{(d-1)\beta-\beta_0+\beta_1^\sigma}=0$
is implied by the isomorphism of Proposition~\ref{p:2iso} and by a dimension argument using the proof of Theorem~11.5 in
\cite{bc}  and Theorems~2.11, 4.8(v) with Corollary~3.14 in \cite{b2}. 
{From} Section~\ref{s:tr} we know that $H^*_{\rm res}(X)\cup H^{*}_{\rm toric}(X)\subset H^*_{\rm res}(X)$. But since
$H^{2d-2}_{\rm res}(X)=0$, the toric part $H^{d-1}_{\rm toric}(X)$ is orthogonal to all other summands  in the middle  cohomology.
Theorem 4.4 in \cite{m} shows that  $R_1(f)_{(q+1)\beta-\beta_0}$ is  orthogonal  to all other summands as well.
The proof of Lemma~\ref{l:van2} below shows that $\omega^i_{A}\cup\omega^j_{B}=0$ 
if $\rho_i,\rho_j\notin\Sigma_X$ lie in two distinct $2$-dimensional cones of $\Sigma_X(2)$. Finally, the projection formula
gives:
$$\omega^i_{A}\cup{\varphi_\tau}_!H_{\rm res}^{d-5}(X\cap V(\tau))={\varphi_\tau}_!(\varphi_\tau^*\omega^i_{A}\cup
H_{\rm res}^{d-5}(X\cap V(\tau)).$$
But it can be seen directly that the restriction $\varphi_\tau^*$ of the \v Cech cocycle $(\omega^i_{A})_I$ is a coboundary.
The theorem is proved.
\end{pf}

\begin{rem} The direct summand $C$ in the above theorem vanishes when $q=0,1,d-1,d-2$. Therefore, we have a complete description
of the middle cohomology in the corresponding Hodge degrees.
\end{rem}

Lemma~\ref{l:van1} tells us that the cup product $\gamma^i_A\cup\omega^j_B$ vanishes in certain cases.
Now we show that this is true in more cases.

\begin{lem}\label{l:van2} Let $X\subset\ps$ be  a  $d$-semiample regular hypersurface
defined by $f\in S_\beta$. 
Then the cup product $\gamma^i_{A}\cup\omega^j_{B}=0$ if $\rho_i,\rho_j\notin\Sigma_X(1)$ lie in two
distinct $2$-dimensional cones of $\Sigma_X(2)$. 
\end{lem}

\begin{pf} We  use the description of the middle cohomology in equation (\ref{e:mid}) and the Poincar\'e nondegenerate pairing
to show that $\gamma^i_{A}\cup\omega^j_{B}=0$ for $\rho_i$ and $\rho_j$ lying in two distinct 2-dimensional cones 
$\sigma^1$ and $\sigma^2$ of $\Sigma_X$. 
Because of this, it is enough to check
 that the cup product of $\gamma^i_{A}\cup\omega^j_{B}$ with all elements in (\ref{e:mid}) vanishes.

Take $[\omega_C]\in H^{d-1}(X)$, corresponding to $C\in S_{({d-p-1})\beta-\beta_0}$,  
in the Hodge component complementary to the one of $\gamma^i_{A}\cup\omega^j_{B}$.
Then
$$\gamma^i_A\cup\omega^j_{B}\cup[\omega_C]=\pm\omega^i_{AC}\cup\omega^j_{B}=
\pm{\varphi_i}_!\widetilde\omega^i_{AC}\cup{\varphi_j}_!\widetilde\omega^j_{B}=
\pm{\varphi_j}_!((\varphi_j^*{\varphi_i}_!\widetilde\omega^i_{AC})\cup\widetilde\omega^j_{B}),
$$
where we use Theorem~\ref{t:nonm}, Proposition~\ref{p:gysin} and the projection formula for Gysin homomorphisms.
By Lemma 5.4 in \cite{m},
there is a commutative diagram:
\begin{equation*}
\minCDarrowwidth1.1cm
\begin{CD}
H^{d-3}(X\cap D_i)@>{\varphi_i}_!>>H^{d-1}(X)\\
 @V\varphi_{ij}^*VV @V\varphi_j^*VV \\
H^{d-3}(X\cap D_i\cap D_j)@>\alpha\cdot{\varphi_{ji}}_!>> H^{d-1}(X\cap D_j), 
\end{CD}
\end{equation*}
where $\varphi_{ij}:X\cap D_i\cap D_j\hookrightarrow X\cap D_i$ is the inclusion map and $\alpha$
is some constant.
On the other hand, $\varphi_{ij}^*\widetilde\omega^i_{AC}$ vanishes, because the cocycle representing
$\widetilde\omega^i_{AC}$ has a multiple of $\dd x_j$ or $x_j$ in each term of the form $\Omega$. 
Therefore,  $\gamma^i_A\cup\omega^j_{B}\cup[\omega_C]=0$. 

The rest of the elements, which span the middle cohomology, have the form ${\varphi_l}_!(a)$ for some $a\in H^{d-3}(X\cap D_l)$.
The  projection formula gives 
$\gamma^i_A\cup\omega^j_{B}\cup{\varphi_l}_!(a)={\varphi_l}_!(\varphi_l^*(\gamma^i_A\cup\omega^j_{B})\cup a)$.
Hence, it suffices to show that $\varphi_l^*(\gamma^i_A\cup\omega^j_{B})=0$.
In further calculations, for simplicity, we assume that 
$A\in S_{\beta_1^{\sigma^1}}$ and $B\in S_{\beta-\beta_0+\beta_1^{\sigma^2}}$.
We will need to use a refinement
$\widetilde{\cal U}$ of the cover ${\cal U}$, by the 
open sets $\widetilde{U}_k=\{x\in\ps:x_k f_k(x)\ne0\}$ for $k=1,\dots,n$. Since $X$ is regular, these sets
cover the toric variety $\ps$.
In this case, the cup product $\gamma^i_A\cup\omega^j_{B}$ is represented by the \v{C}ech cocycle
$${\Biggl\{\frac{(-1)^{d}AB}{(\prod_{\rho_k\subset\sigma^1} x_k)(\prod_{\rho_k\subset\sigma^2} x_k)}
\biggl(\frac{u_{i_1,j_1}^i}{f_{i_1}}-
\frac{u_{i_0,j_0}^i}{f_{i_0}}\biggr)\lrcorner
\biggl(\frac{K_{i_2}(\partial^j_{k_2}\lrcorner\Omega)}{f_{i_2}}-
\frac{K_{i_1}(\partial^j_{k_1}\lrcorner\Omega)}{f_{i_1}}\biggr)\Biggr\}}_I,$$ 
where the index set $I=\{(i_0,{j_0},k_0),(i_1,{j_1},k_1),(i_2,{j_2},k_2)\}$ corresponds to the refinement
of $\widetilde{\cal U}$, ${\cal U}^{\sigma^1}$ and ${\cal U}^{\sigma^2}$,
and where $u_{i_t,j_t}^i$ denotes $\langle\partial_{i_t}\wedge\partial^i_{j_t},\dd f\rangle$. 
Note that 
$$\frac{\langle\partial_{i_0}\wedge\partial^i_{j_0},\dd f\rangle}{f_{i_0}}\lrcorner
\frac{K_{i_1}(\partial^j_{k_1}\lrcorner\Omega)}{f_{i_1}}=
\frac{K_{i_1}(\partial^j_{k_1}\lrcorner\partial^i_{j_0}\lrcorner\Omega)}{f_{i_1}}
+\frac{\langle\partial^i_{j_0},\dd f\rangle K_{i_1} K_{i_0}(\partial^j_{k_1}\lrcorner\Omega)}
{f_{i_0}f_{i_1}}.$$ 
For $\rho_l$, not lying in the cones $\sigma^1$ and $\sigma^2$,
the restriction $\varphi_l^*$ of the above cocycle vanishes:
if $i$ is among $\{i_0,i_1,i_2\}$, then $\widetilde{U}_i\cap D_i$ is empty;
if $i\notin\{i_0,i_1,i_2\}$,
each term of the cocycle is multiple of $x_l$ or $\dd x_l$ coming from  $\Omega$.
We are left to consider $\rho_l\subset\sigma^1\cup\sigma^2$.
For $\rho_l\subset\sigma^1$, we will show that the restriction $\varphi_l^*$ of the cocycle is  a \v Cech coboundary;
 the other case is similar. 
Compute
\begin{multline*}
\biggl(\frac{u_{i_1,j_1}^i}{f_{i_1}}-
\frac{u_{i_0,j_0}^i}{f_{i_0}}\biggr)\lrcorner
\biggl(\frac{K_{i_2}(\partial^j_{k_2}\lrcorner\Omega)}{f_{i_2}}-
\frac{K_{i_1}(\partial^j_{k_1}\lrcorner\Omega)}{f_{i_1}}\biggr)
\\=
\sum_{\tilde I=I\setminus\{(i_s,j_s,k_s)\}}(-1)^s
\biggl(
\frac{K_{\tilde i_1}(\partial^j_{\tilde k_1}\lrcorner
(\partial^i_{\tilde j_0}-\partial^i_{\tilde j_1})\lrcorner\Omega)}{f_{\tilde i_1}}
+\frac{\langle\partial^i_{\tilde j_0},\dd f\rangle K_{\tilde i_1} K_{\tilde i_0}(\partial^j_{\tilde k_0}\lrcorner\Omega)}
{f_{\tilde i_0}f_{\tilde i_1}}
\biggr).
\end{multline*}
Using this, we can see that the restriction $\varphi_l^*(\gamma^i_A\cup\omega^j_{B})$ 
is represented by a \v Cech coboundary because of the following observations.
The polynomial $\langle\partial^i_{\tilde j_0},\dd f\rangle$ is divisible by $x_l$. 
If $\rho_l\not\subset\sigma^1_{j_t}$, then the restricted open set  $U_{\sigma^1_{j_t}}\cap D_l$ is empty.
If $\sigma^1_{\tilde j_0}$ and  $\sigma^1_{\tilde j_1}$ contain $\rho_l$,
then $K_{\tilde i_1}(\partial^i_{\tilde j_0}-\partial^i_{\tilde j_1})\lrcorner\Omega$ is either 0 or divisible by $x_l$
because of equation (\ref{e:eur}).
Thus, the restriction $\varphi_l^*(\gamma^i_A\cup\omega^j_{B})=0$, and the result follows.
\end{pf}

At this point,
let us summarize our calculations of the cup products  $H^*(X,\wedge^*{\cal T}_X)$
with the middle cohomology $H^*(X,\Omega^{d-1-*}_X)$ for $d$-semiample regular hypersurfaces.
We have the elements in $H^*(X,\wedge^*{\cal T}_X)$ represented by $\gamma_{\_}$, $\gamma^i_{\_}$ (with $\rho_i$ lying in some
$\sigma\in\Sigma_X(2)$ such that
$\rho_i\notin\Sigma_X$), 
and the corresponding elements in $H^*(X,\Omega^{d-1-*}_X)$ represented by $\omega_{\_}$, $\omega^i_{\_}$.
Theorem~\ref{t:mod}  provides $\gamma_{A}\cup\omega_{B}=\omega_{AB}$, while Theorems~\ref{t:nonm}~and~\ref{t:nonm1} have
$\gamma_{A}\cup\omega^i_{B}=\omega^i_{AB}$ and  $\gamma^i_{A}\cup\omega_{B}=\omega^i_{AB}$.
Lemmas~\ref{l:van1}~and~\ref{l:van2} tell us that  the cup product 
$\gamma^i_{A}\cup\omega^j_{B}=0$, for $i\ne j$, unless $\rho_i$ and $\rho_j$ span a 2-dimensional cone of $\Sigma$.
Thus, for the constructed elements in $H^*(X,\wedge^*{\cal T}_X)$ and $H^*(X,\Omega^{d-1-*}_X)$, we are missing
only the cup products $\gamma^i_{A}\cup\omega^j_{B}$ when
$\rho_i$ and $\rho_j$ ($i$ may be equal to $j$) span a cone of $\Sigma$ contained in some  
2-dimensional cone of $\Sigma_X$.

First, we consider the nontrivial cup products $\gamma^i_{A}\cup\omega^j_{B}$ lying in 
$H^{d-1}(X,O_X)$, which is isomorphic to $R_1(f)_{d\beta-\beta_0}$, by Theorem~\ref{t:midn}.
We note here that the inclusion
\begin{equation}\label{e:isomorm}
\mu:R_1(f)_{d\beta-\beta_0}@>\prod_{k=1}^n x_k>>R_0(f)_{d\beta}
\end{equation}
induced by the multiplication is an isomorphism
because the dimensions of the spaces is the same number (of the interior integral points of a polytope $\Delta$ 
corresponding to $\beta$)
by  the isomorphism
$R_1(f)_{d\beta-\beta_0}=H^{0,d-1}(H^{d-1}(X\cap\ttt))$
of   \cite[Theorem~4.4]{m}, by \cite[Section~5.8]{dk}, and  by
  \cite[Theorem~11.5]{bc} with \cite[Corollary~3.14]{b1}.
 The cup product should be represented by a polynomial in the above spaces.

\begin{pr}\label{p:comp} Let $X\subset\ps$ be  a  $d$-semiample regular hypersurface
defined by $f\in S_\beta$, and denote 
$$G^\sigma(f):=\frac{x_{s}f_{s}x_{t}f_{t}\prod_{\rho_k\not\subset\sigma} x_k}{\mult(\sigma)\prod_{\rho_k\subset\sigma} x_k}
\in S_{2\beta+\beta_0-2\beta_1^\sigma}$$
 for  
$\sigma\in\Sigma_X(2)$
spanned by $\rho_s$ and $\rho_t$. 
 Given $A\in S_{(p-1)\beta+\beta_1^{\sigma}}$, $B\in S_{(d-1-p)\beta-\beta_0+\beta_1^{\sigma}}$,
then

{\rm (i)} for  $\rho_i=\rho_{l_k}\notin\Sigma_X$, as in (\ref{e:pic}), contained in $\sigma\in\Sigma_X(2)$:   
$$\gamma^i_{A}\cup\omega^i_{B}=\frac{\mult(\sigma_k+\sigma_{k+1})
[\omega_{\mu^{-1}(ABG^\sigma(f))}]}{\mult(\sigma_k)\mult(\sigma_{k+1})}\text { in }H^{d-1}(X,O_X),$$

{\rm (ii)} for $\rho_i,\rho_j\notin\Sigma_X$ which span 
a $2$-dimensional cone  $\sigma'\in\Sigma$ contained in $\sigma\in\Sigma_X(2)$:
$$\gamma^i_{A}\cup\omega^j_{B}=-\frac{[\omega_{\mu^{-1}(ABG^\sigma(f))}]}{\mult(\sigma')}\text { in }H^{d-1}(X,O_X).$$
\end{pr}

\begin{pf} To simplify the proof we assume that $p=1$.

(i) After a simple modification it follows  that
the cup product $\gamma^i_{A}\cup\omega^i_{B}$ is represented by the  cocycle
$${\Biggl\{\frac{(-1)^{d+((d-3)^2/2)}AB}{(\prod_{\rho_k\subset\sigma} x_k)^2}
\sum_{\tilde{I}=I\setminus\{(i_k,j_k)\}}(-1)^k
\frac{\langle\partial_{\tilde i_0}\wedge\partial^i_{\tilde j_0},\dd f\rangle}{f_{\tilde i_0}} 
\lrcorner
\frac{K_{\tilde i_{d-2}}\cdots K_{\tilde i_{1}}
(\partial^i_{\tilde j_{1}}\lrcorner\Omega)}
{f_{\tilde i_1}\cdots f_{\tilde i_{d-2}}}\Biggr\}}_I,
$$ 
where $I=\{(i_0,j_0),\dots,(i_{d-1},j_{d-1})\}$ is the index set corresponding to the open sets $\widetilde U_{i_k}$
(defined in Lemma~\ref{l:van2}) and $U_{\sigma_{j_k}}$.
Note that 
\begin{multline*}
\frac{AB}{(\prod_{\rho_k\subset\sigma} x_k)^2}
\frac{\langle\partial_{\tilde i_0}\wedge\partial^i_{\tilde j_0},\dd f\rangle}{f_{\tilde i_0}} 
\lrcorner
\frac{K_{\tilde i_{d-2}}\cdots K_{\tilde i_{1}}
(\partial^i_{\tilde j_{1}}\lrcorner\Omega)}
{f_{\tilde i_1}\cdots f_{\tilde i_{d-2}}}
\\
=
\frac{(-1)^{d-1}AB}{(\prod_{\rho_k\subset\sigma} x_k)^2}\biggl(\langle\partial^i_{\tilde j_0},\dd f\rangle
\frac{K_{\tilde i_{d-2}}\cdots K_{\tilde i_{0}}
(\partial^i_{\tilde j_{1}}\lrcorner\Omega)}
{f_{\tilde i_0}\cdots f_{\tilde i_{d-2}}}
+\frac{K_{\tilde i_{d-2}}\cdots K_{\tilde i_{1}}
(\partial^i_{\tilde j_{1}}\lrcorner\partial^i_{\tilde j_{0}}\lrcorner\Omega)}
{f_{\tilde i_1}\cdots f_{\tilde i_{d-2}}}\biggr).
\end{multline*}
The first summand is well defined on the corresponding open set: if $i\in\{\tilde i_0,\dots,\tilde i_{d-2}\}$, then
$x_i\ne0$ on the open set; otherwise, $\langle\partial^i_{\tilde j_0},\dd f\rangle K_{\tilde i_{d-2}}\cdots K_{\tilde i_{0}}
(\partial^i_{\tilde j_{1}}\lrcorner\Omega)$ is a multiple of $(x_i)^2$.
Therefore, the corresponding sum in the above cocycle forms a \v Cech coboundary, and 
$\gamma^i_{A}\cup\omega^i_{B}$ is represented by
$${\Biggl\{\frac{(-1)^{1+((d-3)^2/2)}AB}{(\prod_{\rho_k\subset\sigma} x_k)^2}
\sum_{\tilde{I}=I\setminus\{(i_k,j_k)\}}(-1)^k
\frac{K_{\tilde i_{d-2}}\cdots K_{\tilde i_{1}}
(\partial^i_{\tilde j_{1}}\lrcorner\partial^i_{\tilde j_{0}}\lrcorner\Omega)}
{f_{\tilde i_1}\cdots f_{\tilde i_{d-2}}}\Biggr\}}_I.
$$ 

By Proposition 5.3 in \cite{c2}, the polynomials $x_{r_0}f_{r_0},\dots,x_{r_{d-1}}f_{r_{d-1}}$ do
not vanish simultaneously on $X$ if $e_{r_0},\dots,e_{r_{d-1}}$ are linearly independent. We can always find
such generators so that $e_{r_0}=e_{l_0}$ and $e_{r_1}=e_{l_{n(\sigma)+1}}$ as in (\ref{e:pic}). 
Since the open sets $\{x\in\ps: x_{r_k}f_{r_k}\ne0\}$ cover the toric variety,
we can assume that
the first index in $I$ takes only the ordered values $r_0,\dots,r_{d-1}$.
In this case, it is not difficult to check that the above cocycle is different by a coboundary {from}
$${\Biggl\{\frac{(-1)^{1+((d-3)^2/2)}AB}{(\prod_{\rho_k\subset\sigma} x_k)^2}
\sum_{\tilde{I}=I\setminus\{(i_k,j_k)\}}(-1)^k\alpha_{\tilde{I}}
\frac{K_{\tilde i_{d-2}}\cdots K_{\tilde i_{1}}
(\partial_{l_{k+1}}\lrcorner\partial_{l_{k-1}}\lrcorner\Omega)}
{f_{\tilde i_1}\cdots f_{\tilde i_{d-2}}\mult(\sigma_k)\mult(\sigma_{k+1})}\Biggr\}}_I,$$
where 
$$\alpha_{\tilde{I}}=\left\{
\begin{matrix} 
-1 & \mbox{ if } \tilde i_0=\tilde i_1=r_2, \tilde j_0\le k< k+1\le \tilde j_1 \\
-1 & \mbox{ if } \tilde i_0=r_1, \tilde i_1=r_2, \tilde j_1\ge k+1 \\
1 & \mbox{ if } \tilde i_0=r_0, \tilde i_1=r_2, \tilde j_1\le k \\
0 & \mbox{ in all other cases.}
\end{matrix} \right.
$$
Using the Euler identities in the proof of Proposition~\ref{p:gysin}, the last cocycle converts to
$$\frac{\mult(\sigma_k+\sigma_{k+1})}{\mult(\sigma_k)\mult(\sigma_{k+1})}{\Biggl\{\frac{(-1)^{(d-3)^2/2}
ABx_{l_0}f_{l_0}x_{l_{n(\sigma)+1}}f_{l_{n(\sigma)+1}}}
{\mult(\sigma)(\prod_{\rho_k\subset\sigma} x_k)^2}
\frac{K_{i_{d-1}}\cdots K_{i_0}\Omega}
{f_{i_0}\cdots f_{i_{d-1}}}
\Biggr\}}_I.
$$
This cocycle represents 
$$\frac{\mult(\sigma_k+\sigma_{k+1})
[\omega_{\mu^{-1}(ABG^\sigma(f))}]}{\mult(\sigma_k)\mult(\sigma_{k+1})}$$
in $H^{d-1}(X,O_X)$.

(ii) Similar to the  proof  of the previous part and Lemma~\ref{l:van2}, the cup product $\gamma^i_{A}\cup\omega^j_{B}$ is 
represented by the following cocycle:
$${\Biggl\{\frac{(-1)^{1+((d-3)^2/2)}AB}{(\prod_{\rho_k\subset\sigma} x_k)^2}
\sum_{\tilde{I}=I\setminus\{(i_k,j_k)\}}(-1)^k
\frac{K_{\tilde i_{d-2}}\cdots K_{\tilde i_{1}}
(\partial^j_{\tilde j_{1}}\lrcorner(\partial^i_{\tilde j_{0}}-\partial^i_{\tilde j_{1}})\lrcorner\Omega)}
{f_{\tilde i_1}\cdots f_{\tilde i_{d-2}}}\Biggr\}}_I.
$$ 
The 1-dimensional cones $\rho_i$ and $\rho_j$ span one of the 2-dimensional cones 
$\sigma_k\subset\sigma$ as in in (\ref{e:pic}). 
The  cocycle  differs by a coboundary {from} 
$${\Biggl\{\frac{(-1)^{1+((d-3)^2/2)}AB}{(\prod_{\rho_k\subset\sigma} x_k)^2}
\sum_{\tilde{I}=I\setminus\{(i_k,j_k)\}}(-1)^k\alpha_{\tilde{I}}
\frac{-K_{\tilde i_{d-2}}\cdots K_{\tilde i_{1}}
(\partial_{l_{k}}\lrcorner\partial_{l_{k-1}}\lrcorner\Omega)}
{f_{\tilde i_1}\cdots f_{\tilde i_{d-2}}\mult(\sigma_k)^2}\Biggr\}}_I,$$
where 
$\alpha_{\tilde{I}}$ is the same as in part (i). The Euler identities show that this represents 
$$-\frac{[\omega_{\mu^{-1}(ABG^\sigma(f))}]}{\mult(\sigma_k)}\in H^{d-1}(X,O_X).$$
\end{pf}

The restriction maps $\varphi^*_l$, induced by the inclusions $\varphi_l:X\cap D_l\hookrightarrow X$, give some
information about the  nontrivial cup products $\gamma^i_A\cup\omega^i_B$ in $H^{d-1}(X)$. 
We will use this in Section~\ref{s:cal3} to calculate  nontrivial triple products on the chiral ring
of anticanonical hypersurfaces.

\begin{pr}\label{p:restr} Let $X\subset\ps$ be  a  $d$-semiample regular hypersurface
defined by $f\in S_\beta$, and let, as in (\ref{e:pic}),
$\rho_i=\rho_{l_k}\notin\Sigma_X$ be in some $\sigma\in\Sigma_X(2)$.
 Then, for $A\in S_{p\beta+\beta_1^{\sigma}}$ and $B\in S_{q\beta-\beta_0+\beta_1^{\sigma}}$,

{\rm (i)} $\varphi^*_{l_{k\pm1}}(\gamma^i_A\cup\omega^i_B)=
\pm\varphi^*_{l_{k\pm1}}\bigl(\omega^i_{ABH^{\sigma}_{{i,\pm1}}(f)}\bigr)$,
where $H^{\sigma}_{{i,\pm1}}(f)$ is a polynomial in $S_{\beta-\beta_1^{\sigma}}$ equal to 
$\sqrt{-1}x_{l_{k\pm1}}f_{l_{k\pm1}}/(\mult(\sigma_{k,k\pm1})\prod_{\rho_k\subset\sigma} x_k)$ 
at $x_{l_k}=0$ and $x_{l_{k\pm1}}=0$, where
$\sigma_{s,t}$ denotes the cone spanned by $\rho_{l_s}$ and $\rho_{l_{t}}$.

{\rm (ii)} $\varphi^*_i(\gamma^i_A\cup\omega^i_B)=
\varphi^*_i\bigl(\omega^i_{ABH^\sigma_i(f)}\bigr)$, 
where $H^\sigma_i(f)$ is a polynomial in $S_{\beta-\beta_1^{\sigma}}$ equal to 
$$
\frac{\sqrt{-1}x_{l_{k+1}}f_{l_{k+1}}}{\mult(\sigma_{k,k+1})\prod_{\rho_k\subset\sigma} x_k}-
\frac{\sqrt{-1}x_{l_{k-1}}f_{l_{k-1}}}{\mult(\sigma_{k-1,k})\prod_{\rho_k\subset\sigma} x_k}
$$ 
with $x_{i}=x_{l_{k-1}}=x_{l_{k+1}}=0$.
\end{pr}

\begin{pf} For simplicity, we assume that $A\in S_{\beta_1^{\sigma}}$ and $B\in S_{\beta-\beta_0+\beta_1^{\sigma}}$.
The cup product $\gamma^i_A\cup\omega^i_B$ is represented by the \v Cech cocycle:
\begin{equation}\label{e:cup1}
{\biggl\{\frac{(-1)^{d}AB}{(\prod_{\rho_k\subset\sigma} x_k)^2}
\biggl(\frac{u_{i_1,j_1}^i}{f_{i_1}}-
\frac{u_{i_0,j_0}^i}{f_{i_0}}\biggr)\lrcorner
\biggl(\frac{K_{i_2}(\partial^j_{k_2}\lrcorner\Omega)}{f_{i_2}}-
\frac{K_{i_1}(\partial^j_{k_1}\lrcorner\Omega)}{f_{i_1}}\biggr)\biggr\}}_I,
\end{equation}
where $u_{i_k,j_k}^i$ denotes $\langle\partial_{i_k}\wedge\partial^i_{j_k},\dd f\rangle$, 
and $I=\{(i_0,j_0),(i_1,j_1),(i_2,j_2)\}$ is the index set.
Compute
\begin{multline*}
\frac{\langle\partial_{i_0}\wedge\partial^i_{j_0},\dd f\rangle}{f_{i_0}}\lrcorner
\frac{K_{i_1}(\partial^i_{j_1}\lrcorner\Omega)}{f_{i_1}}=
\frac{K_{i_1}(\partial^i_{ j_1}\lrcorner\partial^i_{j_0}\lrcorner\Omega)}{f_{i_1}}
+\frac{\langle\partial^i_{j_0},\dd f\rangle K_{i_1} K_{i_0}(\partial^i_{j_1}\lrcorner\Omega)}
{f_{i_0}f_{i_1}}\\
=\frac{K_{i_0}(\partial^i_{ j_1}\lrcorner\partial^i_{j_0}\lrcorner\Omega)}{f_{i_0}}
+\frac{\langle\partial^i_{j_1},\dd f\rangle K_{i_1} K_{i_0}(\partial^i_{j_0}\lrcorner\Omega)}
{f_{i_0}f_{i_1}}
+
\frac{K_{i_1}K_{i_0}(\partial^i_{j_1}\lrcorner\partial^i_{j_0}\lrcorner(\dd f\wedge\Omega)}{f_{i_0}f_{i_1}}
\\
\equiv 
\frac{K_{i_0}(\partial^i_{ j_1}\lrcorner\partial^i_{j_0}\lrcorner\Omega)}{f_{i_0}}
+\frac{\langle\partial^i_{j_1},\dd f\rangle K_{i_1} K_{i_0}(\partial^i_{j_0}\lrcorner\Omega)}
{f_{i_0}f_{i_1}},
\end{multline*}
where,  as in Lemma~\ref{l:cobi1}, we used $\dd f\wedge\Omega\equiv0$ modulo multiples of $f$.
Hence, the cup product $\gamma^i_A\cup\omega^i_B$ is represented by the \v Cech cocycle
\begin{equation*}
{\biggl\{\frac{(-1)^{d}AB}{(\prod_{\rho_k\subset\sigma} x_k)^2}
\sum_{\tilde I=I\setminus\{(i_s,j_s)\}}(-1)^s
{\biggl(\frac{K_{\tilde i_0}(\partial^i_{\tilde j_1}\lrcorner\partial^i_{\tilde j_0}\lrcorner\Omega)}{f_{\tilde i_0}}
+\frac{\langle\partial^i_{\tilde j_1},\dd f\rangle K_{\tilde i_1} K_{\tilde i_0}(\partial^i_{\tilde j_0}\lrcorner\Omega)}
{f_{\tilde i_0}f_{\tilde i_1}}
\biggr)\biggr\}}}_I.
\end{equation*}

For part (i), consider the restriction $\varphi^*_{l_{k\pm1}}$ of this cocycle.
Note that 
the open set $U_{\sigma_j}\cap D_{l_{k\pm1}}$ is empty, if $\sigma_j$ does not contain $\rho_{l_{k\pm1}}$,
and that $\partial^i_j=0$, if the corresponding cone $\sigma_j$ does not contain $\rho_i$.
Using this and $\partial^i_j\wedge\partial^i_j=0$, we get that
the restriction $\varphi^*_{l_{k\pm1}}$ of $\gamma^i_A\cup\omega^i_B$ is  represented by
\begin{equation}\label{e:rest}
\pm\sqrt{-1}{\biggl\{\frac{(-1)^{d+(1/2)}AB}{(\prod_{\rho_k\subset\sigma} x_k)^2}
\sum_{\tilde I=I\setminus\{(i_s,j_s)\}}(-1)^s
\frac{x_{l_{k\pm1}}f_{l_{k\pm1}}}{\mult(\sigma_{k,k\pm1})} 
\frac{K_{\tilde i_1} K_{\tilde i_0}(\partial^i_{\tilde j_0}\lrcorner\Omega)}
{f_{\tilde i_0}f_{\tilde i_1}}
\biggr\}}_I,
\end{equation}
where the index set $I$ corresponds to the restricted open cover ${\cal U}^\sigma|_{X\cap D_{l_{k\pm1}}}$, 
and where $\sigma_{k,k\pm1}$ is the cone generated by $e_{l_k}$ and 
$e_{l_{k\pm1}}$. Notice that this cocycle is similar to the restriction 
$\varphi^*_{l_{k\pm1}}(\omega^i_{C})$ for some polynomial $C$. The problem here is that
$x_{l_{k\pm1}}f_{l_{k\pm1}}$ is not necessarily divisible by $\prod_{\rho_k\subset\sigma} x_k$. 
So some work is required to get the correct polynomial.
Let $X$ be linearly equivalent to a torus invariant divisor $D=\sum_{k=1}^nb_k D_k$ with 
the associated polytope $\Delta=\Delta_D$ given by the conditions $b_l+\langle m,e_l\rangle\ge0$ . Then we can write 
$f=\sum_{m\in\Delta\cap M}a_m x^{D(m)}$, where $x^{D(m)}$ denotes $\prod_{l=1}^n x_l^{b_l+\langle m,e_l\rangle}$. 
Note that 
$$x_{l_{k\pm1}}f_{l_{k\pm1}}=\sum_{m\in\Delta\cap M}a_m(b_{l_{k\pm1}}+\langle m,e_{l_{k\pm1}}\rangle) x^{D(m)}.$$
If $b_{l_{k\pm1}}+\langle m,e_{l_{k\pm1}}\rangle=0$, then the corresponding monomial $x^{D(m)}$ is not present in 
$x_{l_{k\pm1}}f_{l_{k\pm1}}$. On the other hand, if  $b_{l_{k\pm1}}+\langle m,e_{l_{k\pm1}}\rangle>1$, then
the multiple of the corresponding monomial $x^{D(m)}$ in (\ref{e:rest}) vanishes, since 
$\partial^i_{\tilde j_0}=\mp x_{l_{k\pm1}}\partial_{l_{k\pm1}}$ or 0. 
By the argument in the proof of Lemma~\ref{l:non}, $b_{l_{k\pm1}}+\langle m,e_{l_{k\pm1}}\rangle=1$ implies that
$b_l+\langle m,e_l\rangle>0$ for all $\rho_l\subset\sigma$ such that $\rho_l\notin\Sigma_X(1)$.
 If $b_i+\langle m,e_i\rangle>1$, the multiple of
the monomial $x^{D(m)}$ in (\ref{e:rest}) forms a coboundary. 
Therefore, only the monomials $x^{D(m)}$ in $x_{l_{k\pm1}}f_{l_{k\pm1}}$, satisfying 
$b_{l_{k\pm1}}+\langle m,e_{l_{k\pm1}}\rangle=1$
and $b_i+\langle m,e_i\rangle=1$, have a nontrivial contribution in the \v Cech cocycle (\ref{e:rest}). 
For all such  monomials, it follows {from} the relations of the cone generators in the proof of Proposition~\ref{p:gysin}
 that $b_{s}+\langle m,e_{s}\rangle>0$ with $s=l_0,l_{n(\sigma)+1}$. Hence,
the monomials are divisible by 
 $\prod_{\rho_k\subset\sigma} x_k$. Thus, 
$\varphi^*_{l_{k\pm1}}(\gamma^i_A\cup\omega^i_B)=\pm\varphi^*_{l_{k\pm1}}\bigl(\omega^i_{ABH^\sigma_{i,\pm1}(f)}\bigr)$,
where $H^\sigma_{i,\pm1}(f)$ is the polynomial
 $$\sum_m\frac{\sqrt{-1}\sum_{m}a_m x^{D(m)}}{\mult(\sigma_{k,k\pm1})\prod_{\rho_k\subset\sigma} x_k}$$
with the sum over all $m\in\Delta\cap M$, satisfying the equalities 
$b_{l_{k\pm1}}+\langle m,e_{l_{k\pm1}}\rangle=1$
and $b_i+\langle m,e_i\rangle=1$. This is the same as 
$\sqrt{-1}x_{l_{k\pm1}}f_{l_{k\pm1}}/(\mult(\sigma_{k,k\pm1})\prod_{\rho_k\subset\sigma} x_k)$
evaluated at $x_i=0$ and $x_{l_{k\pm1}}=0$.

In part (ii) we will need to use the refinement
$\widetilde{\cal U}$ of the cover ${\cal U}$ defined in Lemma~\ref{l:van2}. 
{From} (\ref{e:cup1}) we get that the cup product $\gamma^i_A\cup\omega^i_B$ is  represented by the 
\v Cech cocycle
$$
{\biggl\{\frac{(-1)^{d}AB}{(\prod_{\rho_k\subset\sigma} x_k)^2}
\sum_{\tilde I=I\setminus\{(i_s,j_s)\}}(-1)^s
{\biggl(\frac{K_{\tilde i_1}(\partial^i_{\tilde j_1}\lrcorner\partial^i_{\tilde j_0}\lrcorner\Omega)}{f_{\tilde i_1}}
+\frac{\langle\partial^i_{\tilde j_0},\dd f\rangle K_{\tilde i_1} K_{\tilde i_0}(\partial^i_{\tilde j_1}\lrcorner\Omega)}
{f_{\tilde i_0}f_{\tilde i_1}}
\biggr)\biggr\}}}_I,
$$
where the index set $I=\{(i_0,j_0),(i_1,j_1),(i_2,j_2)\}$ corresponds to the refinement of 
$\widetilde{\cal U}$ and ${\cal U}^\sigma|_X$.
Notice 
\begin{multline*}
\frac{K_{i_1}(\partial^i_{ j_1}\lrcorner\partial^i_{j_0}\lrcorner\Omega)}{f_{i_1}}-
\frac{K_{i_2}(\partial^i_{ j_1}\lrcorner\partial^i_{j_0}\lrcorner\Omega)}{f_{i_2}}
\\
=
\frac{\langle\partial^i_{j_0},\dd f\rangle K_{i_2} K_{i_1}(\partial^i_{j_1}\lrcorner\Omega)}
{f_{i_1}f_{i_2}}-
\frac{\langle\partial^i_{j_1},\dd f\rangle K_{i_2} K_{i_1}(\partial^i_{j_0}\lrcorner\Omega)}
{f_{i_1}f_{i_2}}-
\frac{ K_{i_2} K_{i_1}(\partial^i_{j_1}\lrcorner\partial^i_{j_0}\lrcorner(\dd f\wedge\Omega)}
{f_{i_1}f_{i_2}}
\\
\equiv
\frac{\langle\partial^i_{j_0},\dd f\rangle K_{i_2} K_{i_1}(\partial^i_{j_1}\lrcorner\Omega)}
{f_{i_1}f_{i_2}}-
\frac{\langle\partial^i_{j_1},\dd f\rangle K_{i_2} K_{i_1}(\partial^i_{j_0}\lrcorner\Omega)}
{f_{i_1}f_{i_2}},
\end{multline*}
since $\dd f\wedge\Omega\equiv0$ modulo multiples of $f$.
Using this, we compute
\begin{align*}
&\frac{AB}{(\prod_{\rho_k\subset\sigma} x_k)^2} 
\sum_{\tilde I=I\setminus\{(i_s,j_s)\}}(-1)^s
\biggl(\frac{K_{\tilde i_1}(\partial^i_{\tilde j_1}\lrcorner\partial^i_{\tilde j_0}\lrcorner\Omega)}{f_{\tilde i_1}}
+\frac{\langle\partial^i_{\tilde j_0},\dd f\rangle K_{\tilde i_1} K_{\tilde i_0}(\partial^i_{\tilde j_1}\lrcorner\Omega)}
{f_{\tilde i_0}f_{\tilde i_1}}\biggr)
\\
&=
\sum_{\tilde I=I\setminus\{(i_s,j_s)\}}(-1)^s
\frac{AB\langle\partial^i_{\tilde j_0},\dd f\rangle K_{\tilde i_1} K_{\tilde i_0}((\partial^i_{\tilde j_1}-
\partial^i_{\tilde j_0})\lrcorner\Omega)}
{(\prod_{\rho_k\subset\sigma} x_k)^2f_{\tilde i_0}f_{\tilde i_1}}
\\
&+
\frac{AB K_{i_2}((\partial^i_{j_0}\wedge\partial^i_{j_1}-\partial^i_{j_0}\wedge\partial^i_{j_2}+
\partial^i_{j_1}\wedge\partial^i_{j_2})\lrcorner\Omega)}{(\prod_{\rho_k\subset\sigma} x_k)^2f_{i_2}}
\\
&+
\frac{AB\langle\partial^i_{j_1},\dd f\rangle}{(\prod_{\rho_k\subset\sigma} x_k)^2}
\biggl(\frac{K_{i_2} K_{i_0}(\partial^i_{j_0}\lrcorner\Omega)}
{f_{i_0}f_{i_2}}
-\frac{K_{i_1} K_{i_0}(\partial^i_{j_0}\lrcorner\Omega)}
{f_{i_0}f_{i_1}}
-\frac{K_{i_2} K_{i_1}(\partial^i_{j_0}\lrcorner\Omega)}{f_{i_1}f_{i_2}}\biggr)
\\
&+
\frac{AB(\langle\partial^i_{j_0},\dd f\rangle+\langle\partial^i_{j_1},\dd f\rangle)}{(\prod_{\rho_k\subset\sigma} x_k)^2}
\biggl(\sum_{\tilde I=I\setminus\{(i_s,j_s)\}}(-1)^s 
\frac{K_{\tilde i_1} K_{\tilde i_0}(\partial^i_{\tilde j_0}\lrcorner\Omega)}
{f_{\tilde i_0}f_{\tilde i_1}}
\biggr).
\end{align*}
It is not difficult to see  that the first summand produces a coboundary.
The open set $U_{\sigma_j}\cap D_i$ is empty unless  $\sigma_j$ contains $\rho_i$. Therefore,
applying the restriction $\varphi^*_i$, we can assume that the second component of the index $(i_0,j_0)$ takes only values
$l_{k-1}$ or $l_{k+1}$. In this case, $\partial^i_{j_0}\wedge\partial^i_{j_1}-\partial^i_{j_0}\wedge\partial^i_{j_2}+
\partial^i_{j_1}\wedge\partial^i_{j_2}$ (and the corresponding summand in the cocycle) vanishes. 
The third summand also ends up contributing zero:
\begin{multline*}
\frac{AB\langle\partial^i_{j_1},\dd f\rangle}{(\prod_{\rho_k\subset\sigma} x_k)^2}
\biggl(\frac{K_{i_2} K_{i_0}(\partial^i_{j_0}\lrcorner\Omega)}
{f_{i_0}f_{i_2}}-\frac{K_{i_1} K_{i_0}(\partial^i_{j_0}\lrcorner\Omega)}
{f_{i_0}f_{i_1}}
-\frac{K_{i_2} K_{i_1}(\partial^i_{j_0}\lrcorner\Omega)}{f_{i_1}f_{i_2}}\biggr)
\\
=\frac{AB\langle\partial^i_{j_1},\dd f\rangle}{(\prod_{\rho_k\subset\sigma} x_k)^2}
\biggl(\frac{\langle\partial^i_{j_0},\dd f\rangle K_{i_2} K_{i_1}K_{i_0}\Omega}
{f_{i_0}f_{i_1}f_{i_2}}+
\frac{ K_{i_2} K_{i_1}K_{i_0}\partial^i_{j_0}\lrcorner(\dd f\wedge\Omega)}
{f_{i_0}f_{i_1}f_{i_2}}\biggr)
\\
\equiv
\frac{AB\langle\partial^i_{j_0},\dd f\rangle\langle\partial^i_{j_1},\dd f\rangle}{(\prod_{\rho_k\subset\sigma} x_k)^2}
\frac{ K_{i_2} K_{i_1}K_{i_0}\Omega}
{f_{i_0}f_{i_1}f_{i_2}}.
\end{multline*}
If $i$ is among $\{i_0,i_1,i_2\}$, then this restricts to an empty set since $\widetilde{U}_i\cap D_i$ is empty. 
In the opposite case, this gives 0 under the restriction since $x_i$ or $\dd x_i$ is present in each term of  $\Omega$.
Thus, the restriction $\varphi^*_i(\gamma^i_A\cup\omega^i_B)$ is represented by the cocycle
$$
{\biggl\{\frac{(-1)^d AB(\langle\partial^i_{j_0},\dd f\rangle+\langle\partial^i_{j_1},\dd f\rangle)}
{(\prod_{\rho_k\subset\sigma} x_k)^2}
\sum_{\tilde I=I\setminus\{(i_s,j_s)\}}(-1)^s 
\frac{K_{\tilde i_1} K_{\tilde i_0}(\partial^i_{\tilde j_0}\lrcorner\Omega)}
{f_{\tilde i_0}f_{\tilde i_1}}\biggr\}}_I,
$$
where the index set $I$ now corresponds  to the open cover ${\cal U}^\sigma|_{X\cap D_i}$.
However, the last calculation shows that if $j_0$ coincides with $j_1$, then the expression in the above cocycle vanishes
on the given open set.
Hence, this cocycle is the same as
$$
{\biggl\{\frac{(-1)^{d} AB}
{(\prod_{\rho_k\subset\sigma} x_k)^2}
\biggl(\frac{x_{l_{k-1}}f_{l_{k-1}}}{\mult(\sigma_{k})}- 
\frac{x_{l_{k+1}}f_{l_{k+1}}}{\mult(\sigma_{k+1})}\biggr)
\sum_{\tilde I=I\setminus\{(i_s,j_s)\}}(-1)^s 
\frac{K_{\tilde i_1} K_{\tilde i_0}(\partial^i_{\tilde j_0}\lrcorner\Omega)}
{f_{\tilde i_0}f_{\tilde i_1}}\biggr\}}_I.
$$

By the arguments similar to part (i), one can show that this coincides with the restriction
$\varphi^*_i\bigl(\omega^i_{ABH^\sigma_i(f)}\bigr)$, 
where $H^\sigma_i(f)$ is equal to 
$$\frac{\sqrt{-1}x_{l_{k+1}}f_{l_{k+1}}}{\mult(\sigma_{k,k+1})\prod_{\rho_k\subset\sigma} x_k}- 
\frac{\sqrt{-1}x_{l_{k-1}}f_{l_{k-1}}}{\mult(\sigma_{k-1,k})\prod_{\rho_k\subset\sigma} x_k}$$
at $x_{l_k}=x_{l_{k-1}}=x_{l_{k+1}}=0$. 
\end{pf}

\begin{rem}\label{r:hs} In the anticanonical case $\beta=\beta_0$ the polynomials $H^\sigma_{i,\pm1}(f)$,  $H^\sigma_i(f)$ of the 
above proposition  can be written 
in a simpler form. Let $D=\sum_{k=1}^n D_k$ be the anticanonical divisor with  the associated polytope $\Delta:=\Delta_D$.
For $f=\sum_{m\in\Delta\cap M}a_m x^{D(m)}$ and $\sigma\in\Sigma_X(2)$, 
denote 
$$H^\sigma(f):=\sqrt{-1}\sum_{m\in\sigma^\perp\cap \Delta\cap M}a_m\frac{x^{D(m)}}{\prod_{\rho_k\subset\sigma} x_k}.$$
Then $H^\sigma_{i,\pm1}(f)=H^\sigma(f)/\mult(\sigma_{k,k\pm1})$ and  
$$H^\sigma_i(f)=\biggl(\frac{1}{\mult(\sigma_{k+1})}-\frac{1}{\mult(\sigma_{k})}\biggr)H^\sigma(f).$$
\end{rem}

\section{The chiral ring for anticanonical hypersurfaces}\label{s:cal3}

Here, we will apply the results of the previous sections to explicitly describe  a subring of the 
chiral ring $H^*(X,\wedge^*{\cal T}_X)$, coming {from} the graded pieces of  $R_1(f)$ and $R^\sigma_1(f)$, for  
semiample anticanonical regular hypersurfaces. 
By Proposition~\ref{p:cal}, such hypersurfaces are Calabi-Yau.
The description of the chiral ring is complete for Calabi-Yau threefolds.

Let $\ps$ be a complete simplicial toric variety, and let $X\subset\ps$ 
be a big and nef regular hypersurface defined by $f\in S_\beta$. {From} Theorem~\ref{t:midn}, we know the following part of the middle 
cohomology of $X$:
$$[\omega_{\_}]\oplus(\oplus_i\omega^i_{\_}):
R_1(f)_{(*+1)\beta-\beta_0}
\oplus\Biggl(\bigoplus\begin{Sb}\sigma
\in\Sigma_X(2)\end{Sb}\bigl(R^\sigma_1(f)_{*\beta-\beta_0+\beta_1^\sigma}\bigr)
^{n(\sigma)}\Biggr)\hookrightarrow  H^{d-1-*,*}(X).$$

Now suppose that $\beta$ is the anticanonical degree $\beta_0$. 
In this case, the isomorphism (\ref{f:simeq}) and Theorems~\ref{t:mod},~\ref{t:nonm} give us:

\begin{thm}\label{t:chiral}
Let $X\subset\ps$ be  a  semiample anticanonical regular hypersurface
defined by $f\in S_\beta$.  
Then there is a natural inclusion
$$\gamma_{\_}\oplus(\oplus_i\gamma^i_{\_}):
R_1(f)_{*\beta}
\oplus\Biggl(\bigoplus\begin{Sb}\sigma
\in\Sigma_X(2)\end{Sb}\bigl(R^\sigma_1(f)_{(*-1)\beta+\beta_1^\sigma}\bigr)
^{n(\sigma)}\Biggr)\hookrightarrow H^*(X,\wedge^*{\cal T}_X),$$
where  the sum $\oplus_i\gamma^i_{\_}$ is over $\rho_i\subset\sigma\in\Sigma_X(2)$  such that 
$\rho_i\notin\Sigma_X$ and    $n(\sigma)$ is the number of such cones. Also, 
$R^\sigma_1(f)_{(q-1)\beta+\beta_1^\sigma}=0$ for $q=0,d-1$.
\end{thm}

\begin{rem} The map given by 
$\gamma_{\_}\oplus(\oplus_i\gamma^i_{\_})$ is an isomorphism onto $H^q(X,\wedge^q{\cal T}_X)$ if $q=0,1,d-2,d-1$ and $d\ne1,3$.
In particular, for semiample anticanonical regular hypersurfaces of dimension 3,
 we get a complete description of the chiral ring.
\end{rem}

We claim that the part of $H^*(X,\wedge^*{\cal T}_X)$ given in the above theorem is a subring.
Let us describe the product structure on  this  part. First, note that the ring $\oplus_p H^p(X,\wedge^p{\cal T}_X)$ is 
commutative.
Theorems~\ref{t:hom},~\ref{t:ncom}, Lemmas~\ref{l:van},~\ref{l:van2} and equation  (\ref{f:simeq})
give us all information about the ring structure except for 
 the products $\gamma^i_{A}\cup\gamma^j_{B}$ when
$\rho_i$ and $\rho_j$  span a cone of $\Sigma$ contained in some  
2-dimensional cone $\sigma\in\Sigma_X$. 
For such $\rho_i$ and $\rho_j$, we first show that $\gamma^i_{A}\cup\omega^j_{B}$ is in the part of the middle cohomology
represented by  $[\omega_{\_}]\oplus(\oplus_k\omega^k_{\_})$.
It is easy to see that $\gamma^i_{A}\in H^*(X,\wedge^*{\cal T}_X)$ 
can be lifted to  $\gamma$ in $H^*(\ps,\wedge^*{\cal T}_\ps)$ with respect to the maps of Lemma~\ref{l:proj}.
Using this lemma, for $h\in H^{d-1}(\ps)$, we have
$$\gamma^i_{A}\cup\omega^j_{B}\cup i^*h=\pm\gamma^i_{A}\cup i^*h\cup\omega^j_{B}=\pm i^*(\gamma\cup h)\cup\omega^j_{B}=0$$
since the toric part is orthogonal to the residue part in the middle cohomology.
Similarly,  $\widetilde{\varphi_\tau^*}\gamma^i_{A}=\eta^*\widetilde{\gamma^i_{A}}$ for a corresponding
$\widetilde{\gamma^i_{A}}\in H^*(X\cap V(\tau),\wedge^*{\cal T}_{X\cap V(\tau)})$, $\tau\in\Sigma(2)$. Therefore, for 
$h'\in H^{d-5}(X\cap V(\tau))$,
$$\gamma^i_{A}\cup\omega^j_{B}\cup {\varphi_\tau}_!h'={\varphi_\tau}_!(\varphi_\tau^*(\gamma^i_{A}\cup\omega^j_{B})\cup h'=
{\varphi_\tau}_!(\widetilde{\gamma^i_{A}}\cup\varphi_\tau^*\omega^j_{B}\cup h')=
\pm{\varphi_\tau}_!(\widetilde{\gamma^i_{A}}\cup h')\cup\omega^j_{B}=0$$
where we use the projection formula for Gysin homomorphisms and Theorem~\ref{t:midn}.
Hence, by the same Theorem and because of the nondegenerate pairing on the middle cohomology,  
the cup product   $\gamma^i_{A}\cup\omega^j_{B}$ lies in the space given by $[\omega_{\_}]\oplus(\oplus_k\omega^k_{\_})$.
By the isomorphism (\ref{f:simeq}), the cup product $\gamma^i_{A}\cup\gamma^j_{B}$ is in the part of the chiral ring 
described in Theorem~\ref{t:chiral}. Thus, this part is a subring of $H^*(X,\wedge^*{\cal T}_X)$.

Since $X$ is Calabi-Yau, we have natural isomorphisms
$$H^{d-1}(X,\wedge^{d-1}{\cal T}_X)\cong H^{d-1}(X,O_X)\cong H^0(X,\Omega^{d-1}_X)\cong{\Bbb C}.$$
The cup product on the middle cohomology induces a nondegenerate pairing on the chiral ring and its subring represented
by $\gamma_{\_}\oplus(\oplus_k\gamma^k_{\_})$. 
Therefore, one can recover the product  structure of the subring, knowing the triple products  on this subring.
Because of Lemmas~\ref{l:van}~and~\ref{l:van2} it suffices to
consider the product of three elements 
$\gamma^i_{A}\cup\gamma^j_{B}\cup\gamma^l_{C}\in H^{d-1}(X,\wedge^{d-1}{\cal T}_X)$ 
in the cases $i=j=l$  and $i=j$ with $l$ such that $\rho_i$, $\rho_l$ span a 
2-dimensional cone of $\Sigma$ contained in some $\sigma\in\Sigma_X(2)$.
For this, compute
\begin{multline*}
(\gamma^i_{A}\cup\gamma^i_{B}\cup\gamma^i_{C}\cup [\omega_1])\cup [\omega_1]
=\varepsilon\gamma^i_{A}\cup\omega^i_{B}\cup\omega^i_{C}=
\varepsilon\gamma^i_{A}\cup\omega^i_{B}\cup{\varphi_i}_!\widetilde\omega^i_{C}
\\=
\varepsilon{\varphi_i}_!(\varphi^*_i(\gamma^i_{A}\cup\omega^i_{B})\cup\widetilde\omega^i_{C})
=
\varepsilon{\varphi_i}_!((\varphi^*_i\omega^i_{ABH^\sigma_i(f)})\cup\widetilde\omega^i_{C})
\\=
\varepsilon\omega^i_{ABH^\sigma_i(f)}\cup{\varphi_i}_!\widetilde\omega^i_{C}
=
\varepsilon\omega^i_{ABH^\sigma_i(f)}\cup\omega^i_{C}=
(\gamma^i_{ABH^\sigma_i(f)}\cup\omega^i_{C})\cup [\omega_1]
\\=
\frac{\mult(\sigma_k+\sigma_{k+1})
[\omega_{\mu^{-1}(ABCH^\sigma_i(f)G^\sigma(f))}]\cup [\omega_1]}{\mult(\sigma_k)\mult(\sigma_{k+1})},
\end{multline*}
where we used Propositions~\ref{p:gysin},~\ref{p:restr},~\ref{p:comp} and the projection formula for Gysin homomorphisms,
and where $\varepsilon$ is a sign depending on the degree of $C$.
Similarly, in the other case (as in (\ref{e:pic}), $\rho_i=\rho_{l_k}$):
\begin{multline*}
(\gamma^i_{A}\cup\gamma^i_{B}\cup\gamma^{l_{k\pm1}}_{C}\cup [\omega_1])\cup [\omega_1]
=\varepsilon\gamma^i_{A}\cup\omega^i_{B}\cup\omega^{l_{k\pm1}}_{C}=
\varepsilon\gamma^i_{A}\cup\omega^i_{B}\cup{\varphi_{l_{k\pm1}}}_!\widetilde\omega^{l_{k\pm1}}_{C}
\\=
\varepsilon{\varphi_{l_{k\pm1}}}_!(\varphi^*_{l_{k\pm1}}(\gamma^i_{A}\cup\omega^i_{B})\cup\widetilde\omega^{l_{k\pm1}}_{C})
=
\pm\varepsilon{\varphi_{l_{k\pm1}}}_!((\varphi^*_{l_{k\pm1}}\omega^i_{ABH_{i,\pm1}^\sigma(f)})\cup\widetilde\omega^{l_{k\pm1}}_{C})
\\=
\pm\varepsilon\omega^i_{ABH^\sigma_{i,\pm1}(f)}\cup\omega^{l_{k\pm1}}_{C}
=\pm(\gamma^i_{ABH^\sigma_{i,\pm1}(f)}\cup\omega^{l_{k\pm1}}_{C})\cup [\omega_1]
\\
=
\mp\frac{[\omega_{\mu^{-1}(ABCH^\sigma_{i,\pm1}(f)G^\sigma(f))}]\cup [\omega_1]}{\mult(\sigma_{k,k\pm1})}.
\end{multline*}

Since there is  an isomorphism $\cup[\omega_1]:H^{d-1}(X,O_X)\cong H^{d-1}(X,\Omega^{d-1}_X)$,
{from} the above calculation  we get an explicit product structure on the chiral ring.

\begin{thm}
Let $X\subset\ps$ be  a  semiample anticanonical regular hypersurface
defined by $f\in S_\beta$.  Then, under the identifications of Theorem~\ref{t:chiral}, we have

{\rm(i)} $\gamma_A\cup\gamma_B=\gamma_{AB}$,

{\rm(ii)}  $\gamma_{A}\cup\gamma^i_{B}=\gamma^i_{AB}$,

{\rm(iii)} $\gamma^i_{A}\cup\gamma^j_{B}=0$, $i\ne j$, unless $\rho_i$ and $\rho_j$ span a cone of $\Sigma$
contained in a $2$-dimensional cone of   $\Sigma_X$,

{\rm(iv)} for  $\rho_i=\rho_{l_k}\notin\Sigma_X$, as in (\ref{e:pic}), contained in $\sigma\in\Sigma_X(2)$ and
$A,B\in R^\sigma_1(f)_{(*-1)\beta+\beta_1^\sigma}$
such that $AB\in R^\sigma_1(f)_{(d-3)\beta+2\beta_1^\sigma}$,   
$$\gamma^i_{A}\cup\gamma^i_{B}=\frac{\mult(\sigma_k+\sigma_{k+1})
\gamma_{\mu^{-1}(ABG^\sigma(f))}}{\mult(\sigma_k)\mult(\sigma_{k+1})}\text {\quad in\quad} H^{d-1}(X,\wedge^{d-1}{\cal T}_X),$$
where the map $\mu$ and  $G^\sigma(f)\in S_{3\beta-2\beta_1^\sigma}$ are defined in 
(\ref{e:isomorm})  and Proposition~\ref{p:comp},

{\rm(v)}  for $\rho_i,\rho_j\notin\Sigma_X$ which span 
a $2$-dimensional cone  $\sigma'\in\Sigma$ contained in $\sigma\in\Sigma_X(2)$ and $A,B$ as in {\rm(iv)},
$$\gamma^i_{A}\cup\gamma^j_{B}=-\frac{\gamma_{\mu^{-1}(ABG^\sigma(f))}}{\mult(\sigma')}\text {\quad in\quad}
H^{d-1}(X,\wedge^{d-1}{\cal T}_X),$$

{\rm(vi)}
for  $\rho_i=\rho_{l_k}$  as in (\ref{e:pic}) and $A,B,C\in R^\sigma_1(f)_{(*-1)\beta+\beta_1^\sigma}$
such that $ABC\in R^\sigma_1(f)_{(d-4)\beta+3\beta_1^\sigma}$,
 $$\gamma^i_{A}\cup\gamma^i_{B}\cup\gamma^i_{C}=\frac{(\mult(\sigma_k)-\mult(\sigma_{k+1}))\mult(\sigma_k+\sigma_{k+1})
\gamma_{\mu^{-1}(ABCH^\sigma(f)G^\sigma(f))}}{(\mult(\sigma_k)\mult(\sigma_{k+1}))^2},$$
where  $H^\sigma(f)\in S_{\beta-\beta_1^\sigma}$ is defined in Remark~\ref{r:hs}, 

{\rm(vii)} for $\rho_i=\rho_{l_k}$ as in (\ref{e:pic})  and $A,B,C$ as in {\rm(vi)},
$$\gamma^i_{A}\cup\gamma^i_{B}\cup\gamma^{l_{k\pm1}}_{C}=\mp\frac{\gamma_{\mu^{-1}(ABCH^\sigma(f)G^\sigma(f))}}
{\mult(\sigma_{k,k\pm1})^2},$$
where $\sigma_{k,k\pm1}$ denotes the cone spanned by $\rho_i$ and $\rho_{l_{k\pm1}}$.
\end{thm}

\begin{rem} If  the multiplicities of the 2-dimensional cones of the fan $\Sigma$, lying inside
a cone of $\Sigma_X(2)$, are equal to 1,
then $\gamma^i_{A}\cup\gamma^i_{B}\cup\gamma^i_{C}=0$ in part (vi) of the above theorem.
 In particular, this holds for the minimal Calabi-Yau hypersurfaces in \cite{b1}.
\end{rem}

\end{document}